\documentclass[]{scrartcl}
\usepackage{amsfonts}
\usepackage{amsmath}
\usepackage{amssymb}
\usepackage{amsthm}
\usepackage{titlesec}
\usepackage{tikz}
\usepackage{mathtools}
\usepackage{subcaption}
\titleformat{\subsection}[runin]{\normalfont\bfseries}{\thesubsection.}{3pt}{}
\titleformat{\subsubsection}[runin]{\normalfont\bfseries}{\thesubsubsection.}{3pt}{}
\usepackage[all]{xy}
\usepackage{float}
\usepackage{color}

\DeclarePairedDelimiter\floor{\lfloor}{\rfloor}

\newcommand{\sigmahat}{\widehat{\sigma}}
\newcommand{\Pihat}{\widehat{\Pi}}
\newcommand{\semisuspension}{\widetilde{\Sigma}}
\newcommand{\toskel}{\phi^{\Pihat}}
\newcommand{\fromskel}{\phi^{\Pi}}
\newcommand{\flagtop}{\Upsilon_{[\sigma, \widehat{1}]}}
\newcommand{\cdtop}{\Phi_{[\sigma, \widehat{1}]}}

\newcommand{\checkcd}{\check{\varPhi}}
\newcommand{\lowint}[1]{[\widehat{0}, #1]}
\newcommand{\localcd}[1]{\ell^\varPhi_{#1}}
\newcommand{\localab}[1]{\ell^\varPsi_{#1}}
\newcommand{\localflag}[1]{\ell^\varUpsilon_{#1}}   
\newcommand{\Pyr}{\text{Pyr}}
\newcommand{\stard}[1]{\text{star}_{\partial \Delta_d}(\Delta_{#1})}
\newcommand{\upint}[1]{[#1, \widehat{1}]}

\definecolor{cof}{RGB}{219,144,71}
\definecolor{pur}{RGB}{186,146,162}
\definecolor{greeo}{RGB}{91,173,69}
\definecolor{greet}{RGB}{52,111,72}

\theoremstyle{definition}

\newtheorem{definition}{Definition}[section]
\newtheorem{theorem}[definition]{Theorem}
\newtheorem{lemma}[definition]{Lemma}
\newtheorem{corollary}[definition]{Corollary}
\newtheorem{proposition}[definition]{Proposition}
\newtheorem{example}[definition]{Example}

\begin{document}
\pagenumbering{roman}

\begin{titlepage}
\centering

\vspace*{3cm}

{\huge  \textbf{Subdividing the $cd$-index of Eulerian Posets}}

\vspace{1cm}
by
\vspace{1cm}

Patrick Dornian

\vspace{3cm}

A thesis

presented to the University Of Waterloo

in fulfillment of the requirement for the degree of

Master of Mathematics

in

Combinatorics and Optimization

\vspace{3cm}
Waterloo, Ontario, Canada, 2016

\vspace{1cm}

\textcopyright Patrick Dornian 2016

\end{titlepage}

\begin{center}
{\large \textbf{Author's Declaration}}
\end{center}
I hereby declare that I am the sole author of this thesis. This is a true copy of the thesis, including any required final revisions, as accepted by my examiners.\\

I understand that my thesis may be made electronically available to the public.
\pagebreak

\begin{center}
{\large \textbf{Abstract}}
\end{center}

This thesis aims to give the reader an introduction and overview of the $cd$-index of a poset, as well as establish some new results. We give a combinatorial proof of Ehrenborg and Karu's $cd$-index subdivision decomposition for Gorenstein* complexes and extend it to a wider class of subdivisions. In doing so, we define a \emph{local cd-index} that behaves analogously to the well studied local $h$-vector. We examine known $cd$-index and $h$-vector bounds, and then use the local $cd$-index to bound a particular class of polytopes with the $cd$-index of a stacked polytope.  We conclude by investigating the $h$-vector and local $h$-vector of posets in full generality, and use an algebra morphism developed by Bayer and Ehrenborg to demonstrate the structural connection between the $cd$-index subdivision decomposition and the local $h$-vector subdivision decomposition.

\pagebreak

\begin{center}
{\large \textbf{Acknowledgements}}
\end{center}

I would first like to thank my supervisor Eric Katz for his mathematical guidance, unwavering support, and poor jokes. \\

Next, I would like to thank my readers Kevin Purbhoo and David Wagner for sacrificing their valuable time in an attempt to parse this tome.\\

Special thanks to my office mates and friends Cameron, Garnet and Anirudh for creating the best work environment I've ever been a part of (and likely the best I'll ever find).\\

My parents have spent the last two years having to deal with me ranting about incomprehensible math as a substitute for small talk. Thanks for putting up with me. You supported me through my entire education and raised me with the capacity to see this through. The early math education that you gave me was the start of this long road. The love you provided was invaluable. The money and care packages were pretty good too.\\

My Waterloo family consists of my loving girlfriend Lilly Zheng and our cats Hodor and Versace. All three provided essential inspiration.\\

Finally, thanks to all my friends for keeping me sane. In no particular order, the following answered my open facebook invitation for people to acknowledge. I'm almost certainly missing a few names, so I preemptively apologize (Sorry).\\ 

{\footnotesize

\begin{minipage}{.5\textwidth}
Carlo Arcovio

Ryan Hancock

Aidan Waite

Jonathan Dornian

Jaqueline Li

Alexander Morash

Kirsten Hattori

Melissa Angyalfi

Bashar Jabbour

Reagan Elly

Jamie Waugh

Andrew Cottle

Kat Dornian

Bill Irons

Kevin Chapman

Jason LeGrow

Steven Sun

Stephanie Raphael

Erin Perri

Darren Stalker

Prashanth Madhi

Evan Ferguson

Pavel Shuldiner

\end{minipage}%
\begin{minipage}{.5\textwidth}
Tavian Barnes

Catherine Maggiori 

Adam Gomes

Tevin Straub

Sakib Imtiaz

Andrijana Nesic

Shawn Puthukkeril

Mike Daw

Zach Neshevich

Dain Galts

Lily Wang

Francis Williams

Ravi Goundalkar

Steven Ye

Charlie Payne

Nolan Shaw

Amy Dornian

Greg Gregory Greggington Lewis

Benjamin Graf

Ella Weatherilt

Erika Angyalfi

Rachel Dornian

Nick Pulos

\end{minipage}

}

\pagebreak

\begin{center}
{\large \textbf{Dedication}}

\vspace{8cm}

In loving memory of Paul Walker and Rob Ford

\vspace{8cm}
\emph{it's been a long day}\\
\emph{without you my friend}\\
\emph{and i'll tell you all about it}\\
\emph{when i see you again......}\\
\end{center}
\pagebreak

\tableofcontents
\pagebreak

\pagenumbering{arabic}
\section{Introduction}

The $f$-vector and $h$-vector of a simplicial complex are classical equivalent tools for encoding its face numbers. Analogously, the flag $f$-vector and the flag $h$-vector were developed by Bayer and Billera \cite{bayerds} in 1985 to enumerate flags in posets. If the poset is Eulerian, we may efficiently encode these vectors into a $cd$-index, a non-commutative generating function developed by Bayer, Fine and Klapper. \cite{bayercd}. The motivating example for an Eulerian poset is typically the face lattice of a convex polytope.

Though the $cd$-index is efficient, its behavior is difficult to characterize. Given a simple geometric transformation on a polytope, the corresponding action on the $cd$-index is often hard to express. Its coefficients are non-negative for large classes of posets, but finding a natural combinatorial interpretation of them is an unresolved problem. Stanley uses a shelling argument to show that the $cd$-index of a regular shellable $CW$-sphere is non-negative \cite{stanleycd}. We paraphrase this result and demonstrate some applications of it as it motivates many of the later techniques we use. 

Using commutative algebra, Karu established non-negativity for complete and quasi-convex fans, as well as for Gorenstein* posets \cite{karufans}. Ehrenborg and Karu extended this result to near-Gorenstein* posets soon after \cite{borgkaru}. In the same article, they develop a decomposition theorem for the $cd$-index of subdivisions of Gorenstein* posets by using the theory of sheaves on fans and posets. We relax their assumptions to a new, slightly broader model of \emph{strongly Eulerian subdivisions} to produce a natural combinatorial model for the decomposition, in the process defining a \emph{local cd-index}. This local $cd$-index behaves analogous to the local $h$-vector introduced by Stanley \cite{stanleylocal}.

We then investigate known bounds on the $cd$-index and the $h$-vector. Using the local $cd$-index, we demonstrate a new bound: that the $cd$-index of a stacked polytope is an upper bound for the $cd$-index of simplicial spheres that may be triangulated by a shellable complex. Finally, we then examine the connections between the $cd$-index subdivision decomposition and the general local $h$-vector decomposition of a subdivision by using techniques developed by Bayer and Ehrenborg in  \cite{bayerborg}.  

\clearpage
\section{The ab-index and cd-index of a poset}

We will first state definitions and notation on partially ordered sets.

\begin{definition}
A \emph{partially ordered set} or \emph{poset} is a set $P$ with a binary relation $\leq$ that satisifies the following three axioms for all elements $s,t,u \in P$:

\begin{enumerate}
\item Reflexivity: $t \leq t$.
\item Antisymmetry: If $s \leq t$ and $t \leq s$ then $s = t$.
\item Transitivity: If $s \leq t$ and $t \leq u$ then $s \leq u$.
\end{enumerate} 

\end{definition}
Given a poset $P$, we call a subset $C \subseteq P$ a \emph{chain} of $P$ if the elements of $C$ are totally ordered. That is,

$$C = x_0 < x_1 < \dots < x_k$$

for some $x_i \in P$. Alternatively, these may also be referred to as \emph{flags}. A chain is \emph{maximal} if it is not contained in any other chain of $P$. 

\begin{definition}
Given a poset $P$ with $s \leq t$, we say that $t$ \emph{covers} $s$ and write $s \lessdot t$ if there exists no $u \in P$ such that $s < u < t$.
\end{definition}

\begin{definition}
Given a poset $P$ and $s,t \in P$, we define the \emph{interval} $[s,t]$ by

$$ [s,t] = \{x \, | \, s \leq x \leq t, x \in P \} \, .$$
\end{definition}

Note that $[s,t]$ is also a well defined poset with order relation inherited from $P$.\\

We use $\widehat{0}$ and $\widehat{1}$ to denote elements (if they exist) $\widehat{0}, \widehat{1} \in P$ such that for all $s \in P$, $\widehat{0} \leq s$ and $s \leq \widehat{1}$. We say a chain is \emph{degenerate} if it contains $\widehat{0}$ or $\widehat{1}$ (because any chain may be trivially extended by adding a maximal or minimal element).

\begin{definition}
Suppose every maximal chain in $P$ has length $n$. We then say that $P$ is \emph{graded} with \emph{rank} $n$.
\end{definition}

\begin{definition}
If a poset $P$ is graded, there exists a unique rank function $\rho_P: \, P \rightarrow \{0,1, \dots , n\}$ such that:
\begin{enumerate}
\item $\rho_P(s) = 0$ if $s$ is a minimal element of $P$.
\item If $s \lessdot t$, then $\rho_P(t) = \rho_P(s) + 1$.
\end{enumerate}
\end{definition}

When context is clear, we will omit the subscript. If $\rho(s) = k$, we say that $s$ has rank $k$. Given a nonempty interval $[s,t]$ we define the \emph{length} of $[s,t]$ to be $\rho(s,t) = \rho(t) - \rho(s)$. \\

Note that if $P$ is graded of rank $n$ and has a $\widehat{0}$ and a $\widehat{1}$ we necessarily have that $\widehat{0}$ is of rank $0$ and $\widehat{1}$ is of rank $n$. We call elements of rank $1$ \emph{atoms} of $P$ and elements of rank $n-1$ to be \emph{coatoms}.

\begin{definition}

We say that a finite graded poset $P$ with $\widehat{0}$ and $\widehat{1}$ is \emph{Eulerian} if every nonempty interval $[s,t]$ contains an equal number of elements of odd rank and even rank.
\end{definition}

Eulerian posets generalize the properties of the Euler characteristic of convex polytopes. We say that a poset is \emph{lower Eulerian} if all its intervals are Eulerian and it has a $\widehat{0}$. Lower Eulerian posets are motivated by polyhedral complexes. Note that an Eulerian poset is trivially also lower Eulerian. \\

For the remainder of this paper, we will let $P$ denote a finite graded poset of rank $n$ with $\widehat{0}$ and $\widehat{1}$ unless otherwise specified. Note that every maximum length chain in P  will have the form $\hat{0} = x_0 < x_1 < ... < x_n = \hat{1}$, with $\rho(x_i) = i$. Since any maximum length chain must begin with $\hat{0}$ and end with $\hat{1}$, their presence is often taken as a given in our calculations. Unless otherwise specified, all chains are non-degenerate. We use $[n]$ to denote the set $[n] := \{1, 2, ..., n\}$.
 
\begin{definition}

For $S \subseteq \{0, 1, \dots, n \}$ we let $P_S$ denote the \emph{S-rank selected subposet} \\ of P,
 
 \[P_S = \{x \in P : \rho(x) \in S\}\]
 
 \end{definition} 
 
\begin{definition}
The \emph{flag f-vector} of $P$ is defined by the function  $\alpha_P: 2^{[n]} \rightarrow \mathbb{Z}$ \\
 where $\alpha_P (S) =$ the number of maximal chains of $P_S$. 
 \end{definition}
 
 Where there is no confusion about $P$, we will often simply denote it by $\alpha$. We may think of $\alpha$ as a vector indexed by all possible subsets of $[n]$.
 
\begin{definition} The \emph{flag h-vector} of $P$ is defined by the function $\beta_P: 2^{[n]} \rightarrow \mathbb{Z}$ where
 \[\beta_P(S) = \sum_{T\subseteq S} (-1)^{\#(S \setminus T)} \alpha_P(T)\]
 \end{definition}
 
 Equivalently, by inclusion-exclusion we have the following more elegant expression.
 
 \[\alpha(S) = \sum_{T \subseteq S} \beta(T) \]
 
 Our next goal is to encode the flag $f$ and $h$-vector into  generating functions. To correspond to a given rank selection $S \subseteq [n]$, we define the noncommutative characteristic monomial $u_S = u_1 u_2 ... u_n$ in the variables $a$ and $b$ by
 
 \[u_i = \begin{cases}
 a& \text{if $i \notin S$} \\
 b& \text{if $i \in S$}
 \end{cases}\]
 
 For example, if $n = 5$ and $S = \{1, 2, 5\}$ we have $u_S = bbaab$. This lets us define the following generating function.
 
\begin{definition} The \emph{ab-index} for a graded poset $P$ of rank $n + 1$ with $\hat{0}$ and $\hat{1}$ is the generating function defined by
 
 \[\Psi_P (a,b) = \sum_{S \subseteq [n]} \beta_P (S) u_S \, . \]
 \end{definition}
 
\begin{definition} The \emph{flag polynomial} for a graded poset $P$ of rank $n + 1$ with $\hat{0}$ and $\hat{1}$ is defined by
 
  \[\Upsilon_P (a, b) = \sum_{S \subseteq [n]} \alpha_P (S) u_S \]
\end{definition}
  
It is easily seen that the flag polynomial and the $ab$-index are equivalent by a linear change of variables. We can transform between then by the following relationship \cite{stanleycd}.

\begin{align*}
\Upsilon_P(a,b) &= \Psi_P(a + b, b), \\
 \Psi_P(a, b) &=\Upsilon_P(a - b,b).
\end{align*}

When convenient, we may omit the $(a,b)$ argument or use the poset $P$ itself as the argument of the function if context is clear. The following alternative notation is also useful.\\

Given a poset $P$ of rank $n+1$, consider some non-degenerate chain $C = x_1 < x_2 < \dots < x_k$. We define the following functions that take chains from $P$ to $ab$-words. First, we set

$$
\alpha^C = u_1 u_2 \dots u_n
$$
where 

 \[u_i = \begin{cases}
 a& \text{if $\nexists$ some $j$ such that $\rho(x_j) = i$} \\
 b& \text{if $\exists$ some $j$ such that $\rho(x_j) = i$}
 \end{cases}\]
 
 and
 
 $$
 \beta^C = u_1 u_2 \dots u_n
 $$
 
 where
 
  \[u_i = \begin{cases}
  a - b& \text{if $\nexists$ some $j$ such that $\rho(x_j) = i$} \\
  b& \text{if $\exists$ some $j$ such that $\rho(x_j) = i$}
  \end{cases} \, .\]
  
  In other words, $\alpha^C$ encodes the ranks of chain $C$ as an $ab$-word where ranks that are present are weighted with a $b$ while missing ranks are weighted with an $a$.
  
  $\beta^C$ encodes $C$ so that ranks present are weighted with $b$, and ranks missing are weighted with $(a-b)$. 
  
  \begin{proposition}

   Given a graded poset $P$, the following two equations hold (and may serve as alternate definitions of the $ab$-index and flag polynomial if desired).
  
  $$\Upsilon_P = \sum_{\substack{C \\ C \text{ a chain of } P}} \alpha^C $$

  $$\Psi_P = \sum_{\substack{C \\ C \text{ a chain of } P}} \beta^C \, .$$
  
  \end{proposition}
  
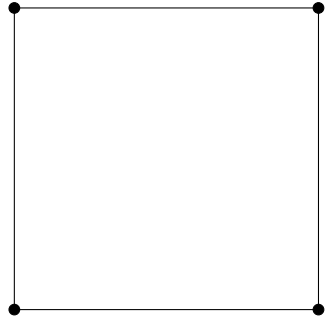
\begin{figure}[h!]
\centering 
\caption{The $ab$-index and flag enumerator of the face poset of a square.}
\begin{tikzpicture}[scale=2]
\draw (0,0)--(2,0)--(2,2)--(0,2)--(0,0);
\draw [fill] (0,0) circle [radius = 1pt];
\draw [fill] (2,0) circle [radius = 1pt];
\draw [fill] (2,2) circle [radius = 1pt];
\draw [fill] (0,2) circle [radius = 1pt];
\end{tikzpicture}
\begin{align*}
\Upsilon_\Gamma &= aa + 4ba + 4ab + 8bb\\
\Psi_\Gamma &= aa + 3ba + 3ab + bb
\end{align*}
\end{figure}

Note that in the following example, our $ab$-index is symmetric. That is, it is invariant under the involution that switches $b$'s with $a$'s and visa versa. This is not a coincidence. Such a statement holds for all Eulerian posets \cite{stanleyecvol1}. 
\begin{lemma}
Let $P$ be an Eulerian poset of rank $n + 1$. Let $S \subseteq [n]$ and denote $\overline{S} = [n] \setminus S$. We then have $\beta_P(S) = \beta_P(\overline{S})$.
\end{lemma}

The symmetry of the $ab$-index is a generalization of the symmetry of the $h$-vector of simplicial polytopes and spheres. That is, it is a generalization of the Dehn-Sommerville equations to flags of a poset \cite{bayerds}. In particular, this allows a particular substitution for the $ab$-index of an Eulerian poset.

\begin{theorem} \emph{(Bayer, Fine, Klapper)} \cite{bayercd} Let $P$ be an Eulerian poset. Then there exists a non-commutative polynomial in the variables $c$ and $d$ denoted $\Phi_P(c,d)$ that satisfies $\Phi_P(a + b, ab + ba) = \Psi_P(a,b)$. 
\end{theorem}

We call this polynomial the $cd$-index. For clarity in particular arguments we may refer to the $cd$-index of $P$ as any of $\Phi_P (c,d)$, $\Phi_P$, or $\Phi(P)$ depending on context.

\begin{example}
In the above figure of a square with face complex $\Gamma$, we have that $\Phi_\Gamma = c^2 + 2d$.
\end{example}

The $cd$-index lacks the natural symmetry of the $ab$-index, but represents a vastly more concise encoding of the flag data. The linear span of flag $f$ or $h$-vectors of graded rank $n+1$ posets has dimension at most $2^n$, corresponding to all possible length $n$ $ab$-words. However, the linear span of flag $f$ or $h$-vectors of Eulerian posets of rank $n+1$ has dimension at most $\mathcal{F}_{n+1}$ where $\mathcal{F}_{n+1}$ is the $(n+1)^{\text{st}}$ Fibonacci number. This corresponds to the number of $cd$ words of degree $n$ where $c$ has degree $1$ and $d$ has degree $2$. Stanley demonstrates that this bound is tight \cite{stanleyecvol1}.\\

In many cases, the $cd$-index of a poset has non-negative coefficients. Determining classes of posets for which the $cd$-index is non-negative is an active area of research. In the next chapter, we discuss Stanley's proof of its non-negative for a certain class of shellable spheres (including polytopes). In Chapter Four, we discuss Karu's extension of this result to Gorenstein* posets.

\clearpage
\section{Shelling the $cd$-index}

The following section develops the theory of shelling and summarizes the relevant results from Chapter Two of \emph{Flag f-vectors and the cd-index} by Stanley \cite{stanleycd}. Many of Stanley's techniques motivate our later results so they are included here for reference. \\

In layman's terms, a shelling is a linear order the facets of a complex such that when `drawn' sequentially, each facet intersects the previous facets in a continues segment (with one exception in the one-dimensional case, but this definition suffices for most geometric intuition). It is a useful tool for proving results inductively is the main technique behind early results on the non-negativity of the $cd$-index.\\

All polytopes mentioned past this point are convex. No non-convex polytopes are discussed in this thesis. In lieu of writing 'convex' repeatedly, we will refer to 'polytopes' with the implicit understanding that we are referring to the convex case.

\begin{definition}
 A \emph{polytopal complex} $\mathcal{C} \subseteq \mathbb{R}^d$ is a finite set of polytopes such that

\begin{enumerate}
\item The empty face is present. That is, $\emptyset \in \mathcal{C}$.
\item If $P \in \mathcal{C}$, then all proper faces of $P$ are also in $C$. 
\item For any $P, Q \in \mathcal{C}$, if $P \cap Q = F$ where $F \neq \emptyset$, then $F$ is a face of both $P$ and $Q$.  
\end{enumerate}
\end{definition}

We may analogously define simplicial complexes by replacing all references above to polytopes with simplices. It is easy to see that any complex has a well defined partial order of its faces under inclusion. We will typically be considering the face complex of a polytope $P$: that is, the set of all its faces under inclusion. Another common construction is the boundary complex $\partial P$, the set of all faces of $P$ excluding the maximal face of $P$ itself. When context is clear we will use $P$ to denote both the polytope itself, as well as the underlying poset of its face complex. When context is less clear, we will use the notation $F(P)$ for the partial order. \\

Note that the face poset of a complex is graded with a well defined rank function induced by dimension. Recall that the empty face $\emptyset$, the minimal face of the complex is typically noted as having dimension $-1$. By definition of minimal element, we must have $\rho(\emptyset) = 0$ and it follows from an easy inductive argument that $\rho(F) = \dim F + 1$ for any $F \in P$. That is, vertices are of dimension 0, but rank 1, edges are of dimension 1 but rank 2 and so on. This shift in dimension can lead to endless arithmetic errors if one is careless, so it is important to keep in mind when generalizing theorems on complexes to posets. \\

 The dimension of a complex $\mathcal{C}$ is the maximum dimension of a face in $\mathcal{C}$. We say that a complex is \emph{pure} if all maximal faces under inclusion have the same dimension.\\

There are many slight variations on the definition of shelling depending on the authors motives. The following is presented by Ziegler \cite{ziegler} and suffices for most combinatorial purposes.

\begin{definition}
Let $\Gamma$ be a pure $d$-dimensional polytopal complex. We say $\Gamma$ is \emph{shellable} if either $\Gamma$ is $0$-dimensional, or there exists a linear order (called a \emph{shelling}) of its facets $F_1, F_2, \dots F_k$ that satisfies the following recursive definition: 

\begin{enumerate}
\item The boundary complex $\partial F_1$ has a shelling.
\item For $1 < j \leq k$, we have

$$F_j \cap \left( \bigcup_{i=1}^{j-1} F_i \right) = G_1 \cup G_2 \dots G_r$$

for some shelling $G_1, G_2, \dots G_r, \dots G_s$ of $\partial F_j$ with $1 \leq r \leq s$. In other words, the intersection between $F_j$ and the union of the previous facets in the shelling is the beginning segment of some shelling of $\partial F_j$.
\end{enumerate}
\end{definition}

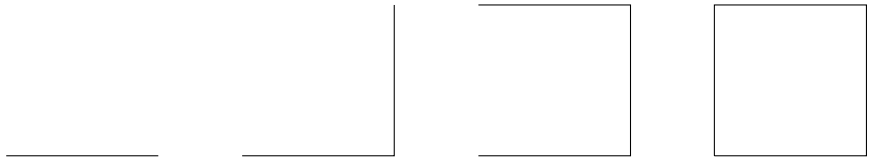
\begin{figure}[h!]
\centering
\begin{tikzpicture}
\draw (0,0) -- (2,0);
\end{tikzpicture}
\, \, \, \, \,
\begin{tikzpicture}
\draw (0,0) -- (2,0)--(2,2);
\end{tikzpicture}
\, \, \, \, \,
\begin{tikzpicture}
\draw (0,0) -- (2,0)--(2,2)--(0,2);
\end{tikzpicture}
\, \, \, \, \, 
\begin{tikzpicture}
\draw (0,0)--(2,0)--(2,2)--(0,2)--(0,0);
\end{tikzpicture}

\caption{We may shell a polygon by `walking' around its edges.}
\end{figure}

Shellings have been widely studied and are a useful inductive tool on complexes and polytopes. In particular, any result that can be proven for a shellable complex holds for all polytopes because of the following result of Bruggesser and Mani \cite{brugmani}.

\begin{theorem}
 Let $P$ be a polytope. Then its boundary complex $\partial P$ is shellable. \\

In particular, they show that every polytope may be shelled using a construction called a \emph{line shelling} \cite{brugmani}.
\end{theorem}

\begin{figure}
\centering

\begin{minipage}{0.5\textwidth}
\centering
\begin{tikzpicture}[thick,scale=.4]
\coordinate (p1) at (0,0);
\coordinate (p2) at (7,0);
\coordinate (p3) at (0,6);
\coordinate (p4) at (7,6);
\coordinate (p5) at (4,2);
\coordinate (p6) at (11,2);
\coordinate (p7) at (4,8);
\coordinate (p8) at (11,8);

\begin{scope}[thick,dashed,]
\draw (p1) -- (p5) ;
\draw (p5) -- (p6);
\end{scope}

\draw (p1) -- (p2) -- (p6);

\end{tikzpicture} 
\end{minipage}%
\begin{minipage}{0.5\textwidth}
\centering
\begin{tikzpicture}[thick,scale=.4]
\coordinate (p1) at (0,0);
\coordinate (p2) at (7,0);
\coordinate (p3) at (0,6);
\coordinate (p4) at (7,6);
\coordinate (p5) at (4,2);
\coordinate (p6) at (11,2);
\coordinate (p7) at (4,8);
\coordinate (p8) at (11,8);

\begin{scope}[thick,dashed,]
\draw (p1) -- (p5) -- (p7);
\draw (p5) -- (p6);
\end{scope}

\draw (p1) -- (p2) -- (p6);
\draw (p7) -- (p8);
\draw (p6) -- (p8);

\end{tikzpicture}
\end{minipage}
\vspace{.5cm}

\begin{minipage}{0.5\textwidth}
\centering
\begin{tikzpicture}[thick,scale=.4]
\coordinate (p1) at (0,0);
\coordinate (p2) at (7,0);
\coordinate (p3) at (0,6);
\coordinate (p4) at (7,6);
\coordinate (p5) at (4,2);
\coordinate (p6) at (11,2);
\coordinate (p7) at (4,8);
\coordinate (p8) at (11,8);

\begin{scope}[thick,dashed,]
\draw (p1) -- (p5) -- (p7);
\draw (p5) -- (p6);
\end{scope}

\draw (p1) -- (p2) -- (p6);
\draw (p3) -- (p7) -- (p8);
\draw (p1) -- (p3);

\draw (p6) -- (p8);

\end{tikzpicture}
\end{minipage}%
\begin{minipage}{0.5\textwidth}
\centering
\begin{tikzpicture}[thick,scale=.4]
\coordinate (p1) at (0,0);
\coordinate (p2) at (7,0);
\coordinate (p3) at (0,6);
\coordinate (p4) at (7,6);
\coordinate (p5) at (4,2);
\coordinate (p6) at (11,2);
\coordinate (p7) at (4,8);
\coordinate (p8) at (11,8);

\begin{scope}[thick,dashed,,opacity=0.6]
\draw (p1) -- (p5) -- (p7);
\draw (p5) -- (p6);
\end{scope}

\draw (p1) -- (p2) -- (p6);
\draw (p4) -- (p3) -- (p7) -- (p8);
\draw (p1) -- (p3);
\draw (p2) -- (p4);
\draw (p6) -- (p8);

\draw[fill=cof,opacity=0.4] (p1) -- (p3) -- (p4) -- (p2)-- (p1);

\end{tikzpicture}
\end{minipage}

\vspace{.5cm}

\begin{minipage}{0.5\textwidth}
\centering
\begin{tikzpicture}[thick,scale=.4]
\coordinate (p1) at (0,0);
\coordinate (p2) at (7,0);
\coordinate (p3) at (0,6);
\coordinate (p4) at (7,6);
\coordinate (p5) at (4,2);
\coordinate (p6) at (11,2);
\coordinate (p7) at (4,8);
\coordinate (p8) at (11,8);

\begin{scope}[thick,dashed,,opacity=0.6]
\draw (p1) -- (p5) -- (p7);
\draw (p5) -- (p6);
\end{scope}

\draw (p1) -- (p2) -- (p6);
\draw (p3) -- (p4) -- (p8) -- (p7) --cycle;
\draw (p1) -- (p3);
\draw (p2) -- (p4);
\draw (p6) -- (p8);

\draw[fill=cof,opacity=0.4] (p1) -- (p3) -- (p4) -- (p2)-- (p1);
\draw[fill=greeo,opacity=0.4] (p2) -- (p4) -- (p8) -- (p6)-- (p2);

\end{tikzpicture}
\end{minipage}%
\begin{minipage}{0.5\textwidth}
\centering
\begin{tikzpicture}[thick,scale=.4]
\coordinate (p1) at (0,0);
\coordinate (p2) at (7,0);
\coordinate (p3) at (0,6);
\coordinate (p4) at (7,6);
\coordinate (p5) at (4,2);
\coordinate (p6) at (11,2);
\coordinate (p7) at (4,8);
\coordinate (p8) at (11,8);

\begin{scope}[thick,dashed,,opacity=0.6]
\draw (p1) -- (p5) -- (p7);
\draw (p5) -- (p6);
\end{scope}

\draw (p1) -- (p2) -- (p6);
\draw (p3) -- (p4) -- (p8) -- (p7) --cycle;
\draw (p1) -- (p3);
\draw (p2) -- (p4);
\draw (p6) -- (p8);

\draw[fill=cof,opacity=0.4] (p1) -- (p3) -- (p4) -- (p2)-- (p1);
\draw[fill=pur,opacity=0.4] (p3) -- (p4) -- (p8) -- (p7)-- (p3);
\draw[fill=greeo,opacity=0.4] (p2) -- (p4) -- (p8) -- (p6)-- (p2);

\end{tikzpicture}
\end{minipage}
\caption{A shelling of a rectangular prism.}
\end{figure}
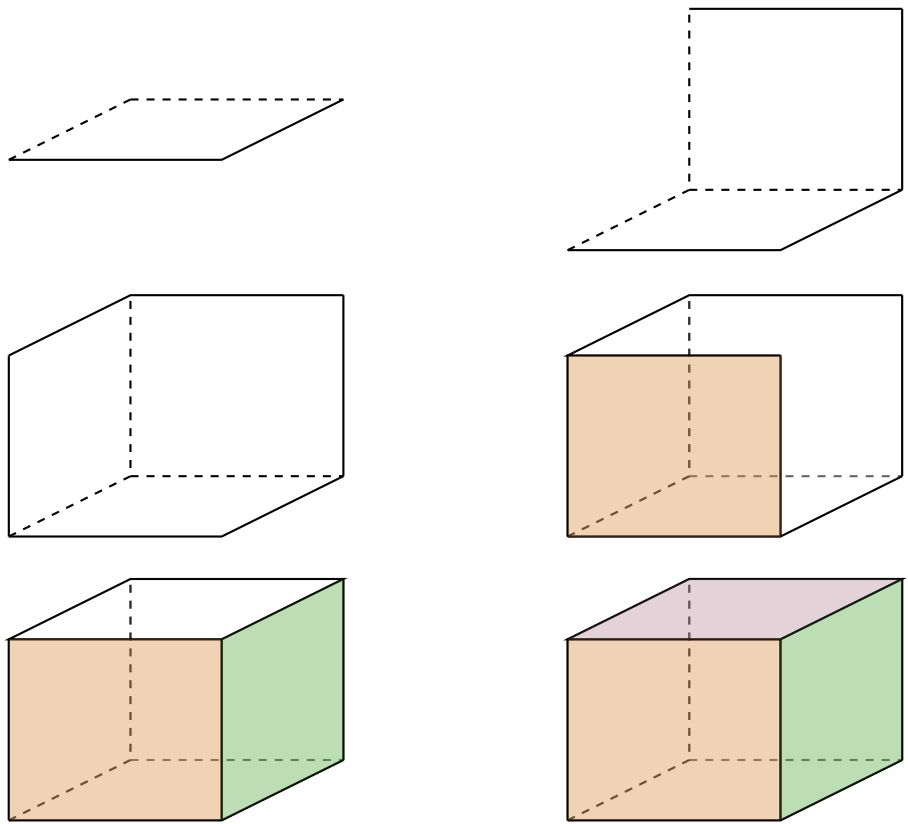

In proving results on polytopes, it is often of use to generalize to the looser topological framework of $CW$-spheres. The following definitions coincide with Munkres \cite{munkres}.

\begin{definition}
A space is a \emph{cell} of dimension $n$ if it is homeomorphic to the closed ball $B^n$. It is an \emph{open cell} if it is homeomorphic to Int$B^n$.
\end{definition}

\begin{definition}
A \emph{CW complex} is a topological space $X$ partitioned into open cells $e_\alpha$ such that:

\begin{enumerate}
\item $X$ is Hausdorff.
\item For each open $n$-cell $e_\alpha$, there exists a continuous map $f_\alpha: B^n \rightarrow X$ that induces a homeomorphism between Int$B^n$ and $e_\alpha$ and carries $\partial B^n$ into a finite union of open cells of dimension less than $n$.

\item $X$ is equipped with the weak topology: A set $A \subseteq X$ is open (or closed) if and only if $A \cap \bar{e}_\alpha$ is open (or closed) in $\bar{e}_\alpha$ for all $\alpha$.
\end{enumerate}
\end{definition}

The final condition is redundant when dealing with complexes of only a finite number of elements (as we are in this paper) but is included for completeness. \\

If $f_\alpha$ is a homeomorphism for all $\alpha$, we say that the complex is a \emph{regular} $CW$-complex (that is, no identifications are made on the boundary of $e_\alpha$). \\

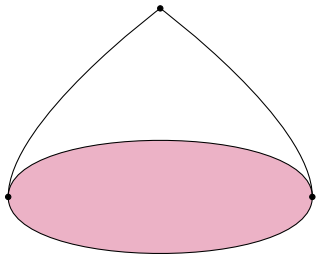
\begin{figure}
\centering
\begin{tikzpicture}
\filldraw [fill=purple!30!white, draw=black ] (-2,0) .. controls (-2,1) and (2,1) .. (2,0) .. controls (2,-1) and (-2,-1) .. (-2, 0) ;
\draw (-2,0).. controls (-2, 1) and (0,2.5)..(0,2.5)..controls (0,2.5) and (2,1)..(2,0);

\draw [fill] (2,0) circle [radius = 1pt];
\draw [fill] (-2,0) circle [radius = 1pt];
\draw [fill] (0, 2.5) circle [radius = 1pt];
\end{tikzpicture}

\caption{A regular $CW$-complex with one face, four edges and three vertices.}
\end{figure}

Regular $CW$-complexes are perhaps the most general framework for a topological space that still maintain combinatorial details. Any regular $CW$-complex $\Gamma$ of dimension $n$ also has a graded poset structure of rank $n+1$ under the order relation  that given cells $e_x , e_y \in \Gamma$, $e_x \leq e_y$ if and only if $e_x \subseteq \partial e_y$ topologically or $e_x = e_y$. We refer to this as the \emph{face poset} of $\Gamma$. We may interpret cells as both topological spaces and as the elements of a poset, and we hope that context ensures clarity. As in the case of polytopes, we refer to the maximal cells of a $CW$-complex as \emph{facets}. It follows that any polytopal complex $\Gamma$ is a regular $CW$-complex with identical poset structure by identifying all $d$-faces $\sigma \in \Gamma$ with an open $d-$cell $e_\sigma$ homeomorphic to its interior and maintaining all inclusion relations. From here on following we will consider faces of any polytopal/simplicial complex $\Gamma$ to be open cells of a $CW$-complex when topologically convenient.\\
  
  If $\Gamma$ is a regular $CW$-complex such that its underlying space $|\Gamma|$ is homeomorphic to a sphere, we will call $\Gamma$ a $CW$-\emph{sphere}. The face poset of $\Gamma$ is Eulerian \cite{stanleyeulerian}.  Additionally, all polytopes are $CW$-spheres. We may analogously define shelling on regular $CW$-complexes and it follows that any result that holds for a shellable $CW$-sphere must hold for a polytope.\\
  
  The following poset operations are introduced by Stanley in \cite{stanleycd} and will be of use to us in both this and further chapters. Though defined on arbitrary posets, they are geometrically motivated. 
  
\begin{definition}
 Given posets P and Q with $\widehat{0}$ and $\widehat{1}$, we define the \emph{join} $P * Q$ to be the poset on the set 
  \[P * Q := P \backslash \{\hat{1}\} \cup (Q \backslash \{\hat{0}\} \]
  
  with $x \leq y$ in $P * Q$ if either
  \begin{enumerate}
  \item $x \leq y$ in $P \backslash \{\widehat{1}\}$
  \item $x \leq y$ in $Q \backslash \{\widehat{0}\}$
  \item $x \in P \backslash \{\widehat{1}\}$ and $y \in Q \backslash \{\widehat{0}\}$
  \end{enumerate}
  \end{definition}
  
  We may think of this as placing $Q$ on top of $P$ in the ordering. We remove the maximal element of $P$, the minimal element of $Q$ and then set every element in $Q$ to be greater than those in $P$ with the pre-existing relations unchanged. It is easy to show that if $P$ and $Q$ are Eulerian, so is $P * Q$.
  
 \begin{figure}[h!]
 \begin{minipage}{.33\textwidth}
 \centering
 \begin{tikzpicture}
 \draw [fill] (0,2) circle [radius = 1pt];
 \draw [fill] (-.5,1) circle [radius = 1pt];
 \draw [fill] (.5,1) circle [radius = 1pt];
 \draw [fill] (0, 0) circle [radius = 1pt];
 
 \draw (0,2) -- (-.5,1)-- (0,0) -- (.5,1) -- (0,2);
 
 \end{tikzpicture}
\end{minipage}%
\begin{minipage}{.33\textwidth}
\centering
 \begin{tikzpicture}
  \draw [fill] (0,3) circle [radius = 1pt];
  \draw [fill] (-1,2) circle [radius = 1pt];
  \draw [fill] (0,2) circle [radius = 1pt];
  \draw [fill] (1, 2) circle [radius = 1pt];
  \draw [fill] (-1,1) circle [radius = 1pt];
  \draw [fill] (0,1) circle [radius = 1pt];
  \draw [fill] (1,1) circle [radius = 1pt];
  \draw [fill] (0, 0) circle [radius = 1pt];
  
  \draw (0,3) -- (-1,2);
  \draw (0,3) -- (0,2);
  \draw (0,3) -- (1,2);
  
  \draw (-1,2) -- (-1,1);
  \draw (-1,2) -- (0,1);
  
  \draw (0,2) -- (0,1);
  \draw (0,2) -- (1,1);
  
  \draw (1,2) -- (1,1);
  \draw (1,2)--(-1,1);
  
  \draw (0,0) -- (-1,1);
  \draw (0,0) -- (0,1);
  \draw (0,0) -- (1,1);
 \end{tikzpicture}
\end{minipage}%
\begin{minipage}{.33\textwidth} 
\centering
 \begin{tikzpicture}
   \draw [fill] (.5,3) circle [radius = 1pt];
   \draw [fill] (-.5,3) circle [radius = 1pt];
   \draw [fill] (0,3.5) circle [radius = 1pt];
   \draw [fill] (-1,2) circle [radius = 1pt];
   \draw [fill] (0,2) circle [radius = 1pt];
   \draw [fill] (1, 2) circle [radius = 1pt];
   \draw [fill] (-1,1) circle [radius = 1pt];
   \draw [fill] (0,1) circle [radius = 1pt];
   \draw [fill] (1,1) circle [radius = 1pt];
   \draw [fill] (0, 0) circle [radius = 1pt];

   \draw (-1,2) -- (-1,1);
   \draw (-1,2) -- (0,1);
   
   \draw (0,2) -- (0,1);
   \draw (0,2) -- (1,1);
   
   \draw (1,2) -- (1,1);
   \draw (1,2)--(-1,1);
   
   \draw (0,0) -- (-1,1);
   \draw (0,0) -- (0,1);
   \draw (0,0) -- (1,1);
   
   \draw (0.5,3) -- (-1, 2);
   \draw (0.5,3) -- (0, 2);
   \draw (0.5,3) -- (1, 2);
   
   \draw (-0.5,3) -- (-1, 2);
   \draw (-0.5,3) -- (0, 2);
   \draw (-0.5,3) -- (1, 2);
   
   \draw (0,3.5) -- (0.5, 3);
   \draw (0,3.5) -- (-0.5, 3);
  \end{tikzpicture}
  \end{minipage}
  \caption{Hasse diagrams of the boolean algebras $B_2, B_3$ and the join $B_3 * B_2$ respectively.} 
 \end{figure}
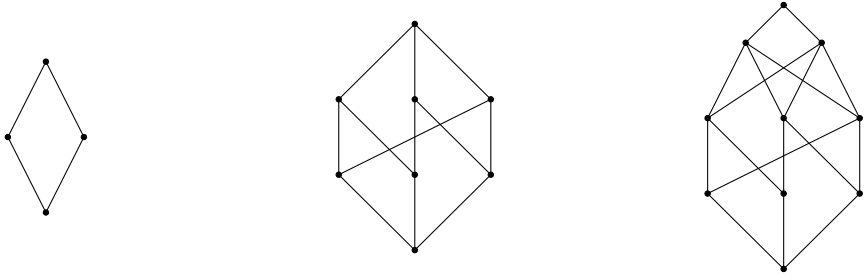
 
 \begin{lemma}

\cite{stanleycd} Let $P$ and $Q$ be Eulerian posets. Then the $cd$-index of $P*Q$ is 

$$\Phi_{P*Q} = \Phi_P \cdot \Phi_Q \, .$$

 \end{lemma}

\emph{Proof:} Chain counting lets us easily see that $\Upsilon_{P*Q}(a,b) = \Upsilon_P (a,b) \cdot \Upsilon_Q(a,b)$. Since we may obtain the $cd$-index by a series of variable substitutions from the flag polynomial, the result follows. \qed

 \begin{definition}

Let $B_2$ denote the boolean algebra of rank two. Let $P$ be a poset. Then we call $P * B_2$ the \emph{suspension} of $P$ and denote it as $\Sigma P$.
 \end{definition}
 
 The geometric intution for the suspension operator is as follows. Suppose $P$ has a geometric realization as the face poset of a regular $(n-1)$-dimension $CW$-sphere. Consider $\partial P := P \setminus \{\widehat{1}\}$. $\Sigma P$ contains exactly $\partial P$ (as its elements of rank $\leq n$), plus the two facets of rank $n+1$ introduced from $B_2$, along with the maximal element. Since these facets contain every face of $\partial P$ by definition of $\Sigma P$, we may view $\Sigma P$ as the complex induced by embedding $\partial P$ into the equator of an $n$-sphere with the two new facets corresponding to its northern and southern hemispheres. See Figure 6 for clarification.
  
 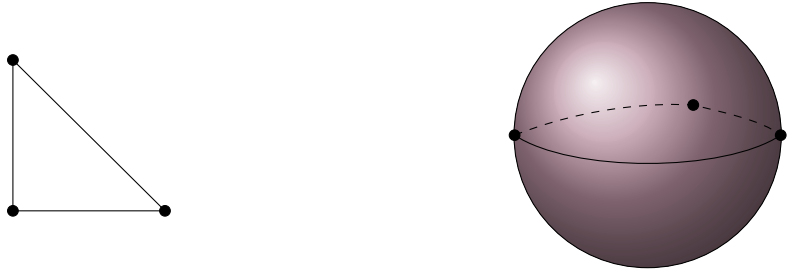
\begin{figure}
 \begin{minipage}{.5\textwidth}
 \centering
 \begin{tikzpicture}
 \draw [fill] (0,0) circle (2pt);
 \draw [fill] (2,0) circle (2pt);
 \draw [fill] (0,2) circle (2pt);
 \draw (0,0) -- (2,0) -- (0,2) -- cycle;
 \end{tikzpicture}
 \end{minipage}%
 \begin{minipage}{.5\textwidth}
 \centering
 \begin{tikzpicture}
 \shade[ball color=pur] (0,0) circle (50pt);
 \draw (0,0) circle (50pt);
 \draw [fill] (-1.75,0) circle (2pt);
 \draw [fill] (1.75,0) circle (2pt);
 \draw [fill] (.6, 0.4) circle (2pt);
 
 \draw (-1.75,0) .. controls (-1,-.5) and (1, -.5).. (1.75,0);
 \draw[dashed] (-1.75,0) ..controls (-.5, 0.5) and (.6, .4).. (.6, .4) .. controls (1,.3) and (1.5, .2).. (1.75,0);
 \end{tikzpicture}
 \end{minipage}
 \caption{$CW$-complexes associated to a triangle and its suspension, respectively.}
 \end{figure}

 \begin{corollary}

 Given an Eulerian poset $P$, we have that
 
$$\Phi_{\Sigma P} = \Phi_P \cdot c \, .$$
 \end{corollary}

\emph{Proof:} We calculate $\Phi_{B_2} = c$ and apply Lemma 3.7. \qed
 
 \begin{definition}$ \text{ }$
\begin{enumerate}
\item We say a poset $P$ is \emph{near-Eulerian} if $P = Q \setminus \{x\}$ where $Q$ is a rank $n+1$ Eulerian poset and $x$ is a coatom of $Q$ (that is, $x$ has rank $n$).

\item Given a near-Eulerian poset P, we may uniquely recover the associated Eulerian poset $Q$ by adding a new element $x$ such that $y \lessdot_Q x$ for any element $y \in P$ where $[y, \widehat{1}]_P$ contains exactly three elements. We call Q the \emph{semisuspension} of $P$ and denote it as $Q := \semisuspension P$.
\end{enumerate}
\end{definition}
The canonical example of a near-Eulerian poset is the face complex of a polytope with exactly one facet removed. We may think of the semisuspension as `capping' the polytope off with the missing facet. Equivalently, given a regular $n-CW$-ball $\Gamma$, $\semisuspension \Gamma$ may be viewed as embedding $\Gamma$ in a hemisphere of the $n$-sphere.\\

\begin{figure}[h!]
\centering
 \begin{minipage}{.5\textwidth}
 \centering
\begin{tikzpicture}[thick,scale=.4]
\coordinate (p1) at (0,0);
\coordinate (p2) at (7,0);
\coordinate (p3) at (0,6);
\coordinate (p4) at (7,6);
\coordinate (p5) at (4,2);
\coordinate (p6) at (11,2);
\coordinate (p7) at (4,8);
\coordinate (p8) at (11,8);
\coordinate (c) at (6,3.5);
\draw (p1) -- (p5) ;
\draw (p5) -- (p6);
\draw (p1) -- (p2) -- (p6);

\fill[red,opacity=.5] (p1) -- (p2) -- (p6) -- (p5)--cycle;
\end{tikzpicture}

\caption*{$P_1(S)$}
\end{minipage}%
\begin{minipage}{.5\textwidth}
\centering
\begin{tikzpicture}[thick,scale=.4]
\coordinate (p1) at (0,0);
\coordinate (p2) at (7,0);
\coordinate (p3) at (0,6);
\coordinate (p4) at (7,6);
\coordinate (p5) at (4,2);
\coordinate (p6) at (11,2);
\coordinate (p7) at (4,8);
\coordinate (p8) at (11,8);
\coordinate (c) at (6,3.5);
\draw (p1) -- (p5) ;
\draw (p5) -- (p6);
\draw (p1) -- (p2) -- (p6);

\fill[red,opacity=.5] (p1) -- (p2) -- (p6) -- (p5)--cycle;
\fill[pur,opacity=.8] plot[smooth] coordinates{(p1) (c) (p6)} (p6) -- (p2) -- (p1);

\begin{scope}[thick,dashed, opacity=.5]
\draw plot[smooth] coordinates{(p2) (c)(p5)};
\draw plot[smooth] coordinates{(p1) (c) (p6)};
\end{scope}
\end{tikzpicture}
\caption*{$\semisuspension P_1(S)$}
\end{minipage}
\caption{A square $S$ and its semisuspension. We adjoin a single new purple $2$-face adjacent to all edges and vertices. Think of it as a membrane lying above the pink face.}
\end{figure}

Since the $ab$-index and $cd$-index are only defined for posets with $\widehat{1}$, it is of use to have an operation that adjoins $\widehat{1}$. 

\begin{definition}
 Given a complex $\Omega$, we let $P_1 (\Omega) = \Omega \cup \{\widehat{1}\}$ be the face poset of $\Omega$ with a maximal $\widehat{1}$ adjoined.
\end{definition}

We also may apply this operation to posets that already contain a $\widehat{1}$ with the newly adjoined $\widehat{1}$ taking precedence as the maximal element.

\begin{definition}
Let $Q$ be an Eulerian poset of rank $n$. Then $P_1(Q)$ is a near-Eulerian poset of rank $n+1$.
\end{definition}

\emph{Proof:} Let $\widehat{1}_Q$ and $\widehat{1}_P$ denote the maximal elements of $Q$ and $P_1(Q)$ respectively. Consider $\Sigma Q$ and recall that it is Eulerian. Removing either coatom of $\Sigma Q$ gives a poset isomorphic to $P_1(Q)$, so $P_1(Q)$ is near-Eulerian. This is easy to see via the Hasse diagram above. $\qed$\\
\begin{figure}
\begin{minipage}{.5\textwidth}
\centering

\begin{tikzpicture}
\draw [fill] (0,3) circle [radius = 2pt];
\draw [fill] (0,4) circle [radius = 2pt];
\draw (0,4) -- (0,3);
\draw [dashed] (0,3) -- (0,2);
\draw (0,0) ellipse [x radius = 1, y radius = 2];
\draw (0,0) node{$Q \setminus \{\widehat{1}\}$};
\draw (0,3) node[anchor=east]{$\widehat{1}_Q$};
\draw (0,4) node[anchor=east]{$\widehat{1}_P$};
\end{tikzpicture}
\caption{$P_1(Q)$} 

\end{minipage} %
\begin{minipage}{.5\textwidth}
\centering

\begin{tikzpicture}
\draw [fill] (.5,3) circle [radius = 2pt];
\draw [fill] (-.5,3) circle [radius = 2pt];
\draw [fill] (0,4) circle [radius = 2pt];
\draw (0,4) -- (.5,3);
\draw (0,4) -- (-.5,3);
\draw [dashed] (0.5,3) -- (0.5, 1.7);
\draw [dashed] (-0.5,3) -- (-0.5, 1.7);
\draw (0,0) ellipse [x radius = 1, y radius = 2];
\draw (0,0) node{$Q \setminus \{\widehat{1}\}$};
\draw (0,4) node[anchor=east]{$\widehat{1}$};
\end{tikzpicture}
\caption{$\Sigma Q$} 
\end{minipage}
\end{figure}
  
We now introduce a modified definition of shelling for regular $CW$-complexes that is required for Stanley's result. 
\begin{definition}
\cite{stanleycd} Let $\Omega$ be an Eulerian $n$ dimensional  regular $CW$-complex. We say that $\Omega$ is \emph{spherically shellable} ($S$-shellable for short) if either $\Omega$ is empty, or there exists a linear order of its facets $\sigma_1, \sigma_2 \dots \sigma_k$ that satisfies the following recursive definition: 
\begin{enumerate}
\item The boundary complex $\partial \sigma_1$ has a shelling.
\item For $1 < i \leq k -1$, let $\Gamma_i$ be the closure of the subcomplex of $\partial \sigma_i$ that is not contained in the previously shelled facets. That is,

\begin{equation}
\Gamma_i:= cl[\partial \sigma_i \setminus \left(\bigcup_{j=1}^{i-1} \overline{\sigma_j} \cap \overline{\sigma_i}  \right)] \, . 
\end{equation}
We then need that $\Gamma_i$ is a near-Eulerian $(n-1)$-$CW$-complex and that $\semisuspension \Gamma_i$ is $S$-shellable with initial facet $\tau$ (where $\tau$ is the new facet introduced by the semisuspension).
\end{enumerate}
\end{definition}

Any line shelling \cite{brugmani} of a polytope satisfies this definition. So any result on $S$-shellings holds for polytopes. However, not all $S$-shellings are shellings (via our original definition), and visa versa \cite{stanleycd}. We may also inductively find that this definition implies that each $\Gamma_i$ is a ball and that $\Omega$ is a sphere.  \\ 

We spend the remainder of this chapter giving an overview of Stanley's result that the $cd$-index of an $S$-shellable $CW$-complex is non-negative. This will soon be bogged down in technicalities, so let us state a brief overview here of our ultimate goal. We fix the following notation for the remainder of the chapter.\\

 Let $\Omega$ be an $n$-dimensional regular $CW$-sphere with an $S$-shelling $\sigma_1, \sigma_2, \dots \sigma_k$.\\
 
 For $1 \leq i \leq k - 1$ let $\Omega_i = \semisuspension P_1(\sigma_1 \cup \sigma_2 \cup \dots \cup \sigma_i)$. That is, $\Omega_i$ is the semisuspension of the complex induced by the first $i$ facets of the shelling. Note that $\Omega_{k-1} = P_1(\Omega)$.\\
 
  Given $\Omega_{i-1}$, let $\tau_{i-1}$ be the new face introduced by the semisuspension of $\sigma_1 \cup \dots \cup \sigma_{i-1}$. We may construct $\Omega_i$ from $\Omega_{i-1}$ by subdividing $\tau_{i-1}$ to `carve' $\sigma_i$ from the face. When subdividing $\Omega_{i-1}$ to $\Omega_i$, $\tau_{i-1}$ is subdivided into exactly two new facets: $F_i$ and $\tau_i$. Furthermore, the subcomplex $\Gamma_i = F_i \cap \tau_i$ coincides exactly with the definition of $\Gamma_i$ in $(1)$.\\

\begin{figure}[h!]
\centering
\begin{minipage}{.5\textwidth}
\centering
\begin{tikzpicture}[thick,scale=.6]
\coordinate (p1) at (0,0);
\coordinate (p2) at (7,0);
\coordinate (p3) at (0,6);
\coordinate (p4) at (7,6);
\coordinate (p5) at (4,2);
\coordinate (p6) at (11,2);
\coordinate (p7) at (4,8);
\coordinate (p8) at (11,8);
\coordinate (c) at (6,3.5);
\draw (p1) -- (p5) ;
\draw (p5) -- (p6);
\draw (p1) -- (p2) -- (p6);

\fill[red,opacity=.5] (p1) -- (p2) -- (p6) -- (p5)--cycle;
\fill[pur,opacity=.8] plot[smooth] coordinates{(p1) (c) (p6)} (p6) -- (p2) -- (p1);

\begin{scope}[thick,dashed, opacity=.5]
\draw plot[smooth] coordinates{(p2) (c)(p5)};
\draw plot[smooth] coordinates{(p1) (c) (p6)};
\end{scope}
\end{tikzpicture}
\caption*{$\Omega_1$}
\end{minipage}%
\begin{minipage}{0.5\textwidth}
\centering
\begin{tikzpicture}[thick,scale=.6]
\coordinate (p1) at (0,0);
\coordinate (p2) at (7,0);
\coordinate (p3) at (0,6);
\coordinate (p4) at (7,6);
\coordinate (p5) at (4,2);
\coordinate (p6) at (11,2);
\coordinate (p7) at (4,8);
\coordinate (p8) at (11,8);
\coordinate (c1) at (3, 5);
\coordinate (c2) at (8, 5);

\fill[red, opacity=.5] (p1)--(p5)--(p7)--(p8)--(p6)--(p2)--cycle;

\draw (p1) -- (p5) -- (p7);
\draw (p5) -- (p6);
\draw (p1) -- (p2) -- (p6);
\draw (p7) -- (p8);
\draw (p6) -- (p8);

\fill[pur, opacity=.8] plot[smooth] coordinates{(p1) (c1) (p7)} -- (p8) -- (p6) -- (p2) -- cycle;

\begin{scope}[thick,dashed, opacity=.5]
\draw plot[smooth] coordinates{(p1) (c1) (p7)};
\draw plot[smooth] coordinates{(p2) (c2) (p8)};
\draw plot[smooth] coordinates{(p6) (c2) (c1) (p5)};

\end{scope}

\end{tikzpicture}
\caption*{$\Omega_2$}
\end{minipage}

\vspace{.5cm}

\begin{minipage}{0.5\textwidth}

\centering

\begin{tikzpicture}[thick,scale=.6]
\coordinate (p1) at (0,0);
\coordinate (p2) at (7,0);
\coordinate (p3) at (0,6);
\coordinate (p4) at (7,6);
\coordinate (p5) at (4,2);
\coordinate (p6) at (11,2);
\coordinate (p7) at (4,8);
\coordinate (p8) at (11,8);
\coordinate (c1) at (4,1);
\coordinate (c2) at (6,4);
\coordinate (c3) at (9,3);

\fill[red, opacity=.5] (p1)--(p3) --(p7)--(p8)--(p6)--(p2)--cycle;

\draw (p1) -- (p5) -- (p7);
\draw (p5) -- (p6);

\draw (p1) -- (p2) -- (p6);
\draw (p3) -- (p7) -- (p8);
\draw (p1) -- (p3);

\draw (p6) -- (p8);

\fill[pur, opacity=.8] (p1)--(p3) --(p7)--(p8)--(p6)--(p2)--cycle; 

\begin{scope}[thick,dashed, opacity=.5]
\draw plot[smooth] coordinates{(p2) (c1) (p3)};
\draw plot[smooth] coordinates{(p2) (c2) (p7)};
\draw plot[smooth] coordinates{(p2) (c3) (p8)};

\end{scope}

\end{tikzpicture}
\caption*{$\Omega_3$}
\end{minipage}%
\begin{minipage}{0.5\textwidth}
\centering
\begin{tikzpicture}[thick,scale=.6]
\coordinate (p1) at (0,0);
\coordinate (p2) at (7,0);
\coordinate (p3) at (0,6);
\coordinate (p4) at (7,6);
\coordinate (p5) at (4,2);
\coordinate (p6) at (11,2);
\coordinate (p7) at (4,8);
\coordinate (p8) at (11,8);

\begin{scope}[thick,dashed,,opacity=0.6]
\draw (p1) -- (p5) -- (p7);
\draw (p5) -- (p6);
\end{scope}

\draw[fill=red,opacity=0.5] (p1)--(p3) --(p7)--(p8)--(p6)--(p2)--cycle;

\draw (p1) -- (p2) -- (p6);
\draw (p4) -- (p3) -- (p7) -- (p8);
\draw (p1) -- (p3);
\draw (p2) -- (p4);
\draw (p6) -- (p8);

\draw[fill=pur,opacity=0.5] (p2) -- (p4) -- (p3) -- (p7) -- (p8) -- (p6) -- cycle;

\end{tikzpicture}
\caption*{$\Omega_4$}
\end{minipage}

\vspace{.5cm}

\begin{tikzpicture}[thick,scale=.6]
\coordinate (p1) at (0,0);
\coordinate (p2) at (7,0);
\coordinate (p3) at (0,6);
\coordinate (p4) at (7,6);
\coordinate (p5) at (4,2);
\coordinate (p6) at (11,2);
\coordinate (p7) at (4,8);
\coordinate (p8) at (11,8);

\begin{scope}[thick,dashed,,opacity=0.6]
\draw (p1) -- (p5) -- (p7);
\draw (p5) -- (p6);
\end{scope}

\draw[fill=red,opacity=0.5] (p1)--(p3) --(p7)--(p8)--(p6)--(p2)--cycle;

\draw (p1) -- (p2) -- (p6);
\draw (p4) -- (p3) -- (p7) -- (p8);
\draw (p1) -- (p3);
\draw (p2) -- (p4);
\draw (p6) -- (p8);
\draw (p4) -- (p8);
\draw[fill=pur,opacity=0.5] (p3) -- (p4) -- (p8) -- (p7) -- cycle;

\draw (3.5,0) node[anchor=north]{$\Omega_5$};
\end{tikzpicture}
\caption{The $\Omega_i$'s of a cube with respect to our previous shelling. The face induced by the semisuspension is drawn in purple while the faces of the shelling are drawn in pink.}
\end{figure}

Finally, let

$$\check{\varPhi_i} = \Phi_{\Omega_i} - \Phi_{\Omega_{i-1}}$$

denote the difference in $cd$-index between $\Omega_i$ and $\Omega_{i-1}$. Note that this is well defined because both complexes are Eulerian posets.

\begin{lemma}
With $\Gamma_i$ as defined above, we may express $\checkcd_i$ as

$$\checkcd_i = \Phi_{\semisuspension \Gamma_i} \cdot c - \Phi_{\partial \Gamma_i} \cdot (c^2 - d) \, .$$
\end{lemma}

\emph{Proof:} This may be established by simple chain counting and clever variable substitution, and is explained in full in \cite{stanleycd}. \qed\\

We also may now compute the $cd$-index of a polytope in terms of complexes of strictly smaller dimension, giving us the following decomposition.

\begin{corollary}
Let $\Omega$ be a $CW$-complex with an $S$-shelling $\sigma_1, \sigma_2, \dots \sigma_k$. With notation defined as previous, we may express the $cd$-index of $\Omega$ as

$$\Phi_\Omega = \Phi_{\sigma_1} \cdot c + \sum_{i = 2}^{k-1} [\Phi_{\semisuspension \Gamma_i} \cdot c - \Phi_{\partial \Gamma_i} \cdot (c^2 - d)] \, . $$
\end{corollary}

\emph{Proof:} It is easy to see that
\begin{equation}
\Phi_\Omega = \Phi_{\Omega_{k-1}} = \Phi_{\Omega_1} + \sum_{i = 2}^{k-1} \checkcd_i \, .
\end{equation}

Note that given an Eulerian poset $Q$, $\semisuspension P_1(Q) = \Sigma Q$. So $\Omega_1 =  \Sigma \sigma_1$. Applying Corollary 3.9 to the first term gives us that $\varPhi_{\Omega_1} = \varPhi_{\sigma_1} \cdot c$. From Lemma 3.12, we know that $\checkcd_i = \Phi_{\semisuspension \Gamma_i} \cdot c - \Phi_{\partial \Gamma_i} \cdot (c^2 - d)$. Substitute each of these values into equation (2) and the result follows. \qed\\

This decomposition is enough to characterize the $cd$-index for two and three dimensional polytopes. 

\begin{lemma}
Consider a $2$-polytope $P$ with $n$ vertices. Then $\Phi_P = c^2 + (n-2)d$.\\

\end{lemma}

 \emph{Proof:} There exists some shelling $E_1, E_2, ..., E_n$ of $P$. Note that the face poset of any edge $E$ is isomorphic to the boolean lattice $B_2$.  So we have that $\Phi_{\Omega_1} = \Phi_{E_1}c = c^2$. 
 
 Each subsequent edge $E_2, ..., E_{n-1}$ is adjacent to exactly one other edge in the polygon and intersects it in exactly a single vertex. So for all $i$, $\Gamma_i$ is a single vertex, $\Phi_{\semisuspension P_1(\Gamma_i) } = c$ and $\Phi_{\partial \Gamma_i} = 1$. This gives us
 
 \begin{align*}
  \Phi_P =& \Phi_{E_1} \cdot c + \sum_{i = 2}^{n-1} [\Phi_{\semisuspension \Gamma_i} \cdot c - \Phi_{\partial \Gamma_i}(c^2 - d)]\\
  &= c^2 + \sum_{i = 2}^{n-1}[c^2 - (c^2 - d)]\\
  &= c^2 + (n - 2)d
 \end{align*}
 \qed
   
\begin{lemma}
Let $P$ be a 3-dimensional polytope. Then (with $f_i$ denoting the number of $i$-faces of $P$) we have $\Phi_P = c^3 + (f_0 - 2)dc + (f_2 - 2)cd$.
\end{lemma}

\emph{Proof:} There exists some shelling $F_1, F_2, .., F_k$ of $P$ where $k = f_2$. Let $v_i$ denote the number of vertices of $F_i$. Applying Lemma 3.16 to $F_1$ gives us that

\begin{equation}
\Phi_P = (c^2 + (v_1 - 2)d)c + \sum_{i=2}^{k-1}[\Phi_{\widetilde{\Sigma}\Gamma_i}c - \varPhi_{\partial \Gamma_i}(c^2 - d)] \, . 
\end{equation}
Note that $\partial \Gamma_i$ for $2 \leq i \leq k - 1$ is $0$-dimensional. Furthermore, $\Gamma_i$ is 1-dimensional and homemorphic to a ball. Specifically, $\Gamma_i$ is a path of edges.  So $\partial \Gamma_i$ is just two points and $\Phi_{\partial \Gamma_i} = c$. Let $v_{\Gamma_i}$ be the number of vertices of $\Gamma_i$. Then $\semisuspension \Gamma_i$ is a polygon with $v_{\Gamma_i}$ vertices and $\Phi_{\semisuspension \Gamma_i} = c^2 + (v_{\Gamma_i} - 2)d$. Substituting into equation (2), we may compute the following:

\begin{align*}
\Phi_P =& c^3 + (v_1 - 2)dc + \left(\sum_{i = 2}^{k-1}[c^2 + (v_{\Gamma_i} - 2)d]c - c(c^2 - d)\right)\\
=& c^3 + (v_1 - 2)dc + \left(\sum_{i=2}^{k-1} (v_{\Gamma_i} - 2)dc + cd \right)
\end{align*}

The key observation for our next step is that $v_1 + \sum_{i = 2}^{k-1}(v_{\Gamma_i} -2) = f_0$. This is because we may partition the vertices of $P$ into the sets of those in $F_1$, and those in the interior of each $\Gamma_i$ (since each $\Gamma_i$ meets the complex $F_1 \cup \dots \cup F_{i-1}$ in exactly its two boundary vertices). Recall that $k = f_2$ and we may then write

$$\Phi_P = c^3 + (f_0 - 2)dc + (f_2 - 2)cd \, .$$ \qed

\begin{theorem}
\emph{(Stanley)} \cite{stanleycd} Let $\Omega$ be an $S$-shellable regular $CW$-$n$-sphere. Then the coefficients of the $cd$-index of $\Omega$ are non-negative. That is, $\varPhi_\Omega \geq 0$. 
\end{theorem}

\emph{Proof:} This requires a finnicky double induction on both the number of facets of $\Omega$ and the dimension. See \cite{stanleycd} for proof.

\begin{corollary}
The $cd$-index of a polytope has non-negative coefficients. 
\end{corollary}

These results give us a pretty good picture of the non-negativity of the $cd$-index of geometric complexes. For the further case of non-geometric posets (as well as a wider class of $CW$-spheres), we turn to the theory of Gorenstein* posets and Gorenstein simplicial complexes. 

\clearpage
\section{Gorenstein Complexes}

We may use the underlying chain structure of a poset to associate to ita simplicial complex. That is, the chains of a poset $P$ can themselves be considered as a poset under inclusion and there exists a simplicial complex whose face poset is isomorphic. This lets us describe arbitrary posets in geometric terms and allows us to apply a toolbox of algebraic topology to the study of posets.

\begin{definition}

 Given a poset $P$, we define a simplicial complex $\Delta(P)$ (called the \emph{order complex} of $P$) as follows:
 
 \begin{itemize}
 \item The vertices of $\Delta(P)$ are exactly the elements of $P$, excluding $\widehat{0}$ and $\widehat{1}$ (if present).
 \item The faces of $\Delta(P)$ correspond bijectively with the nondegenerate chains of $P$. That is, for any chain $\widehat{0} < x_0 < x_1 < \dots x_k < \widehat{1}$ in $P$, there exists a $k$-simplex in $\Delta(P)$ with vertices $x_0, x_1, \dots x_k$ and visa versa.
 \end{itemize}
 
 \end{definition}
 
Note that if $P$ is the face poset of some complex, then $\Delta(P)$ is isomorphic to the face poset of the barycentric subdivision of $P \setminus \{\widehat{1}\}$. We also have that a poset $P$ is Eulerian if and only if $\Delta(P)$ is Eulerian \cite{borgkaru} \\

Also note that this definition of order complex implies a shift between dimension and rank. If $P$ is of rank $n + 1$ and contains $\widehat{1}$ and $\widehat{0}$, then $\Delta(P)$ has dimension $n-1$.\\

We also first need to define the operations of the \emph{star} and \emph{link} of a simplicial complex. One must take additional care when describing these operations for infinite complexes, but for the finite case the following definitions will suffice.    

\begin{definition}

Let $\Gamma$ be a pure $n$-dimensional simplicial complex, and let $\sigma \in \Gamma$ be an $i$-face. 

\begin{itemize}
\item The \emph{star} of $\sigma$ is the $n$-dimensional simplicial subcomplex induced by the faces that contain $\sigma$. That is,

$$\text{Star}_\Gamma (\sigma) = \{x \in \Gamma \, | \, \exists F \in \Gamma \text{ such that } \sigma \leq F, \text{ and } x \leq F\} \, . $$

\item The \emph{link} of $\sigma$ is the $(n-i-1)$-dimensional simplicial subcomplex given by all the faces of the Star$_\Gamma (\sigma)$ that do not intersect $\sigma$. That is, 

$$\text{Link}_\Gamma (\sigma) = \{x \in \text{Star}_\Gamma (\sigma) \, | \, \sigma \nleq x\} \, .$$
\end{itemize}

\end{definition}

This is easier to gather intuition for visually: see figure 9 for clarification. When the underlying complex is clear, we may omit the subscript.\\

\begin{figure}

\begin{minipage}{.5\textwidth}

\centering
\begin{tikzpicture}
\coordinate (p1) at (0,0);
\coordinate (p2) at (0,2);
\coordinate (p3) at (0,4);
\coordinate (p4) at (2,8);
\coordinate (p5) at (6,8);
\coordinate (p6) at (4,4);
\coordinate (p7) at (2,6);

\draw[fill=red,opacity=0.4] (p3)--(p4)--(p6)--cycle;

\draw [fill] (p1) circle [radius = 2pt]; 
\draw [fill] (p2) circle [radius = 2pt]; 
\draw [fill,red] (p3) circle [radius = 2pt]; 
\draw [fill,red] (p4) circle [radius = 2pt]; 
\draw [fill] (p5) circle [radius = 2pt]; 
\draw [fill,red] (p6) circle [radius = 2pt]; 
\draw [fill,red] (p7) circle [radius = 2pt]; 

\draw (p1) -- (p2) -- (p3);
\draw (p4) -- (p5) --(p6) -- (p1);
\draw (p2) -- (p6);
\draw [red] (p3) -- (p4);
\draw [red](p3) -- (p6);
\draw [red](p7) -- (p3);
\draw [red] (p7) -- (p4);
\draw [red] (p7) -- (p6);
\draw [red] (p4) -- (p6);

\draw (p7) node[anchor=east]{$v$};

\end{tikzpicture}
\end{minipage}%
\begin{minipage}{.5\textwidth}

\centering
\begin{tikzpicture}
\coordinate (p1) at (0,0);
\coordinate (p2) at (0,2);
\coordinate (p3) at (0,4);
\coordinate (p4) at (2,8);
\coordinate (p5) at (6,8);
\coordinate (p6) at (4,4);
\coordinate (p7) at (2,6);

\draw [fill] (p1) circle [radius = 2pt]; 
\draw [fill] (p2) circle [radius = 2pt]; 
\draw [fill,red] (p3) circle [radius = 2pt]; 
\draw [fill,red] (p4) circle [radius = 2pt]; 
\draw [fill] (p5) circle [radius = 2pt]; 
\draw [fill,red] (p6) circle [radius = 2pt]; 
\draw [fill] (p7) circle [radius = 2pt]; 

\draw (p1) -- (p2) -- (p3);
\draw (p4) -- (p5) --(p6) -- (p1);
\draw (p2) -- (p6);
\draw [red] (p3) -- (p4);
\draw [red](p3) -- (p6);
\draw (p7) -- (p3);
\draw (p7) -- (p4);
\draw (p7) -- (p6);
\draw [red] (p4) -- (p6);
\draw (p7) node[anchor=east]{$v$};
\end{tikzpicture}
\end{minipage}
\caption{The star and link of $v$, respectively.}
\end{figure}
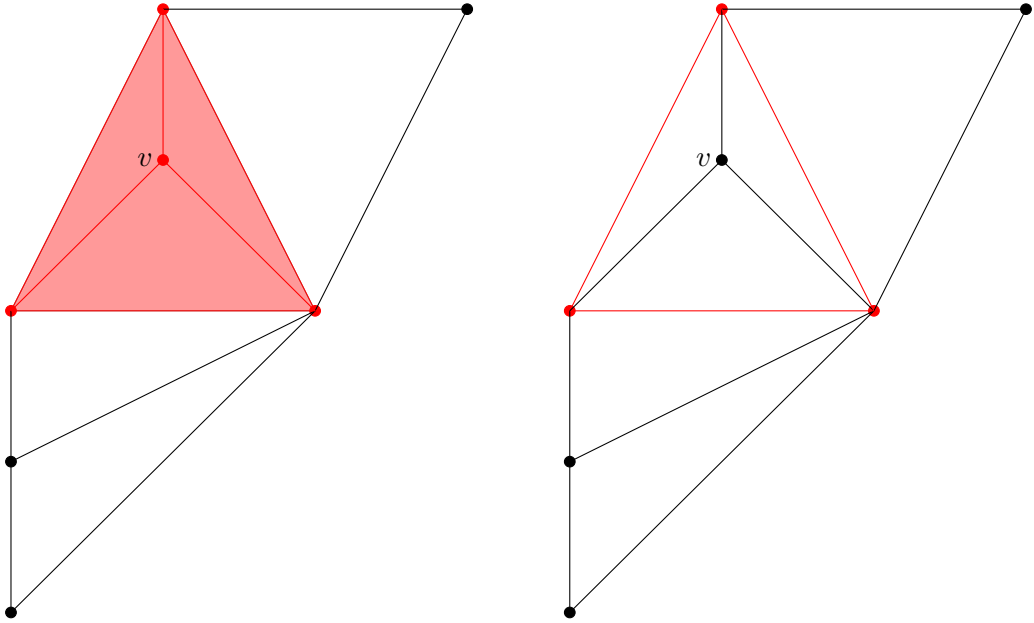

The link is a combinatorially useful structure. The face poset of $\text{Link}_\Gamma (\sigma)$ is isomorphic to the subposet $\{x \, | \,\sigma \leq x \}$ in $\Gamma$. That is, there is a correspondence between $\widehat{0} \in \text{Link}_\Gamma (\sigma)$ and $\sigma \in \{x \, | \,\sigma \leq x \}$, the zero dimensional faces of $\text{Link}_\Gamma (\sigma)$ and the $i + 1$ dimensional faces of $\{x \, | \,\sigma \leq x \}$, and so on. In other words, the link of $\sigma$ gives us a low dimensional picture of the faces that contain it by shifting those faces down $i + 1$ dimensions.\\

Gorenstein complexes are a class of well behaved simplicial complexes. They generalize simplicial spheres by being defined to have the same homology groups.  

\begin{definition}
We say that an $n$-dimensional simplicial complex $\Delta$ is \emph{Gorenstein} if $\Delta$ is a real homology sphere of dimension $n$. Explicitly, this means that we have that the real reduced homology of $\Delta$ is

\[\widetilde{H_i}(\Delta, \mathbb{R}) = \left\{
     \begin{array}{lr}
       \mathbb{R} & \text{if } i = n\\
       0 & \text{otherwise}
     \end{array}
   \right.\]
   and for any simplex $x \in \Delta$ with $\dim (x) = m$, we have that the reduced homology of its link is 

\[\widetilde{H_i}(\text{Link}_{\Delta} (x), \mathbb{R}) = \left\{
     \begin{array}{lr}
       \mathbb{R} & \text{if } i = n - m - 1\\
       0 & \text{otherwise}
     \end{array}
   \right. .\]
   \end{definition}

\begin{example}
 The boundary of any simplicial polytope or simplicial sphere is a Gorenstein complex. 
\end{example}

\begin{example}
 Let $\Delta$ be an $n$-dimensional Gorenstein complex with some $i$-face $\sigma$. Then Link$(\sigma)$ is a Gorenstein complex of dimension $(n-i-1)$.
 \end{example}

We say that a poset $P$ is a \emph{Gorenstein*} complex if $\Delta(P)$ is a Gorenstein complex. Note that if $P$ is Eulerian of rank $n + 1$, then $\Delta(P)$ has dimension $n-1$.\\

To be near-Gorenstein is to invoke a similar condition on a homology ball \cite{borgkaru}. 

\begin{definition}
 We say that a pair of simplicial complexes $(\Delta, \partial \Delta)$ is \emph{near-Gorenstein} if $\Delta$ is a real homology ball of dimension $n$ with boundary $\partial \Delta$. That is,

\begin{itemize}
\item The complex $\partial \Delta$ is Gorenstein of dimension $n - 1$.
\item For every $x \in \Delta$ of dimension $m$, the reduced homology of its link is

\[\widetilde{H_i}(\text{Link}_{\Delta} (x), \mathbb{R}) = \left\{
     \begin{array}{lr}
       \mathbb{R} & \text{if } i = n - m - 1 \text{ and } x \notin \partial \Delta \\
       0 & \text{otherwise}
     \end{array}
   \right. .\]

\end{itemize}
\end{definition}

Note that this implies that the link of any element of $ \partial \Delta$ has trivial homology. This includes $\Delta$ since by convention the empty face has Link$(\emptyset) = \Delta$.The motivating example for our purposes is that any full dimensional simplicial subdivision of an $n$-dimensional polytope $P$ (that is, a subdivision of its maximal face $|P|$ into an $n$-dimensional simplicial complex) is near-Gorenstein. \\

\begin{figure}[h!]
\centering
\begin{tikzpicture}[scale=2]
\draw (0,0)--(1,0)--(3,0)--(3,1)--(1,2)--(0,1)--(0,0);
\draw (0,1)--(1,0); 
\draw (1,0) --(2,1)--(1,2);
\draw (3,0)--(2,1)--(3,1);
\draw (1,2)--(1,0);
\draw [fill] (0,0) circle [radius = 1pt];
\draw [fill] (1,0) circle [radius = 1pt];
\draw [fill] (3,0) circle [radius = 1pt];
\draw [fill] (0,1) circle [radius = 1pt];
\draw [fill] (2,1) circle [radius = 1pt];
\draw [fill] (3,1) circle [radius = 1pt];
\draw [fill] (1,2) circle [radius = 1pt];
\end{tikzpicture}
\caption{A near-Gorenstein complex.}
\end{figure}
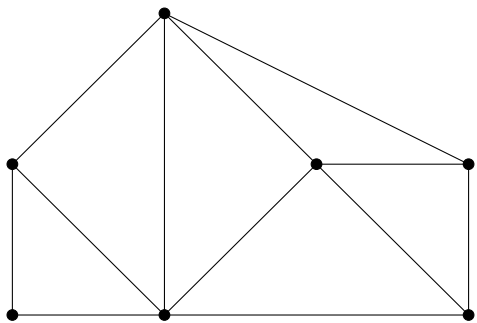

To define the analogous condition for posets, we must establish a notion of boundary for a poset.

\begin{definition}
Suppose $P$ is a near-Eulerian poset. Recall that we may uniquely obtain its semisuspension $\semisuspension P$ by adding a new coatom $\tau$ such that $x \leq_{\semisuspension P} \tau$ for any $x \in P$ such that $[x, \widehat{1}]_P$ is a three element chain. We then define the \emph{boundary} of $P$ to be 

$$\partial P:= P_1 ([\widehat{0}, \tau) ) \,.$$ 

This is well defined because $\partial P \subset P$. Also note that $\partial P$ is Eulerian. We refer to any element $x \in P \setminus \partial P$ as an element of (or lying in) the \emph{interior} of $P$.
\end{definition}

\begin{definition}
Let $Q$ be an Eulerian poset. Then we may write $\partial Q := Q \setminus \{\widehat{1}\}$.
\end{definition}

This coincides with the polytope definition of a boundary complex, but is not a notation used in literature. We include it here as a formality because it helps some of our summations look more nice.
 
 \begin{example}
 Consider the sets $P = \{\emptyset, \{1\}, \{2\}, \{3\}, \{1, 2\}, \{2, 3\}, \{1, 2, 3\} \}$  partially ordered under inclusion (where we have $\widehat{0} = \emptyset$ and $\widehat{1} = \{1, 2, 3\}$). Then $\partial P = \{ \{1\}, \{3\} \}$. 
 \end{example}

\begin{figure}[h!]

\begin{minipage}{.6\textwidth}
\centering
\begin{tikzpicture}[scale=.6]
\tikzstyle{every node}=[font=\small]
\draw [fill] (0,0) circle [radius = 2pt];
\draw [fill] (-4,4) circle [radius = 2pt];
\draw [fill] (4,4) circle [radius = 2pt];
\draw [fill] (-2,2) circle [radius = 2pt];
\draw [fill] (2,2) circle [radius = 2pt];
\draw (-4,4)-- (0,0) -- (4,4);
\draw (-4,4) node[anchor=south]{$\{1\}$};
\draw (4,4) node[anchor=south]{$\{3\}$};
\draw (0,0) node[anchor=north]{$\{2\}$};
\draw (2,2) node[anchor=west]{$\{2, 3\}$};
\draw (-2,2) node[anchor=east]{$\{1, 2\}$};
\draw (-3,3) node[anchor=east]{$\{1\} \leq \{1,2\}$};
\draw (-1,1) node[anchor=east]{$\{2\} \leq \{1,2\}$};
\draw (1,1) node[anchor=west]{$\{2\} \leq \{2,3\}$};
\draw (3,3) node[anchor=west]{$\{3\} \leq \{2,3\}$};
\end{tikzpicture}
\caption{$\Delta(P)$}
\end{minipage} %
\hfill\vline\hfill %
\begin{minipage}{.2\textwidth}
\centering
\begin{tikzpicture}[scale=.6]

\draw [fill] (-1,0) circle [radius = 2pt];
\draw [fill] (1,0) circle [radius = 2pt];
\draw (-1,0) node[anchor=south] {$\{1\}$};
\draw (1,0) node[anchor=south] {$\{3\}$};
\end{tikzpicture}
\caption{$\Delta (\partial P)$.}
\end{minipage}
\end{figure}

\begin{definition}
We say a pair of posets $(P, \partial P)$ is \emph{near-Gorenstein*} if $(\Delta (P), \Delta(\partial P))$ is a near-Gorenstein complex.
\end{definition}

Note that the above definition implies that $\partial P$ is a Gorenstein* poset. We may see Gorenstein* and near-Gorenstein* posets as a restriction of Eulerian and near-Eulerian posets to those with particularly well behaved order complexes. \\

We may decompose the $ab$-index of a near-Eulerian poset $P$ into a term that counts the chains of its boundary, and a remainder term. This is a particularly nice way to approach its chains because $\partial P$ is Eulerian.

\begin{lemma}
Given a near-Eulerian poset $P$ with boundary $\partial P$, we may write its $ab$-index as 

\[\Psi_{P} = \localab{P}\ + \Psi_{\partial P} \cdot a \, , \]

where $\localab{P} := \varPsi_{\semisuspension P} - \varPsi_{\partial P} \cdot (a + b)$. Furthermore, $\localab{P}$ is expressible in the variables $c = a + b$ and $d = ab + ba$. Note that $\localab{P}$ is implicitly defined by $P$.

\end{lemma}

\emph{Proof:} Consider the semisuspension $\semisuspension P$. $\semisuspension P$ adds a single new coatom $\tau$. The semisuspension contains chains of exactly two types: those contained in $P$,  and chains in $\partial P$ followed  by $\tau$. This gives us the following equations of the flag polynomial:

\begin{align*}
\varUpsilon_{\semisuspension P} &= \varUpsilon_{P} + \varUpsilon_{\partial P} \cdot b \, . \\
\varUpsilon_{P} &= \varUpsilon_{\semisuspension P} - \varUpsilon_{\partial P} \cdot b \, . 
\end{align*}
Substituting $a - b$ for $a$ gives us the $ab$-index, which may then be manipulated as follows.

\begin{align*}
\Psi_{P}  &= \Psi_{\semisuspension P} - \Psi_{\partial P} \cdot b  \\
 &= (\Psi_{\semisuspension P} - \Psi_{\partial P} \cdot (a + b)) + \Psi_{\partial P} \cdot a  \\
 &= \localab{P} +  \Psi_{\partial P} \cdot a  \, .
\end{align*}

 We omit the $a$ and $b$ arguments of the function $\localab{P}$ when context is clear and it is notationally convenient.\\
 
  To see that $\localab{P}$ may be expressed in $c = a + b$ and $d = ab + ba$, note that both $\semisuspension P$ and $\partial P$ are Eulerian and thus there exist $cd$-indexes $\varPhi_{\semisuspension}$ and $\varPhi_{\partial P}$. Therefore the polynomial $\localcd{P} (c,d):= \varPhi_{\semisuspension P} - \varPhi_{\partial P} \cdot c$ is well defined and we have 
  $$ \localcd{P} (a + b, ab + ba) = \localab{P}(a,b) \, . \qed$$

 We call $\localab{P}(a,b)$ and $\localcd{P}(c,d)$ the \emph{local $ab$-index} and \emph{local $cd$-index} of $P$ respectively. This result is also described in \cite{borgkaru} particularly for Gorenstein* posets, but we deviate from their terminology and notation for reasons that will become clear after chapter 8. Note that the definition implies for that the trivial single element poset $\{\widehat{0} \}$ we have $\localcd{\widehat{0}} = 1$. \\

We present the following definition analogous to the flag polynomial of a poset (Definition 2.10).

\begin{definition}
 We define the \emph{local flag polynomial} of $P$ by 

$$\localflag{P} (a,b) = \localab{P} (a + b, b) \, .$$
\end{definition}

In the expression $\Psi_{P} = \localab{P}\ + \Psi_{\partial P} \cdot a$, both $\localab{P}$ and $\Psi_{\partial P}$ may be rewritten in terms of $c$ and $d$. So $\Psi_{P}$ is `almost' a polynomial in $c$ and $d$. We define the following in accordance with Ehrenborg and Karu $\cite{borgkaru}$.

\begin{definition}
 Let $P$ be a near-Eulerian poset of rank $n+1$. Then we define the $cd$-index of $P$ to be

$$\varPhi_P := \localcd{P} + \varPhi_{\partial P} \, .$$
\end{definition}

Note that this is a non-homogeneous polynomial. The terms of $\localcd{P}$ are degree $n$, while the terms of $\varPhi_{\partial P}$ are rank $n-1$ (where $c$ contributes 1 to degree and $d$ contributes 2).

Analyzing the local flag polynomial gives us a more clear combinatorial picture of what exactly the local $ab/cd$ index is enumerating. 
 
 \begin{lemma}
  Given a near-Eulerian poset $P$, we have that $\localflag{P} = \varUpsilon_P - \varUpsilon_{P_1 (\partial P)}$.
  \end{lemma}
 
 \emph{Proof:}
\begin{align*}
\localab{P} &= \varPsi_{\semisuspension P} - \varPsi_{\partial P} \cdot (a + b) \\
\localflag{P} &= \varUpsilon_{\semisuspension P} - \varUpsilon_{\partial P} \cdot (a + 2b)\\
&= (\varUpsilon_{P} + \varUpsilon_{\partial P} \cdot b) - \varUpsilon_{\partial P} \cdot (a + 2b)\\
&= \varUpsilon_{P} - \varUpsilon_{\partial P} \cdot (a + b)\\
&= \localflag{P} = \varUpsilon_P - \varUpsilon_{P_1 (\partial P)} \, . \qed
\end{align*}

We may think of the local flag polynomial as a difference operator between $P$ and $P_1(\partial P)$. If $\Pi$ is the face poset of a near-Gorenstein complex, $P_1 (\partial \Pi)$ has a natural topological interpretation as the face poset of the $CW$-complex obtained from $\partial \Pi$ by attaching a maximal face to its boundary. In the next chapter, we examine the concept of the subdivision of a poset that is analogous to a topological subdivision. In general, we find that $P$ is a subdivision of $P_1(\partial P)$, and that $\localflag{P}$ counts the chains introduced by subdividing $P_1 (\partial P)$ to $P$.

\begin{example}
 Consider the following two dimensional near-Gorenstein complex $\Pi$. Let $\Pi' = P_1(\partial \Pi)$. Then, with $\Pi$ and $\Pi'$ as Figures 15 and 16, we have 

\begin{align*}
\localflag{\Pi} &= \Upsilon_\Pi - \Upsilon_{\Pi'} = 2aba + 2aab + 4bba + 4bab + 4abb + 8bbb.\\
\localab{\Pi} &= 2aba + 2aab + 2bba + 2bab.\\
\localcd{\Pi} &= 2cd.
\end{align*}
\end{example}
\begin{figure}[h!]
\centering
\caption{$\Pi$ \\ $\Upsilon_\Pi = aaa + 5baa + 7aba + 3aab + 14bba + 9bab + 9abb + 18bbb $}
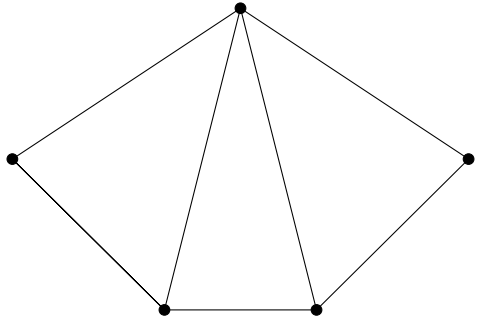
\begin{tikzpicture}[scale=2]
\draw (0,1)--(1,0)--(2,0)--(3,1)--(1.5,2)--(0,1)--(1,0);
\draw (1,0) -- (1.5,2)--(2,0);
\draw [fill] (0,1) circle [radius = 1pt];
\draw [fill] (1,0) circle [radius = 1pt];
\draw [fill] (2,0) circle [radius = 1pt];
\draw [fill] (1.5,2) circle [radius = 1pt];
\draw [fill] (3,1) circle [radius = 1pt];

\end{tikzpicture}
\end{figure}

\begin{figure}[h!]
\centering
\caption{$\Pi'$ \\ $\Upsilon_{\Pi'} = aaa + 5baa + 5aba + aab + 10bba + 5bab + 5abb + 10bbb $}
\begin{tikzpicture}[scale=2]
\draw (0,1)--(1,0)--(2,0)--(3,1)--(1.5,2)--(0,1);
\draw [fill] (0,1) circle [radius = 1pt];
\draw [fill] (1,0) circle [radius = 1pt];
\draw [fill] (2,0) circle [radius = 1pt];
\draw [fill] (1.5,2) circle [radius = 1pt];
\draw [fill] (3,1) circle [radius = 1pt];

\end{tikzpicture}
\end{figure}

The specification to Gorenstein* and near-Gorenstein* posets is of use because they generalize $CW$-spheres and $CW$-balls respectively. Ehrenborg and Karu proved non-negativity results on the $cd$-index of both Gorenstein* and near-Gorenstein* posets that are natural extensions of Stanley's result on the non-negativity of the $cd$-index of shellable CW-spheres. 

\begin{theorem}
\cite{karufans} The $cd$-index of a Gorenstein* poset has non-negative integer coefficients.
\end{theorem}

\begin{theorem}
\cite{bayerborg} The $cd$-index of a near-Gorenstein* poset has non-negative integer coefficients.
\end{theorem}

The proof of these results requires heavy algebraic machinery that is beyond the scope of this thesis. However, non-negativity is a particularly useful characterization for establishing bounds on the $cd$-index. We will examine some consequences of it in later chapters.  

\begin{corollary}
 The local $cd$-index of a near-Gorenstein* poset has non-negative integer coefficients.

\end{corollary}

\emph{Proof:} Let $P$ be a near-Gorenstein* poset of rank $n + 1$. Then $\varPhi_P = \localcd{P} + \varPhi_{\partial P}$. Note that $\localcd{P}$ is homogeneous of degree $n$  (where each $c$ contributes 1 and each $d$ contributes 2 to the degree) and $\varPhi_{\partial P}$ is homogenous of degree $n-1$. Since they have no shared terms and $\varPhi_P$ has non-negative coefficients, it follows that $\localcd{P}$ must have non-negative coefficients. $\qed$
\clearpage
\section{Skeletal Decomposition of a Subdivision}

We first present a traditional topological definition of subdivisions.

\begin{definition}
 Given face posets of regular CW complexes $\Pi$ and $\widehat{\Pi}$, we say that $\widehat{\Pi}$ is a \emph{subdivision} of $\Pi$ if there exists a surjective, order preserving map $\phi : \widehat{\Pi} \rightarrow \Pi$ such that for all $F \in \Pi$ the following conditions hold:

\begin{enumerate}
\item $\widehat{F} := \phi^{-1}([\widehat{0}, F])$ is a CW complex of the same dimension as $F$.
\item $\phi^{-1}(F)$ consists of all the interior cells of $\widehat{F}$.
\item $|F| = |\widehat{F}|$
\end{enumerate}
\end{definition}

If $\Pi$ and $\widehat{\Pi}$ are Gorenstein, with $\widehat{F}$ near Gorenstein for all $F \in \Pi$, then $\widehat{\Pi}$ is a \emph{Gorenstein subdivision} of $\Pi$. As an alternative notation not dependent on $\widehat{\Pi}$, given a subdivision $\phi : \Gamma \rightarrow \Pi$, for any $F \in \Pi$ we let $\Gamma_F := \phi^{-1}([\widehat{0}, F])$ denote the restriction of $\Gamma$ to $F$. 

\begin{example}
In Figures 14 and 15, $\Pi$ is a subdivision of $\Pi'$.
\end{example}

There are many generalizations of topological subdivisions to arbitrary posets. The following is of particular use to us. 

\begin{definition}
Given lower Eulerian posets $P$ and $\widehat{P}$ of rank $n$, we say $\widehat{P}$ is a \emph{strong Eulerian subdivision} of $P$ if there exists a surjective, order preserving map $\phi: \widehat{P} \rightarrow P$ such that for all $\sigma \in P$, 

\begin{itemize}
\item $P_1(\phi^{-1}[\widehat{0}, \sigma])$ is a near-Eulerian poset.
\item $\phi^{-1}[\widehat{0}, \sigma]$ has rank $\rho_P (\sigma)$.
\end{itemize}

\end{definition}

We may specialize even further by restricting the Eulerian condition to Gorenstein* \cite{borgkaru}.

\begin{definition}
Let $P$ and $\widehat{P}$ be Gorenstein* posets of rank $n$. Suppose we have a surjective, order preserving map $\phi: \widehat{P} \rightarrow P$. We say $\phi$ is a \emph{Gorenstein*} subdivision map if for every $\sigma \in P$, 

\begin{itemize}
\item $(\phi^{-1}[\widehat{0}, \sigma], \phi^{-1}[\widehat{0}, \sigma))$ is a near-Gorenstein* poset.
\item $\phi^{-1}[\widehat{0}, \sigma]$ has rank $\rho (\sigma)$.
\end{itemize}
\end{definition}

For shorthand, we let $\sigmahat: = P_1(\phi^{-1}[\widehat{0}, \sigma])$ and $\partial \sigmahat = P_1(\phi^{-1}[\widehat{0}, \sigma))$. This coincides with our earlier definition of the boundary of a near-Eulerian poset.

Notational aside: This unfortunately leads to the notation $\widehat{\widehat{1}_\Pi} = P_1(\phi^{-1}[\widehat{0}, \widehat{1}_{\Pi}]) = P_1(\Pi)$. However, in practice such an element is rarely acknowledged because the existence of $\widehat{1}$ in a strongly Eulerian subdivision between $\Pi$ and $\Pihat$ is mostly a formality; note that the definition of strongly Eulerian implies that $\phi(\sigma) = \widehat{1}_{\Pi}$ if and only if $\sigma = \widehat{1}_{\Pihat}$. Furthermore it is easy to calculate that for any such strongly Eulerian subdivision, $\localcd{\widehat{\widehat{1}_\Pi}} = 0$. Since the $cd$-index and $ab$-index do not technically count chains including $\widehat{1}$, its presence is often implicit in literature (to the annoyance of this writer). We aim to clarify its presence at the expense of clarity in notation. \\

 It is easy to see that a simplicial subdivision between Gorenstein complexes is both a Gorenstein* subdivision and a strongly Eulerian subdivision.

\begin{definition}
 Given a [Gorenstein*/simplicial/Eulerian] subdivision $\phi: \Pihat \rightarrow \Pi$ and some face $\sigma \in \Pihat$ we define $\phi(\sigma) = \tau \in \Pi$ to be the \emph{carrier} of $\sigma$ in $\Pi$. We say that $\sigma$ is \emph{carried} by $\tau$.
\end{definition}

In the simplicial case, we may geometrically see $\tau$ as the minimal face of $\Pi$ that topologically contains $\sigma$. That is, $|\sigma| \subseteq |\tau|$.\\

\begin{figure}
\begin{minipage}{.5\textwidth}
\centering
\begin{tikzpicture}
\draw [fill] (0,0) circle [radius = 2pt];
\draw [fill] (4,0) circle [radius = 2pt];
\draw [fill] (2,4) circle [radius = 2pt];
\draw (0,0) -- (2,4) -- (4,0);
\draw [red] (0,0) -- (4,0);
\draw (2,0) node[anchor=south]{$\tau$};
\end{tikzpicture}
\caption*{$\Pi$}
\end{minipage} %
\begin{minipage}{.5\textwidth}
\centering
\begin{tikzpicture}
\draw [fill] (0,0) circle [radius = 2pt];
\draw [fill] (4,0) circle [radius = 2pt];
\draw [fill] (2,4) circle [radius = 2pt];
\draw [fill] (2,0) circle [radius = 2pt];
\draw [fill] (1,2) circle [radius = 2pt];
\draw [fill] (3,2) circle [radius = 2pt];
\draw [fill] (2,2) circle [radius = 2pt];
\draw (0,0) -- (2,4) -- (4,0);
\draw (2,2) -- (0,0);
\draw (2,2) -- (1,2);
\draw (2,2) -- (2,4);
\draw (2,2) -- (3,2);
\draw (2,2) -- (4,0);
\draw (2,2) -- (2,0);
\draw (0,0) -- (2,0);
\draw [red] (2,0) -- (4,0);
\draw (3,0) node[anchor=south]{$\sigma$};

\end{tikzpicture}
\caption*{$\widehat{\Pi}$}
\end{minipage}
\caption{$\widehat{\Pi}$ is a simplicial subdivison of $\Pi$. The face $\sigma \in \widehat{\Pi}$ is carried by $\tau \in \Pi$.}
\end{figure}

Our motivation for the remainder of this chapter, along with most of chapter six, is to find a combinatorial description and elementary proof of the Gorenstein* subdivision $cd$-index decomposition developed by Ehrenborg and Karu in \cite{borgkaru}. Restricting to a special case of Gorenstein simplicial complexes, we find an elegant solution via a process we call `skeletal decomposition' of a subdivision. This gives a clear geometric decomposition when applied to face posets, and generalizes the formula to strongly Eulerian subdivisions.

\begin{definition}
Let $\phi$ be a strong Eulerian subdivision of rank $n$ Eulerian posets $\phi: \widehat{\Pi} \rightarrow \Pi$. We define the \emph{skeletal posets} $\Pi_i$ for $0 \leq i \leq n$ as follows:

\begin{itemize}
 \item $\Pi_n ;= \widehat{\Pi}$ .
 \item $\Pi_{0} ;= \Pi $.
  \item Let $\Pi |_{\geq i + 1}:= \{\sigma \in \Pi: \rho(\sigma) \geq i + 1 \}$ and let $\Pihat |_{\leq i}:= \{\sigma \in \Pihat: \rho( \phi (\sigma)) \leq i \}$. For $0 < i < n$, we define $$\Pi_i := \Pi |_{\geq i + 1} \cup \Pihat |_{\leq i}$$ to be a poset where for $\sigma, \tau \in \Pi_i$ we have the order relation $\sigma \leq \tau$ if and only if exactly one of the following holds:
  
\begin{enumerate}
\item  $\sigma, \tau \in \Pi |_{\geq i + 1}$ and $\sigma \leq \tau$ in $\Pi$.
\item  $\sigma, \tau \in \Pihat |_{\leq i}$ and  $\sigma \leq \tau$ in $\Pihat$.
\item $\sigma \in \Pihat |_{\leq i}, \tau \in \Pi |_{\geq i + 1}$ and $\phi(\sigma) \leq \tau$ in $\Pi$.
\end{enumerate}
  
  \end{itemize}
  \end{definition}
  In other words, for elements exclusively of rank $\geq i + 1$ or exclusively of rank $\leq i$ then their order relation is inherited from $\Pi$ and $\Pihat$ respectively. If $\rho(\sigma) \leq i$ and $\rho(\tau) \geq i + 1$ we have that $\sigma \leq \tau$ if and only if the carrier of $\sigma$ is ordered beneath $\tau$ in $\Pi$. Note that $\Pi_i$ is a well defined poset of rank $n$ with a natural rank function. For any $\sigma \in \Pi_i$, we have
  
 $$\rho(\sigma) = 
 \begin{cases}
\rho_{\Pihat} (\sigma) \text{ if } \sigma \in \Pihat |_{\leq i}\\
\rho_{\Pi} (\sigma) \text{ if } \Pi |_{\geq i+1}
 \end{cases}
 $$ 
  We then define \emph{skeletal maps} $\phi_i: \Pi_{i+1} \rightarrow \Pi_i$ for $0 \leq i \leq n-1$ by the following.

 $$
 \phi_i(\sigma) = \left\{
      \begin{array}{lr}
        \phi(\sigma) & \text{if } \sigma \in \Pihat|_{\leq i} \text{ and } \rho( \phi(\sigma)) = i + 1 \\
        \sigma & \text{otherwise}
      \end{array}
    \right.
 $$

Note that the map only has non-identity outputs for $\sigma \in \Pihat |_{\leq i} \subseteq \Pi_{i+1}$, that is, for $\sigma$ carried by some element of rank $i + 1$ in $\Pi$. 

\begin{lemma}
Let $\Pi, \Pihat, \phi, \phi_i$ and $\Pi_i$ be defined as above. Then $$\phi_0 \circ \phi_1 \circ \dots \circ \phi_{n-2} \circ \phi_{n-1} = \phi \,. $$ 
\end{lemma}
\emph{Proof:} This follows from careful construction of the skeletal posets and maps.\\

 Consider $\sigma \in \Pihat$. Fix the rank of its carrier, $\rho_{\Pi} (\phi(\sigma)) = k$. 

Then by definition, for all $k \leq i \leq n$, $\sigma \in \Pi_i$. Furthermore, for any $\phi_i: \Pi_{i+1} \rightarrow \Pi_i$ with $k \leq i \leq n - 1$, $\phi_i ( \sigma) = \sigma$ because $\sigma$ is carried by an element of rank $k$. So $\phi_k \circ \phi_{k+1} \circ \dots \circ \phi_{n-1} (\sigma) = \sigma$.

By definition of $\phi_{k-1}$, we have that $\phi_{k-1} (\sigma) = \phi(\sigma)$. Similarly, since $\phi(\sigma)$ has rank $k$ we have that $\phi(\sigma) \in \Pi_{i}$ for all $0 \leq i \leq k - 1$ and that $\phi_i(\phi(\sigma)) = \phi(\sigma)$ for all $0 \leq i \leq k - 2$. So $\phi_0 \circ \phi_{1} \circ \dots \circ \phi_{k-2} (\phi(\sigma)) = \phi(\sigma)$. Therefore

\begin{align*}
\phi_0 \circ \phi_1 \circ \dots \circ \phi_{n-2} \circ \phi_{n-1} (\sigma) &= \phi_0 \circ \phi_1 \circ \dots \circ \phi_{k-2} \circ \phi_{k-1}(\sigma)\\
&= \phi_0 \circ \phi_1 \circ \dots \circ \phi_{k-2} (\phi(\sigma))\\
&= \phi(\sigma) \, .
\end{align*} 

Since this holds for all $\sigma \in \Pihat$, it follows that $\phi_0 \circ \phi_1 \circ \dots \circ \phi_{n-2} \circ \phi_{n-1} = \phi$. \qed \\

 It follows that the following commutative diagram holds. We call the process of identifying the skeletal maps and posets of $\phi: \Pihat \rightarrow \Pi$ the \emph{skeletal decomposition} of $\phi$. 

$$
\xymatrix@1{\widehat{\Pi} = \Pi_n\ar[r]^{\phi_{n-1}} \ar@{-->}@/_1pc/[0,4]_\phi & \Pi_{n-1} \ar[r]^{\phi_{n-2}} & \dots \ar[r]^{\phi_1} & \Pi_{1} \ar[r]^{\phi_{0}} & \Pi_{0} = \Pi }
$$

 It is also convenient to define the following maps.
 
 $$\toskel_i : \Pihat \rightarrow \Pi_i \hspace{2cm}\toskel_i = \phi_{n-1} \circ \phi_{n-2} \circ \dots \circ \phi_i$$
 
 $$\fromskel_i : \Pi_i \rightarrow \Pi \hspace{2cm}\fromskel_i = \phi_{i-1} \circ \phi_{i-2} \circ \dots \circ \phi_{0}$$
 
 This gives us
  $$\phi = \phi_{n-1} \circ \dots \circ \phi_{0} = \toskel_i \circ \fromskel_i \text{ for } 0 \leq i \leq n - 1$$.
 
The maps $\toskel_i$ and $\fromskel_i$ are induced by $\Pi_i$. We may think of $\fromskel_i$ as the map that subdivides $\Pi$ to $\Pi_i$, and $\toskel_i$ as the map that subdivides $\Pi_i$ to $\Pihat$. We have not actually proved that $\fromskel_i$ and $\toskel_i$ are formal subdivisions, but that is ultimately irrelevant to our use of them and they are convenient to think of as such. \\ 

Though seemingly obtuse when applied to general posets, this process is rooted in strong geometric intuition. Suppose that $\Pi$ and $\Pihat$ are simplicial $n$-dimensional complexes such that $\Pihat$ is a subdivision of $\Pi$. $\Pi_i$ is obtained from $\Pi$ by restricting its subdivision to $\widehat{\Pi}$ to the elements $\sigma \in \Pi$ with $\rho(\sigma) \leq i$. That is, we subdivide the $(i-1)$-faces of $\Pi$ and leave the higher dimensional faces untouched (recall that in simplicial/CW complex, elements of rank $i$ correspond to faces of dimension $(i-1)$). Additionally, we can see that $\Pi_{i+1}$ is obtained from $\Pi_{i}$ by subdividing exactly the $i$-faces of $\Pi_{i}$. In other words, we can see skeletal decomposition of a subdivision as the process of first subdividing the $0$-faces of $\Pi$, then its $1$-faces, and so on until we subdivide its $n$-faces and achieve $\Pihat$. Significantly, though each $\Pi_i$ is the poset of a geometrically realizable complex, they are not necessarily simplicial complexes (they are best visualized as $CW$-complexes)! 

\begin{example}
Let $\Pi$ be a tetrahedron and let $\widehat{\Pi}$ be achieved by subdividing two of its facets in half, with two new edges meeting at a single new vertex. Note $\Pi_2$ is no longer a simplicial complex since the two facets adjacent to the subdivided edge are now both combinatorially equivalent to a quadrilateral. 
\end{example}
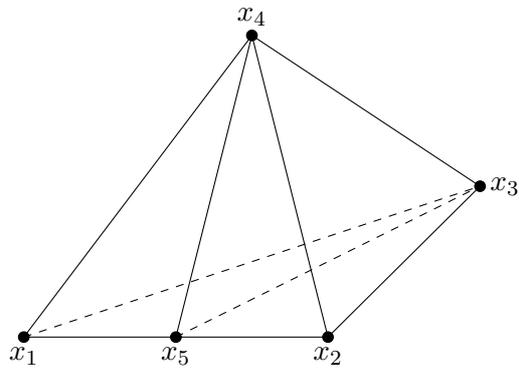
\begin{figure}[h!]
\caption{$\Pi = \Pi_{0} = \Pi_1$}
\centering
\begin{tikzpicture}[scale=2]
\draw (0,0)--(2,0)--(3,1)--(1.5,2);
\draw (0,0)--(1.5,2)--(2,0); 
\draw [dashed] (0,0) --(3,1);
\draw [fill] (0,0) circle [radius = 1pt];
\draw [fill] (2,0) circle [radius = 1pt];
\draw [fill] (3,1) circle [radius = 1pt];
\draw [fill] (1.5,2) circle [radius = 1pt];
\draw (0,0) node[anchor=north]{$x_1$};
\draw (2,0) node[anchor=north]{$x_2$};
\draw (3,1) node[anchor=west]{$x_3$};
\draw (1.5,2) node[anchor=south]{$x_4$};
\end{tikzpicture}
\end{figure}

\begin{figure}[h!]
\caption{$\Pi_2$}
\centering
\begin{tikzpicture}[scale=2]
\draw (0,0)--(2,0)--(3,1)--(1.5,2);
\draw (0,0)--(1.5,2)--(2,0); 
\draw [dashed] (0,0) --(3,1);
\draw [fill] (0,0) circle [radius = 1pt];
\draw [fill] (2,0) circle [radius = 1pt];
\draw [fill] (3,1) circle [radius = 1pt];
\draw [fill] (1.5,2) circle [radius = 1pt];
\draw [fill] (1,0) circle [radius = 1pt];
\draw (0,0) node[anchor=north]{$x_1$};
\draw (1,0) node[anchor=north]{$x_5$};
\draw (2,0) node[anchor=north]{$x_2$};
\draw (3,1) node[anchor=west]{$x_3$};
\draw (1.5,2) node[anchor=south]{$x_4$};
\end{tikzpicture}
\end{figure}

\begin{figure}[h!]
\caption{$\Pi_3 = \Pi_4 = \widehat{\Pi}$}
\centering
\begin{tikzpicture}[scale=2]
\draw (0,0)--(2,0)--(3,1)--(1.5,2);
\draw (0,0)--(1.5,2)--(2,0); 
\draw (1,0)--(1.5,2);
\draw [dashed] (1,0)--(3,1);
\draw [dashed] (0,0) --(3,1);
\draw [fill] (0,0) circle [radius = 1pt];
\draw [fill] (2,0) circle [radius = 1pt];
\draw [fill] (3,1) circle [radius = 1pt];
\draw [fill] (1.5,2) circle [radius = 1pt];
\draw [fill] (1,0) circle [radius = 1pt];
\draw (0,0) node[anchor=north]{$x_1$};
\draw (1,0) node[anchor=north]{$x_5$};
\draw (2,0) node[anchor=north]{$x_2$};
\draw (3,1) node[anchor=west]{$x_3$};
\draw (1.5,2) node[anchor=south]{$x_4$};
\end{tikzpicture}
\end{figure}
\newpage

We may divide the flags in any $\Pi_i$ of a skeletal decomposition to be of three distinct combinatorial categories. 

\begin{definition}
Let $C = F_0 < F_1 < \dots < F_k$ be a flag of $\Pi_i$.

\begin{enumerate}
\item If $F_0, \dots, F_k \in \Pi|_{\geq i + 1}$, we say that $C$ is an \emph{old} flag.

\item If $F_0, \dots, F_k \in \Pihat|_{\leq i}$, we say that $C$ is a \emph{new} flag. If $\rho ( \fromskel_i (F_k)) = m$ (that is, $F_k$ is carried by a rank $m$ element of $\Pihat$) we say that $C$ \emph{switches} at rank $m$. 

\item If $C = F_0 < \dots < F_j < F_{j+1}< \dots < F_k$ where $F_0, \dots, F_j \in \Pihat|_{\leq i}$ and $F_{j+1}, \dots, F_k \in \Pi|_{\geq i + 1}$ we say that $C$ is a \emph{mixed} flag. If $\rho (\fromskel_i) (F_j) = m$ we say that $C$ \emph{switches} at rank $m$.
\end{enumerate} 
\end{definition}

In other words, a flag of $\Pi_i$ switches at $m$ if its maximal element in $\Pihat|_{\leq i}$ is carried by a rank $m$ element of $\Pi$. Additionally, we will also say that a flag switches at $\sigma$ if $\sigma \in \Pi$ carries its maximal face in $\Pihat|_{\leq i}$. Note that by construction, any mixed or new flag of $\Pi_i$ must switch at some $m \leq i$ because for any $\sigma \in \Pihat_i |_\leq{i}$ we have  $ \rho(\phi(\sigma)) \leq i$.

\begin{example}
Consider $\Pi_2$ of figure 19. We will refer to faces by their vertex sets. $\{x_5\} < \{x_2, x_5\}$ is a new chain. $\{x_5\} < \{x_2, x_5\} < \{x_1, x_2, x_4, x_5\} $ is a mixed chain that switches at rank $2$. 
\end{example}

Note that even though $\{x_1, x_2, x_4, x_5\}$ contains a vertex not present in $\Pi$, we still count it as an `old' face of $\Pi_i$ because its interior is identical to the interior of the face $\{x_1, x_2, x_4\}$ of $\Pi$ (recall that in a $CW$-complex, faces are defined as open cells). \\

For the remainder of this paper, when not specified we will assume some arbitrary subdivision between $n$-dimensional Eulerian posets $\phi: \Pihat \rightarrow \Pi$. \\

Note that we may naturally restrict any subdivision map $\phi: \Pihat \rightarrow \Pi$ to a face $F \in \Pi$. To make this operation well defined on $\Pi_i$ for all $i$, we introduce the following family of restriction maps.

\begin{definition}
For all $0 \leq i \leq n$ and all $F \in \Pi$ we define the restriction of $\Pihat_i$ to $F$ by the following map:

$$\chi|^i_F: \Pi_i \rightarrow (\fromskel_i)^{-1}([\widehat{0}, F])$$

$$
\chi|^i_F(\sigma) =
\begin{cases}
        \sigma \text{ if } \sigma \in (\fromskel_i)^{-1}([\widehat{0}, F]),  \\
        \emptyset \text{ otherwise.}
\end{cases}
$$
    
$\chi|^i_F$ is just a characteristic function for the faces of $\Pi_i$ that map to an element of $[\widehat{0}, F]$. In the topological case, given a face $F$ of $\Pi$ the function identifies all faces of $\Pi_i$ that are contained in $F$ (or its boundary) and maps all other faces of $\Pi_i$ to the empty set.
\end{definition}
    
\begin{lemma}
Given a fixed face $F \in \Pi$ and $0 \leq i \leq n -1$, the following diagrams commute.

$$
\xymatrix@=50pt{
\Pihat = \Pi_n \ar[r]^{\toskel_i} \ar[d]_{\chi|^n_F} & 
\Pi_i \ar[d]^{\chi|^i_F} \\
\Pi_n \ar[r]_{\toskel_i} & \Pi_i
}
$$

$$
\xymatrix@=50pt{
\Pi_i\ar[r]^{\fromskel_i} \ar[d]_{\chi|^i_F} & 
\Pi_{0} \ar[d]^{\chi|^{0}_F} \\
\Pi_i \ar[r]_{\fromskel_i} & \Pi_{0} = \Pi
}
$$
\end{lemma}

\emph{Proof:} We may demonstrate the first diagram commutes as follows.\\

Take $\sigma \in \Pihat$. If $\sigma \notin \phi^{-1}([\widehat{0}, F])$ it holds trivially.

So suppose  $\sigma \in \phi^{-1}([\widehat{0}, F])$ . If $\rho( \phi (\sigma) ) \leq i$ we have

$$\toskel_i \circ \chi|^n_F (\sigma) = \toskel_i (\sigma) = \sigma$$

and

$$\chi|^i_F \circ \toskel_i (\sigma) = \chi|^i_F (\sigma) = \sigma$$
and the equation holds. \\

If $\rho( \phi(\sigma)) \geq i + 1$ then

$$\toskel_i \circ \chi|^n_F (\sigma) = \toskel_i(\sigma) =  \phi(\sigma) \in \Pi|_{\geq i + 1} \subset \Pi_i \, .$$

$$\chi|^i_F \circ \toskel_i (\sigma) = \chi|^i_F (\phi (\sigma))$$

Since $\phi (\sigma) \in [\widehat{0},F]$ by definition, we must have $\phi(\sigma) \in (\fromskel_i)^{-1}([\widehat{0}, F])$ and $\chi|^i_F (\phi (\sigma)) = \phi(\sigma)$. So the equation follows and the diagram commutes.

The proof that

$$\chi|^{0}_F \circ \fromskel_i = \fromskel_i \circ \chi|^{i}_F$$

is entirely symmetrical. \qed.
\clearpage
\section{Subdividing the cd-index}

Given a Gorenstein* subdivision $\phi: \Pihat \rightarrow \Pi$, can we express the $cd$-index of $\Pihat$ in terms of the elements of $\Pi$? Ehrenborg and Karu answered this question in the affirmative in Theorem 2.7 of \cite{borgkaru}. They achieve the result by using sheaf theory on fans and poset. 

By using the theory of skeletal decompositions we may establish this result purely combinatorially by careful chain counting, sidestepping the commutative algebra. Additionally, we extend the subdivision decomposition to hold for any strong Eulerian subdivision. We introduce the theorem below and prove it over the remainder of the chapter.

\begin{theorem}
Let $\phi : \widehat{\Pi} \rightarrow \Pi$ be a strong Eulerian subdivision map between rank $n$ posets. Then

\[\Phi_{\widehat{\Pi}} = \sum_{\sigma \in \Pi} \localcd{\widehat{\sigma}} \cdot \Phi_{[\sigma, \widehat{1}]} \]  
\end{theorem}
\emph{Proof:} We claim that for $i \geq 0$,  $\sum_{\substack{\sigma \in \Pi \\ \rho(\sigma) = i}} \localflag{\widehat{\sigma}} \cdot \Upsilon_{[\sigma, \widehat{1}]}$ counts exactly the new chains obtained by subdividing $\Pi_{i-1}$ to $\Pi_i$. That is, 

\begin{equation}
\sum_{\substack{\sigma \in \Pi \\ \rho(\sigma) = i}} \localflag{\widehat{\sigma}} \cdot \Phi_{[\sigma, \widehat{1})} = \Upsilon_{\Pi_i} - \Upsilon_{\Pi_{i-1}} \, .
\end{equation}

Equation (4) will be proven in Lemma 6.2. Assuming its validity (for now), we have 

$$ \sum_{\sigma \in \Pi} \localflag{\widehat{\sigma}} \cdot \Upsilon_{[\sigma, \widehat{1})} = \Upsilon_{\Pi_{0}} + (\Upsilon_{\Pi_1} - \Upsilon_{\Pi_{0}}) + (\Upsilon_{\Pi_2} - \Upsilon_{\Pi_1}) + \dots + (\Upsilon_{\Pi_n} - \Upsilon_{\Pi_{n-1}}) = \Upsilon_{\widehat{\Pi}} \, .$$

Since we may freely convert between the flag polynomial and $cd$ indexes by a variable substitution, it follows that

$$\sum_{\sigma \in \Pi } \localcd{\sigmahat} \cdot \flagtop = \Upsilon_{\Pihat} $$

if and only if

$$\sum_{\sigma \in \Pi } \localcd{\sigmahat} \cdot \cdtop = \Phi_{\Pihat} \, . $$

The theorem follows. \qed\\

We now establish the truth of equation $(4)$.

\begin{lemma} For some fixed $1 \leq i \leq n$, the following holds: 

$$\Upsilon_{\Pi_i} - \Upsilon_{\Pi_{i-1}} = \sum_{\substack{\sigma \in \Pi \\ \rho(\sigma) = i }} \localflag{\sigmahat} \cdot \flagtop = \sum_{\substack{C \text{ a chain of } \Pi_i \\ C \text{ switches at rank }i}} \alpha^C \, \, \, \, \, \, \, \, \, \, - \sum_{\substack{C \text{ a chain of } \Pi_{i-1} \\ C \text{ contains an element of rank }i}} \alpha^C $$

\end{lemma}

\emph{Proof:} We demonstrate that each of $\Upsilon_{\Pi_i} - \Upsilon_{\Pi_{i-1}}$ and $ \sum_{\substack{\sigma \in \Pi \\ \rho(\sigma) = i }} \localflag{\sigmahat} \cdot \flagtop$ equals the right hand side of the equation to prove the lemma.\\ 

\noindent
\text{\textbf{Claim 1:} }
$$
\Upsilon_{\Pi_i} - \Upsilon_{\Pi_{i-1}} = \sum_{\substack{C \text{ a chain of } \Pi_i \\ C \text{ switches at rank }i}} \alpha^C \, \, \, \, \, \, \, \, \, \, - \sum_{\substack{C \text{ a chain of } \Pi_{i-1} \\ C \text{ contains an element of rank }i}} \alpha^C $$

That is, $\Upsilon_{\Pi_i} - \Upsilon_{\Pi_{i-1}}$ records exactly the number of flags of $\Pi_i$ that switch at rank $i$, minus the flags of $\Pi_{i-1}$ that contain a rank $i$-element.\\

First, recall the definition of $\Upsilon_{\Pi_i}$ and note that we may write 

$$ \Upsilon_{\Pi_i} - \Upsilon_{\Pi_{i-1}} = \sum_{C \text{ a chain of } \Pi_i} \alpha^C - \sum_{C \text{ a chain of } \Pi_{i-1}} \alpha^C \, .$$  

Our goal is to cancel chains from this difference until only claim 1 remains.\\

First, we can see that the old flags of $\Pi_i$ are in bijective correspondence with the old flags $F_0 <\dots < F_k $ of $\Pi_{i-1}$ where $\rho_{\Pi_{i-1}} (F_0) \geq i + 1$ because $\phi_{i-1}: \Pi_i \rightarrow \Pi_{i-1}$ acts as the identity for any $\sigma \in \Pi_i$ with $ \rho (\phi(\sigma)) > i$.  It follows that such flags are eliminated in the difference.

Next, we may note a bijective correspondence between the flags of $\Pi_i$ that switch at some $j \leq i - 1$ and the flags of $\Pi_{i-1}$ that switch at $j$ and do not contain an element of rank $i$. Note that any such flag is present in both $\Pi_i$ and $\Pi_{i-1}$ because it decomposes into an initial segment contained in $\Pihat |_{\leq i - 1}$ followed by a (possibly empty) segment in $\Pi |_{\geq i + 1}$. So these flags are eliminated in the difference.

Since all mixed and new flags of $\Pi_i$ switch at some $j \leq i$, the only flags we have not acknowledged of $\Pi_i$ are those that switch exactly at rank $i$. In $P_{i-1}$, the only flags that we have not acknowledged are those that contain an element of rank $i$ (including both new, mixed and old chains). The proof of Claim 1 is a result.\\

\noindent
\textbf{Claim 2:}
$$\sum_{\substack{\sigma \in \Pi \\ \rho(\sigma) = i }} \localflag{\sigmahat} \cdot \flagtop = \sum_{\substack{C \text{ a chain of } \Pi_i \\ C \text{ switches at rank }i}} \alpha^C \, \, \, \, \, \, \, \, \, \, - \sum_{\substack{C \text{ a chain of } \Pi_{i-1} \\ C \text{ contains an element of rank }i}} \alpha^C \, . $$

Recall that for any $\sigma \in \Pi$ with $\sigma$ of rank $i$ and $\sigmahat = P_1(\phi^{-1} \lowint{\sigma})$ we have that $\localflag{\sigmahat} = \Upsilon_{\widehat{\sigma}} - \Upsilon_{P_1(\partial \sigmahat)}$. By construction, since $\lowint{\sigma}$ is a rank $i$ subposet of  $\Pi$, $\phi^{-1} \lowint{\sigma}$ is a rank $i$ subposet of $\Pi_i$ and $ \partial \sigmahat$ is a rank $i$ subposet of $\Pi_{i-1}$. Note that $\partial \sigmahat$ is well defined for all $\sigma \in \Pi$ because $\phi$ is a strong Eulerian subdivision (our boundary operator was only defined on near-Eulerian posets). Expand $\localflag{\sigmahat}$ into a difference of flags.

\begin{equation}
\localcd{\sigmahat} = \Upsilon_{\sigmahat} - \Upsilon_{P_1 (\partial \sigmahat)} = \sum_{C \text{ a chain of } \sigmahat} \alpha^C - \sum_{C \text{ a chain of } P_1 (\partial \sigmahat)} \alpha^C
\end{equation}

Again, we look carefully and see what we can cancel. We may partition the chains of $\sigmahat$ into two sets: Chains contained entirely in $\partial \sigmahat$, and chains $C = C_1 < \tau$ where $\tau$ is carried by $\sigma$. Similarly, we may partition the chains $C$ of $P_1 (\partial \sigmahat)$ into two sets: those whose maximal element $\tau$ has $\rho(\tau) < i$ (and so $C$ is contained entirely in the boundary of $\sigmahat$),  and chains whose maximal element is $\widehat{1}_{\partial \sigmahat}$, the single rank $i$ element in $P_1(\partial \sigmahat)$.

It follows that chains contained strictly in the boundary of $\sigmahat$ are eliminated in the difference and we have

\begin{align*}
\localflag{\sigmahat} &= \sum_{\substack{C = C_1 < \tau \text{ a chain of } \sigmahat \\ \phi(\tau) = \sigma}} \alpha^C \, \, \, \, \, \, \, \, \, \, - \sum_{C = \dots < \widehat{1}_{\partial \sigmahat} \text{ a chain of } P_1(\partial \sigmahat)} \alpha^C \\
&= \sum_{\substack{C = C_1 < \tau \text{ a chain of } \sigmahat \\ \phi(\tau) = \sigma}} \alpha^C \, \, \, \, \, \, \, \, \, \, -  \sum_{C \text{ a chain of } \partial \sigmahat} \alpha^C \cdot b \, .
\end{align*}

Summing over all such $\sigma \in \Pi$ of rank $i$ gives the following: 

\begin{equation}
\sum_{\substack{\sigma \in \Pi \\ \rho(\sigma) = i }} \localflag{\sigmahat} \cdot \flagtop  = \sum_{\substack{\sigma \in \Pi \\ \rho(\sigma) = i }} \sum_{\substack{C = C_1 < \tau \text{ a chain of } \sigmahat \\ \phi(\tau) = \sigma}} \alpha^C  \cdot \flagtop - \sum_{\substack{\sigma \in \Pi \\ \rho(\sigma) = i }} \sum_{C \text{ a chain of } \partial \sigmahat} \alpha^C  \cdot b \cdot \flagtop \, .
\end{equation}

Now consider new and mixed chains of $\Pi_i$ that switch at $i$. Let $C$ be such a chain. Note that we may uniquely decompose $C$ into $C = C_1 < \tau < C_2$ where $\phi(\tau) = \sigma$ for some $\sigma \in \Pi$ with $\rho (\sigma) = i$, $C_1 < \tau$ is a non-empty chain of $\sigmahat$,  and $C_2$ is a possibly empty chain of $(\sigma, \widehat{1}_\Pi)$ (That is, a nondegenerate chain of $[\sigma, \widehat{1}_\Pi]$). Since any chain that switches at rank $i$ switches at some unique $\sigma \in \Pi$ of rank $i$, it follows that

\begin{equation}
\sum_{\substack{C \text{ a chain of } \Pi_i \\ C \text{ switches at rank }i}} \alpha^C = \sum_{\substack{\sigma \in \Pi \\ \rho(\sigma) = i }} \sum_{\substack{C = C_1 < \tau \text{ a chain of } \sigmahat \\ \phi(\tau) = \sigma}} \alpha^C  \cdot \flagtop \, .
\end{equation}

Similarly, consider chains of $\Pi_{i-1}$ that contain an element of rank $i$. Let $C$ be such a chain. We may uniquely decompose it into $C = C_1 < \sigma < C_2$ for some $\sigma \in \Pi$ with $\rho(\sigma)= i$ such that $C_1$ is contained entirely in $\partial \sigmahat$ and $C_2$ is contained in $(\sigma, \widehat{1}_\Pi)$. Summing this decomposition over all such chains gives us

\begin{equation} \sum_{\substack{C \text{ a chain of } \Pi_{i-1} \\ C \text{ contains an element of rank }i}} \alpha^C = \sum_{\substack{\sigma \in \Pi \\ \rho(\sigma) = i }} \sum_{C \text{ a chain of } \partial \sigmahat} \alpha^C  \cdot b \cdot \flagtop \, . 
\end{equation}

Substitute equations (7) and (8) into (6) to find the result

$$\sum_{\substack{\sigma \in \Pi \\ \rho(\sigma) = i }} \localflag{\sigmahat} \cdot \flagtop = \sum_{\substack{C \text{ a chain of } \Pi_i \\ C \text{ switches at rank }i}} \alpha^C \, \, \, \, \, \, \, \, \, \, - \sum_{\substack{C \text{ a chain of } \Pi_{i-1} \\ C \text{ contains an element of rank }i}} \alpha^C \, , $$ 
thus proving Claim 2.\\

Therefore, 

$$ 
\sum_{\substack{\sigma \in \Pi \\ \rho(\sigma) = i}} \localflag{\sigmahat} \cdot \flagtop = \Upsilon_{\Pi_i} - \Upsilon_{\Pi_{i-1}} 
$$

and the lemma holds. \qed\\

This completes the missing claim of theorem 6.1 and it follows that our $cd$-index subdivision decomposition holds!\\

In particular, if we restrict to Gorenstein* subdivisions the following result is an immediate corollary (Stated as a theorem by Ehrenborg and Karu in \cite{borgkaru}). 

\begin{corollary}

Let $\Pi$ and  $\Pihat$ be Gorenstein* posets such with a Gorenstein* subdivision $\phi: \Pihat \rightarrow \Pi$. As a coefficient wise inequality, we then have

$$\varPhi_{\Pihat} \geq \varPhi_{\Pi} \, .$$
\end{corollary} 

\emph{Proof:} Since $\phi$ is a Gorenstein* subdivision, it is also strongly Eulerian and we may apply our subdivision theorem. 

\begin{align*}\Phi_{\widehat{\Pi}} &= \sum_{\sigma \in \Pi} \localcd{\widehat{\sigma}} \cdot \Phi_{[\sigma, \widehat{1}]}\\
&= \localcd{\widehat{0}} \cdot \Phi_{[\widehat{0}, \widehat{1}]} + \sum_{\sigma \in \Pi, \sigma \neq \widehat{0}} \localcd{\widehat{\sigma}} \cdot \Phi_{[\sigma, \widehat{1}]}\\
&= \Phi_{\Pi} + \sum_{\sigma \in \Pi, \sigma \neq \widehat{0}} \localcd{\widehat{\sigma}} \cdot \Phi_{[\sigma, \widehat{1}]} \, .
\end{align*}

Note that $\sigmahat$ is near-Gorenstein* and $[\sigma, \widehat{1}]$ is Gorenstein* for all $\sigma \in \Pi$ by definition of Gorenstein* subdivision and Gorenstein* poset respectively. Since both the $cd$-index and the local $cd$-index are non-negative for Gorenstein* and near-Gorenstein* posets, the result follows. \qed

\begin{example}
 We return to the subdivided tetrahedron of Example 5.8 (pictured in Figures 18, 19 and 20) to demonstrate an application of our subdivision theorem. Let's calculate the $cd$-index of $\widehat{\Pi}$. First consider (as seen in the previous corollary)

$$\Phi_{\widehat{\Pi}} = \Phi_{\Pi} + \sum_{\sigma \in \Pi, \sigma \neq \widehat{0}} \localcd{\widehat{\sigma}} \cdot \Phi_{[\sigma, \widehat{1}]} \, . $$  

By applying our characterization of the $cd$-index of $3$-polytopes (Lemma 3.17), we have $\Phi_{\Pi} = c^3 + 2dc + 2cd$. Also, note that for any $\sigma \in \Pi$ such that $P_1 (\sigma) \cong \sigmahat$ as face posets it easily calculated that $\localcd{\sigmahat} = 0$. A similar calculation shows that $\localcd{\widehat{1}_{\Pi}} = 0$. Topologically, this implies that the only faces $\sigma \in \Pi$ with non-zero local $cd$-index are those that are non-trivially subdivided. There are exactly three such faces: the edge with a single vertex added to its interior, and the two faces adjacent to that edge that have been divided in two. Call the edge $E$ and the faces $F$ (we need not distinguish between the two of them as long as we count both because they have identical chain structure).

\begin{figure}[h]
\centering
\begin{tikzpicture}
\draw (0,0)--(2,0)--(4,0);

\draw [fill] (0,0) circle [radius = 2pt];
\draw [fill] (2,0) circle [radius = 2pt];
\draw [fill] (4,0) circle [radius = 2pt];

\end{tikzpicture}

\caption{$\localab{E} = ba + ab$ and $\localcd{E} = d$.}

\end{figure}

\begin{figure}[h]
\centering
\begin{tikzpicture}
\draw (0,0)--(2,0)--(4,0) --(2,4)--cycle;
\draw (2,0) -- (2,4);
\draw [fill] (0,0) circle [radius = 2pt];
\draw [fill] (2,0) circle [radius = 2pt];
\draw [fill] (4,0) circle [radius = 2pt];
\draw [fill] (2,4) circle [radius = 2pt];

\end{tikzpicture}

\caption{$\localab{F} = aba + aab + bba +bab$ and $\localcd{F} = cd$.}
\end{figure}

We also note that $\Phi_{[E, \widehat{1}]}= c$ and $\Phi_{[F, \widehat{1}]} = 1$. Put it all together and we find 

\begin{align*}
\Phi_{\widehat{\Pi}} &= \Phi_{\Pi} + \sum_{\sigma \in \Pi, \sigma \neq \widehat{0}} \localcd{\widehat{\sigma}} \cdot \Phi_{[\sigma, \widehat{1}]}\\ 
&= \Phi_{\Pi} + \localcd{E} \cdot \Phi_{[E, \widehat{1}]} + 2 \localcd{F} \cdot \Phi_{[F, \widehat{1}]}\\
&= (c^3 + 2dc + 2cd) + dc + 2cd\\ 
&= c^3 + 3dc + 4cd \, .
\end{align*}

This answer is easy to double check. Our classification of the $cd$-index of 3-polytopes only depended on shellability. Since $\Pihat$ is shellable, we may also apply it to $\Pihat$ and note that it gives us the same result for $\Phi_{\Pihat}$.

\end{example}

The subdivision decomposition theorem, combined with the described non-negativity conditions on Gorenstein* posets are a powerful tool for establishing upper and lower bounds on the $cd$-index via natural combinatorial arguments. In the next chapter, we examine the traditional upper and lower bounds on the face numbers of polytopes and look at the analogous $cd$-index bounds.

\clearpage
\section{Upper Bounds, Lower Bounds, and the cd-index of Stacked Polytopes}

Enumerating chains of Eulerian posets is a generalization of the classical theory of enumerating faces of a convex polytope. Though likely familiar to the reader, we present the definitions of $f$-vector and $h$-vector for reference.

\begin{definition}

Given a $CW$ $(d-1)$-complex $\Gamma$, let $f_i = | \{ \sigma \in \Gamma \, | \, \dim \sigma = i\} |$ be the number of $i$-dimensional cells of $\Gamma$. With $f_{-1}$ = 1 enumerating the single empty cell, we write $f(\Gamma) = [f_{-1}, f_0, \dots, f_{d-1}]$ to be the $f$\emph{-vector} of $\Gamma$.
\end{definition}

In practice, we may assume $\Gamma$ is a simplicial complex for the remainder of this chapter, but present the above definition for the sake of full generality. The following definition is restricted to simplicial complexes.

\begin{definition}

Given a $(d-1)$ simplicial complex $\Delta$, we define $h(\Delta) = [h_0, h_1, \dots, h_d]$ to be the $h$\emph{-vector} of $\Delta$ by the relationship

\begin{equation}
\sum_{i=0}^d f_{i-1} (x - 1)^{d-i} = \sum_{i=0}^d h_i x^{d-i}\, .
\end{equation}
\end{definition} 

Extracting coefficients, we may explicitly write

$$ h_k = \sum_{i=0}^k (-1)^{k-i} \binom{d - i}{d-k} f_{i-1} $$
and

$$ f_{k-1} = \sum_{i = 0}^{d-k} \binom{d-i}{k} h_i \, .$$
We will often consider the $h$-\emph{polynomial} written
$$ h(\Delta, x) = h_0 + h_1 x + \dots + h_d x^d \, .$$

In particular, the $h$-vector and $f$-vector of a simplicial complex are the motivation for the flag $h$-vector $\beta_P$ and flag $f$-vector $\alpha_P$ of a poset $P$ as defined in Chapter 2 (Definitions 2.7 and 2.8). In fact, the flag $f$ and $h$ vector of $P$ determine the $f$ and $h$-vector of its order complex $\Delta (P)$.

\begin{lemma}

\cite{stanleyecvol1} Let $P$ be a finite graded poset of rank $n$ with flag $f$-vector $\alpha_P$ and flag $h$-vector $\beta_P$. We may recover the $f$-vector and $h$-vector of its order complex $\Delta(P)$ by the following equations:

\begin{align*}
f_i(\Delta(P)) = \sum_{\substack{|S| = i + 1\\ S \subseteq{[n]}}} \alpha_P (S) \\
h_i(\Delta(P)) = \sum_{\substack{|S| = i \\ S \subseteq{[n]}}} \beta_P (S) \\
\end{align*}
\end{lemma}

Note that the $h$-vector and $f$-vector of a simplicial complex enumerate the same data, since each can be transformed into the other. Though the $f$-vector is more immediately intuitive, the $h$-vector possesses many useful properties when applied to convex simplicial polytopes. Many of the results of the field are much more simple to express when using the $h$-vector rather than the $f$-vector. The most immediately obvious is the \emph{Dehn-Sommerville equations}. 

\begin{theorem}
\emph{(The Dehn-Sommerville Equations)} Let $\partial P$ be the boundary complex of a simplicial $d$-polytope with $h(\partial P) = [h_0, h_1, \dots, h_d]$. Then the following equations hold:

$$ h_i = h_{d-i} \text{ for all } 1 \leq i \leq \floor{d/2} \, .$$
\end{theorem}
A proof of these equations (along with a combinatorial interpretation of the $h$-vector) is given via a shelling argument in $\cite{ziegler}$. An equivalent formulation of these equations using $f$-vectors was given in a classical paper by Sommerville \cite{sommerville}.\\

 The $h$-vector of any simplicial polytope is symmetric. This is a concise presentation of the linear relations between the number of faces of different dimension of a polytope and a powerful characterization. 

The flag $h$-vector has analogous symmetry on Eulerian posets. That is, the $ab$-index of an Eulerian polytope is invariant under the involution that switches $a$'s and $b$'s. This relationship, described as the Dehn-Sommerville relations for Eulerian posets, was originally described by Bayer and Billera \cite{bayerds} in different, but equivalent words. It is presented by Stanley in the language of flag $f$ and $h$-vectors in \cite{stanleyecvol1}.

This motivates understanding the $cd$-index as an encoding of the Dehn-Sommerville relations for posets. Given a theorem on the $h$ or $f$-vector of a polytope or simplicial sphere, we then have the natural question of "is there an analogous result for the $cd$ or $ab$-index of an Eulerian poset?" We will examine two of the most famous bounding theorems of polytope face enumeration: McMullen's upper bound theorem, and the generalized lower bound theorem.\\

The first important problem concerns an upper bound. Given a polytope $P$ of dimension $d$ with $k$ vertices, is there a sensible (hopefully tight) upper bound for for its $f$ or $h$ vector? We find an answer by looking at the class of \emph{cyclic polytopes}.

\begin{definition}
\cite{grunbaum} Consider the function $x: \mathbb{R} \rightarrow \mathbb{R}^d$, $x(t) = (t, t^2, \dots, t^d)$. For any $k \geq d+1$, we define the cyclic polytope of $k$ vertices in dimension $d$ (denoted $C_d(k)$) as the convex hull of $\{x(t_1), x(t_2), \dots, x(t_k) \}$ where $t_1, \dots , t_k$ are arbitrary real numbers such that $t_1 < t_2 < \dots t_k$.
\end{definition}

This definition may inspire skepticism due to the possible different choices of $t_i$. However, $C_d(k)$ has a well defined $f$-vector. That is, no matter the choices of the $t_i$'s, $f(C_d(k))$ is invariant \cite{grunbaum}.

Cyclic polytopes are simplicial \cite{ziegler}. They possess many notable combinatorial features. Under consideration of the topic of this thesis, we will ignore most of them (\cite{grunbaum} and \cite{ziegler} give a much more thorough explanation). Of particular note is that they have a large amount of faces relative to their number of vertices. We would present a diagram here if not for that their complexity makes all but the most trivial examples of them exceedingly difficult to display intelligently projected on a page.\\

Cyclic polytopes were conjectured to provide an upper bound for $f$-vectors by Motzkin in 1957, but a proof for the general case was not established until 1970 by McMullen \cite{ziegler}.
\begin{theorem}
\emph{(McMullen's Upper Bound Theorem)} \cite{mcmullenupperbound} Let $P$ be a $d$-polytope with $k$ vertices. Then for all integers $0 \leq n \leq d$,

\begin{align*}
f_{n-1} (P) &\leq f_{n-1} (C_d(k))\\
h_n(P) &\leq h_n(C_d(k))
\end{align*}

\end{theorem}

In other words, $P$ has as at most as many $n$-faces as the corresponding cyclic polytope of the same dimension and vertex number.\\

We omit McMullen's proof for the sake of brevity. However, it is worth noting that it makes use of both shelling and an $h$-vector argument. \\

Perhaps unsurprisingly, the result generalizes to the $cd$-index, though the correspondence is not trivial. It was not proven until 2000 by Billera and Ehrenborg \cite{billeraborgmontonicity}. Though this result predates Ehrenborg and Karu's decomposition theorem, a key lemma required for the upper bound proof is both generalized (as well as proved more succinctly) by use of a Gorenstein* subdivision decomposition \cite{borgkaru}. We present a sketch of their arguments.\\

First we consider the following lemma presented in \cite{borgreaddy} and attributed to Richard Stanley.

\begin{lemma}

 Let $P_1$ and $P_2$ be $d$-polytopes such that $P_1 \cap P_2 = F$, where $F$ is a single facet of both $P_1$ and $P_2$. Let $P_1 \cup^* P_2 = \partial P_1 \cup \partial P_2 \cup \{\widehat{1}\} \setminus \{F\}$ be the face complex formed by `joining' $P_1$ and $P_2$ at $F$. That is, $P_1 \cup^* P_2$ contains all faces of $P_1$ and $P_2$ except $F$. We may the express the $cd$-index of $P_1 \cup^* P_2$ as 

$$\Phi(P_1 \cup^* P_2) = \Phi(P_1) + \Phi(P_2) - \Phi(F) \cdot c\, .$$
\end{lemma}

\emph{Proof:} This follows from a simple chain counting argument. Any nondegenerate chain of $P_1 \cup^* P_2$ is contained in either $P_1, P_2$ or both of them. The only chains that are contained in both are those that lie in $\partial F$. The only chains of $P_1$ and $P_2$ that are not present in $P_1 \cup^* P_2$ are those whose maximal element is $F$. So we may count the chains of $P_1 \cup^*P_2$ by summing the chains of $P_1$ and $P_2$, then removing one set of the double counted chains, and both sets of the chains whose maximal element is $F$ (since there is a set contributed by both $P_1$ and $P_2$).

Expressed as a flag polynomial, we count chains contained in $\partial F$ by $\Upsilon(F) \cdot a$ and chains with maximal element $F$ by $\Upsilon(F) \cdot b$. Since we remove two sets of the latter such chains, in total we remove $\Upsilon(F) \cdot (a + 2b)$. This gives us

\begin{align*}
\Upsilon(P_1 \cup^*P_2) &= \Upsilon(P_1) + \Upsilon(P_2) - \Upsilon(F) \cdot (a + 2b) \\
\Psi(P_1 \cup^*P_2) &= \Psi(P_1) + \Upsilon(P_2) - \Psi(F) \cdot (a + b) \\
\Phi(P_1 \cup^*P_2) &= \Phi(P_1) + \Phi(P_2) - \Phi(F) \cdot c\, .
\end{align*}
\qed

Joining polytopes together at a face is often a useful construction. The preceding lemma allows us to easily calculate the $cd$-index of such objects.\\

\begin{figure}
\centering
\begin{minipage}{.33\textwidth}
\centering
\begin{tikzpicture}
\draw (0,0)--(2,0)--(2,2);
\draw[red] (2,2) -- (0,2);
\draw (0,2) -- (0,0);
\draw (1,2) node[anchor=south]{$F$};

\draw [fill] (0,0) circle [radius = 2pt];
\draw [fill] (2,0) circle [radius = 2pt];
\draw [fill] (0,2) circle [radius = 2pt];
\draw [fill] (2,2) circle [radius = 2pt];
\end{tikzpicture}
\caption{$P_1$}
\end{minipage}%
\begin{minipage}{.33\textwidth}
\centering
\begin{tikzpicture}
\draw [red] (0,0)--(2,0);
\draw (2,0) -- (1, 2) -- (0,0);

\draw (1,0) node[anchor=south]{$F$};

\draw [fill] (0,0) circle [radius = 2pt];
\draw [fill] (2,0) circle [radius = 2pt];
\draw [fill] (1,2) circle [radius = 2pt];
\end{tikzpicture}
\caption{$P_2$}
\end{minipage}%
\begin{minipage}{.33\textwidth}
\centering
\begin{tikzpicture}
\draw (0,0)--(2,0)--(2,2) -- (1, 4) --(0,2)--cycle;
\draw [dashed, red] (0,2) -- (2,2);

\draw [fill] (0,0) circle [radius = 2pt];
\draw [fill] (2,0) circle [radius = 2pt];
\draw [fill] (0,2) circle [radius = 2pt];
\draw [fill] (2,2) circle [radius = 2pt];

\draw [fill] (1,4) circle [radius = 2pt];
\end{tikzpicture}
\caption{$P_1 \cup^* P_2$}
\end{minipage}
\end{figure}

Before we continue to the aforementioned `key lemma', we must first define a pyramid operation on posets.  

\begin{definition}

Given a poset $P$, we define the \emph{pyramid} of $P$ by the set 

$$\Pyr(P):= P \times \{0,1\}$$
equipped with the partial order that for $(x,\alpha), (y, \beta) \in \Pyr(P)$ we have 

$$(x, \alpha) \leq_{\Pyr(P)} (y, \beta ) \iff x \leq_P y \text{ AND } \alpha \leq \beta \, . $$ 

Furthermore, if $P$ is a graded poset of rank $n$, $\Pyr(P)$ is a graded poset of rank $n+1$ with the rank function

$$\rho_{\Pyr(P)} (x, \alpha) := \rho_P(x) + \alpha \, .$$
\end{definition}

This operation is the natural generalization of taking a pyramid over a polytope. In the case of $P$ being the face poset of a polytope, we can interpret the element $(\widehat{0}, 1)$ as the new vertex induced by the pyramid, and any element $(x, 1)$ as the new face induced by the convex hull of $x$ and the new vertex. Any element of the form $(x, 0)$ corresponds to the old face $x \in P$

\begin{proposition}
 The pyramid operator is non-negative with respect to the $cd$-index. In particular, for any Eulerian poset $P$, 

$$\Phi_P \cdot c \leq \Phi_{\Pyr(P)} \, .$$

We refrain from proving this, but note that it is not too technical. It is best demonstrated by showing that $\Pyr$ may be well defined as a derivation on $cd$-words. See \emph{Coproducts and the cd-index} \cite{borgreaddy} and \emph{Montonicity of the cd-index for polytopes} \cite{billeraborgmontonicity} for details.
\end{proposition}

Finally, we present the key lemma. In \cite{billeraborgmontonicity}, an analogous result specifically on polytopes and their pyramids was used to prove the upper bound theorem for $cd$-indexes. The full result for Gorenstein* lattices was proved by Ehrenborg and Karu in $\cite{borgkaru}$. Recall that a \emph{lattice} is a poset such that any two elements $u$ and $v$ have both a greatest lower bound (denoted $u \wedge v$) and a least upper bound (denoted $u \vee v$). The face poset of any polytope is a Gorenstein* lattice.

\begin{lemma}
\cite{borgkaru} Let $L$ be a Gorenstein* lattice and let $x \in L$ be an element such that $x \neq \widehat{0}$ or $\widehat{1}$. We then have the following coefficientwise inequalities:

\begin{align}
\Phi_L &\geq \Phi_{[\widehat{0}, x]} \cdot \Phi_{\Pyr([x, \widehat{1}])} \\
\Phi_L &\geq \Phi_{\Pyr([\widehat{0}, x])} \cdot \Phi_{[x, \widehat{1}]} \, . 
\end{align}
\end{lemma}

\emph{Proof Sketch:} This may be proved via a direct application of the subdivision decomposition theorem (Theorem 6.1). Recall that the poset join operator $*$ corresponds to cd-index multiplication (Lemma 3.7), so $\Phi([\widehat{0}, x] * \Pyr([x, \widehat{1}]) ) = \Phi([\widehat{0}, x]) \cdot \Phi(\Pyr([x, \widehat{1}]))$. 

To complete the proof, we need to demonstrate two facts. 

\begin{itemize}
\item Show that $[\widehat{0}, x] * \Pyr([x, \widehat{1}]) $ is a Gorenstein* poset.

\item Construct a Gorenstein* subdivision $\phi: L \rightarrow [\widehat{0}, x] * \Pyr([x, \widehat{1}]) $. 
\end{itemize} 

It then follows by the subdivision decomposition theorem (in particular, from Corollary 6.5 we have that a subdivision of a Gorenstein* poset increases the coefficients of the $cd$-index by a non-negative quantity) that $$\Phi(L) \geq \Phi([\widehat{0}, x] * \Pyr([x, \widehat{1}]) ) = \Phi([\widehat{0}, x]) \cdot \Phi(\Pyr([x, \widehat{1}])) \, .$$

We leave the verification of the Gorenstein* condition to Ehrenborg and Karu in \cite{borgkaru} and note that the subdivision may be achieved by the map

$$
\phi(\sigma) =
\begin{cases}
\sigma \text{ if } \sigma < x,\\
(\sigma, 1) \text{ if } \sigma \geq x,\\
(\sigma \vee x, 0) \, \, \text{ otherwise}
\end{cases}
$$
giving us inequality (10). We find inequality (11) by considering the dual of $L$. 
\qed

\begin{corollary}
 For any Gorenstein* lattice $L$ and coatom $x$, we have 

$$\Phi_L \geq \Phi_{\Pyr([\widehat{0}, x])} \, .$$
\end{corollary}

\emph{Proof:} Apply inequality (11) and note that $\Phi_{[x, \widehat{1}]} = 1$ for $x$ a coatom. \qed\\

Though maybe not immediately impressive, pyramids are a powerful poset building block. In particular, given the boolean algebra on $n$ elements $B_n$, we have $\Pyr(B_n) \cong B_{n +1}$ as posets. Ehrenborg and Karu\cite{borgkaru} used these inequalities to easily settle a long standing conjecture of Stanley \cite{stanleycd}.

\begin{theorem}

Let $L$ be a rank $n$ Gorenstein* lattice. Then 
$$\Phi_L \geq \Phi_{B_n} \, .$$

\end{theorem}

That is, the $cd$-index is minimized by boolean algebras over all Gorenstein* lattices. \\

\emph{Proof:} Repeatedly apply our previous corollary to get $\Phi(L) \geq \Phi(\Pyr^{n-1} ([\widehat{0}, x]))$ for some atom $x$ of $P$ and note that $[\widehat{0}, x] \cong B_2$.  \qed\\

The lemma is also a component of Ehrenborg and Billera's proof of the upper bound theorem for the $cd$-index \cite{billeraborgmontonicity}. 

\begin{theorem}

 Let $P$ be a $d$-polytope with $k$ vertices. Then

$$\Phi_P \leq \Phi_{C_d(k)} \, .$$

\end{theorem}

\emph{Proof Sketch:} One can show that by `pulling' vertices of $P$ one at a time and taking the resulting convex hull, we may deform $P$ into a simplicial polytope $Q$ with the same number of vertices. An argument involving Lemma 7.7 and Lemma 7.10 demonstrates that $\Phi_P \leq \Phi_Q$. Stanley notes in $\cite{stanleycd}$ that for a simplicial polytope $Q$, the $cd$-index may be entirely determined by its $h$-vector as a sum of non-negative $cd$-polynomials determined by the dimension of $Q$. Since the $h$-vector is maximized for $C_d(k)$ by McMullen's Upper Bound Theorem, the result follows. \qed\\

So the $h$-vector and the $cd$-index of a polytope align on the upper bound theorem.\\

In general, lower bounds are harder to establish. We first describe a particular polytope important for bounds.

\begin{definition}

 Given a simplicial polytope $P$, we say that a subdivision $\Gamma$ of $P$ is a \emph{triangulation} if the following holds:

\begin{itemize}
\item Every face of $\Gamma$ is a simplex.
\item $\partial \Gamma$ and $\partial P$ are combinatorially isomorphic as simplicial complexes. 
\end{itemize}
\end{definition}
That is, $\Gamma$ subdivides $P$ into simplices such that no new faces are introduced in to the boundary of $P$. 

\begin{definition}
We say that a simplicial $d$-polytope $P$ is a \emph{stacked polytope} if there exists a triangulation $\Gamma$ of $P$ such that for any face $\sigma \in \Gamma$ with $\dim{\sigma} \leq d - 2$, we also have that $\sigma$ is a face of $P$. We say that $P$ is $(r-1)$-stacked if for any $\sigma \in \Gamma$ with $\dim{\sigma} \leq d-r$ we have that $\sigma$ is a face of $P$. When given an $(r-1)$-stacked polytope $P$, we will implicitly refer to $\Gamma$ as the \emph{stack subdivision} of $P$. 

\end{definition} 

Note that a 0-stacked polytope must be a simplex. 1-stacked polytopes are stacked polytopes. Also note that by definition, any $r$-stacked polytope of dimension $d$ is $k$-stacked for any $r \leq k \leq d$. For clarity (unless otherwise mentioned) when we refer to an $r$-stacked polytope $P$ we will assume that $r$ is the smallest such integer for which $P$ is stacked.  

A commonly used (and perhaps more intuitive) equivalent definition is that a $d$-polytope is a stacked polytope if it can be obtained from a $d$-simplex by repeatedly placing a point beyond a facet and taking the convex hull. That is, we can imagine any stacked polytope (and its underlying triangulation) arising by starting with a simplex, and then `stacking' additional simplexes to it sequentially so that they each intersect with exactly one facet of the previous structure. The only interior (that is, faces not in its boundary) faces of its stack subdivision $\Gamma$ in such a construction are of dimension $(d-1)$ or higher.\\

Stacked polytopes are relevant to the Generalized Lower Bound Theorem. This theorem was original presented as a conjecture by McMullen and Walkup in 1971 \cite{mcmullenglbc}.

\begin{theorem}

\emph{(The Generalized Lower Bound Theorem)} Given a simplicial $d$-polytope $P$, the following hold:

\begin{enumerate}
\item $1 = h_0(P) \leq h_1(P) \leq \dots \leq h_{\floor*{\frac{d}{2}}}$

\item For any integer $1 \leq r \leq \floor*{\frac{d}{2}}$, $h_{r-1}(P) = h_r(P)$ if and only if $P$ is $(r-1)$-stacked.
\end{enumerate}
\end{theorem}

 Part 1 of the theorem is implied by the well known $g$-theorem for simplicial polytopes that completely characterizes the $h$-vectors of simplicial polytopes. This particular part of the $g$-theorem was established by the work of Stanley in \cite{stanleyface}.

Part 2 was not fully resolved until much more recently. The implication that an $(r-1)$-stacked polytope $P$ has $h_{r-1}(P) = h_r(P)$ was resolved by McMullen and Walkup when they posed their conjecture \cite{mcmullenglbc}. However, the converse was not proven until 2012 by Murai and Nevo \cite{murainevo}. Extending the theorem to hold for simplicial spheres remains an active area of research. \\

This result immediately implies a lower bound for $h$-vectors analagous to the lower bound given in Theorem 7.12.

\begin{corollary}

The $h$-vector of the $d$-simplex is a lower bound for all $h$-vectors of simplicial $d$-polytopes.
\end{corollary}

\emph{Proof:} \cite{stanleyecvol1} It can be easily calculated that the $h$-vector of the $d$-simplex is $[1, 1, \dots, 1]$ . The result then follows from the lower bound theorem. \qed \\

Stacked polytopes are useful for inducing bounds because of their simple construction. Let $\Delta_d$ denote the face poset of the $d$-simplex. We may simply express the $cd$-index of stacked polytopes in terms of the $cd$-index of simplices. 

\begin{lemma}Let $P$ be a stacked $d$-polytope such that its stack subdivision contains exactly $k$ $d$-simplices. We then have that

$$\Phi(P) = k \cdot \Phi(\Delta_d) - (k - 1) \cdot \Phi(\Delta_{d-1}) \cdot c \, .$$
\end{lemma}

\emph{Proof:} By definition, we form a stacked polytope by sequentially `gluing' simplices together so that each one intersects the previous complex in exactly one facet (that is, a $(d-1)$-simplex). Note that this is exactly the setup of Lemma 7.7. Apply the lemma $(k-1)$ times and the result follows. \qed\\

In fact, this construction implies that the stack subdivision $\Gamma$ of a stacked polytope $P$ is shellable. Each $d$-simplex in the sequence intersects the previous in exactly one facet, whose boundary is a simplex and therefore shellable. It is conjectured that stack subdivision of any $(r-1)$-stacked polytope is shellable for $r \leq (d + 1)/2$ \cite{murainevo}. The `niceness' of the resulting $cd$-index formula can be seen as a consequence of how each simplex intersects the previous ones in exactly one facet.\\

What each simplex in the shelling order intersects the previous in more than one facet? That is, what can we say about the $cd$-index of a simplicial sphere that can be triangulated by a shellable complex in general? We spend the remainder of this chapter attempting to find an answer.

\begin{definition} We say a $(d-1)$-dimensional simplicial sphere is \emph{k-shellably triangulizable} if it may be triangulated as a shellable $d$-ball using exactly $k$ $d$-simplices.
\end{definition}

Note that a stacked $d$-polytope has exactly $kd + 1$ facets if its stack subdivision contains exactly $k$ $d$-simplices. So a stacked $d$-polytope with $kd + 1$ facets is  $k$-shellably triangulizable. \\

The following facts about shelling a simplex will be useful to us. The first is retrieved from \cite{ziegler}.

\begin{proposition} Consider the boundary complex of a $d$-simplex $\partial \Delta_d$. Any order $F_1 \dots F_{d+1}$ of its facets is a shelling. Furthermore, any two shellings of $\partial \Delta_d$ are isomorphic.  
\end{proposition}

\begin{proposition}

Let $G_1, G_2, \dots G_{d+1}$ be a shelling of $\partial \Delta_d$. Then for any $1 \leq r \leq d + 1$ we have that $G_1 \cup \dots \cup G_r$ is isomorphic to $\stard{d-r}$ as a simplicial complex, where $\Delta_{d-r}$ is any $(d-r)$ face of $\partial \Delta_d$. 
\end{proposition}

\emph{Proof:} Any $r$ facets of a $d$-simplex intersect in a single $(d-r)$-simplex \cite{grunbaum}. By definition, $\stard{d-r}$ is the subcomplex induced by all the facets of $\partial \Delta_d$ that contain $\Delta_{d-r}$ and their faces, which is exactly $G_1 \cup \dots \cup G_r$. \qed  

\begin{figure}
\centering
\begin{minipage}{.33\textwidth}
\centering
\begin{tikzpicture}
\draw (0,0) -- (2, 0) -- (1, 2) --cycle;
\end{tikzpicture}
\caption{star$_{\partial \Delta_3} (\Delta_2)$ }
\end{minipage}%
\begin{minipage}{.33\textwidth}
\centering
\begin{tikzpicture}
\draw (0,0) -- (2, 0) -- (1, 2) --cycle;
\draw (1,2) -- (3,1.5) -- (2,0);
\end{tikzpicture}
\caption{star$_{\partial \Delta_3} (\Delta_1)$ }
\end{minipage}%
\begin{minipage}{.33\textwidth}
\centering
\begin{tikzpicture}
\draw (0,0) -- (1, .8) -- (2,0) -- cycle;
\draw (0,0) -- (1, 2)--(1, .8);
\draw (1,.8) -- (1, 2)--(2, 0);
\end{tikzpicture}
\caption{star$_{\partial \Delta_3} (\Delta_0)$ }
\end{minipage}
\end{figure}

\begin{corollary}
In any shelling $F_1, F_2, \dots F_k$ of a simplicial $d$-dimensional complex, $F_j \cap \left( \bigcup_{i = 1}^{j-1} F_i \right) = \text{star}_{\Delta_d}(\Delta_{d-r})$ for some integer $1 \leq r \leq d+1$.
\end{corollary} 

\emph{Proof:} This follows from the definition of shelling. Since $F_j$ is a $d$-simplex,  $F_j \cap \left( \bigcup_{i = 1}^{j-1} F_i \right)$ is the initial segment of a shelling of $\partial \Delta_d$. In our previous proposition we proved this must be some star of a face of $\partial \Delta_d$. \qed\\

Note that for any $r$, $P_1(\stard{r})$ is a near-Gorenstein complex of dimension $n-1$. \\

Consider a shelling $F_1, F_2, \dots F_k$ of a an arbitrary Gorenstein complex. Note that the complex induced by the first $j$ facets of the shelling for $j < k-1$ is near-Gorenstein. We now demonstrate that adding the next facet of the shelling to the complex gives a non-negative increase in the local $cd$-index.

\begin{lemma}Let $S$ be a Gorensten complex with a shelling $F_1, F_2, \dots, F_k$. Let $\sigma_i = P_1(F_1 \cup F_2 \cup \dots \cup F_i)$ be the subcomplex induced by shelling the first $i$ facets. Then, for any integers $1 \leq i \leq k-2$, we have that 

$$\localcd{\sigma_i} \leq \localcd{\sigma_{i+1}} \, .$$

\end{lemma}

In other words, the local $cd$-index increases by a non-negative amount at each stage of a shelling.\\

\emph{Proof:} Recall that $\localflag{\sigma_i} = \Upsilon_{\sigma_i} - \Upsilon_{P_1(\partial \sigma_i)} $. Let $\Gamma := F_{i+1} \cap \sigma_i$. Let $\overline{\Gamma} : = F_{i + 1} \setminus \text{int} (\Gamma)$ be the complement of $\Gamma$ relative to $F_{i+1}$. That is, $\Gamma$ and $\overline{\Gamma}$ are near-Gorenstein complexes such that $\Gamma \cup \overline{\Gamma} = F_{i + 1}$ and $\partial \Gamma = \partial \overline{\Gamma}$.  We will count the new chains introduced in $\localflag{\sigma_{i+1}}$ relative to $\localflag{\sigma_{i}}$ by considering attaching $F_{i+1}$ to $\sigma_i$ along $\Gamma$.

\begin{align}
\localflag{\sigma_{i+1}} - \localflag{\sigma_i} &=
(\Upsilon_{\sigma_{i+1}} - \Upsilon_{P_1(\partial \sigma_{i+1})}) - (\Upsilon_{\sigma_{i}} - \Upsilon_{P_1( \partial \sigma_{i})})\\
&= (\Upsilon_{\sigma_{i+1}} - \Upsilon_{\sigma_i}) - (\Upsilon_{P_1( \partial \sigma_{i+1})} - \Upsilon_{P_1( \partial \sigma_{i})})\\
&= (\Upsilon_{\Gamma} \cdot b + \Upsilon_{\overline{\Gamma}} \cdot (a + b) - \Upsilon_{\partial \Gamma} \cdot a(a + b) ) - (\Upsilon_{\overline{\Gamma}} - \Upsilon_{\Gamma})(a + b)\\
&= \Upsilon_{\Gamma} (a + 2b) - \Upsilon_{\partial \Gamma} \cdot a (a+ b) \, .
\end{align}

In the above equations, the only non-trivial step is from (13) to (14) by means of a chain counting argument. See Figure 29 for reference. The remainder of the calculations follow from algebraic manipulation. We now outline the counting argument. \\

\begin{figure}
\centering
\begin{tikzpicture} 

\draw [fill, gray!30!white] (0,0).. controls (0,1) and (0,2).. (1,2).. controls (2,2) and (2,1).. (2,0) --cycle;

\draw [fill, gray!30!white] (0,0).. controls (0,-1) and (0,-2).. (1,-2).. controls (2,-2) and (2,-1).. (2,0) --cycle;

\draw (0,0).. controls (0,1) and (0,2).. (1,2).. controls (2,2) and (2,1).. (2,0) --cycle;

\draw (0,0).. controls (0,-1) and (0,-2).. (1,-2).. controls (2,-2) and (2,-1).. (2,0) --cycle;
 
\draw (1,0) node[anchor= north]{$\Gamma$};
\draw (1,1) node{$F_{i+1}$};
\draw (1,2) node[anchor=south]{$\overline{\Gamma}$};
\draw (1,-1) node{$\sigma_i$};

\draw [fill] (0,0) circle [radius = 2pt];
\draw [fill] (2,0) circle [radius = 2pt];

\end{tikzpicture}
\caption{A sketch of $\sigma_{i+1}$. Note that the bottom cell represents all of the complex $\sigma_i$, while the top cell represents the single face $F_{i+1}$. $\Gamma = F_{i+1} \cap \sigma_i$ and $\overline{\Gamma} = F_{i + 1} \setminus \text{int} (\Gamma)$.}
\end{figure}
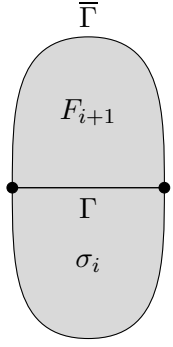

\emph{Claim 1:} $\Upsilon_{\sigma_{i+1}} - \Upsilon_{\sigma_i} =  \Upsilon_{\Gamma} \cdot b + \Upsilon_{\overline{\Gamma}} \cdot (a + b) - \Upsilon_{\partial \Gamma} \cdot a(a + b)$.\\

Since $\sigma_i \subseteq \sigma{i+1}$ as posets, way may enumerate $\Upsilon_{\sigma_{i+1}} - \Upsilon_{\sigma_i}$ by counting exactly the chains of $\sigma_i+1$ not present in $\sigma_i$. Note that any chain that does not contain a face of $F_{i+1}$ must be present in both $\sigma_{i+1}$ and $\sigma_i$ and is eliminated in the difference. Any remaining chain must have an initial segment in $\Gamma$, $\overline{\Gamma}$, or both (if its initial segment is in $\partial \Gamma$ and contains no max dimensional faces of $\Gamma$ or $\overline{\Gamma}$). By basic inclusion-exclusion, we count the $ab$ words corresponding to chains with initial segments in $\Gamma$ and $\overline{\Gamma}$, and then subtract the chains with initial segments in $\partial \Gamma$. \\
 
The only chains of $\sigma_{i+1}$ with an initial segment in $\Gamma$ that are \emph{not} in $\sigma_{i}$ as well are exactly those with $F_{i+1}$ as a maximal element. These may be counted by $\Upsilon_{\Gamma} \cdot b$. The only chains of $\sigma_{i+1}$ with an initial segment in $\overline{\Gamma}$ that are not in $\sigma_{i}$ are of exactly two types. If the initial segment is in the interior of $\overline{\Gamma}$, then the chain is not in $\sigma_i$ by definition. If the initial segment is in $\partial \Gamma$, it must have $F_{i+1}$ as its maximal element. So we enumerate these chains by $\Upsilon_{\overline{\Gamma}}\cdot(a + b) - \Upsilon_{\partial \Gamma} \cdot aa$. Finally, we may enumerate all chains of $\sigma_{i+1}$ with initial segment in $\partial \Gamma$ that are not present in $\sigma_i$ by $\Upsilon_{\partial \Gamma} \cdot ab$ (because they must have $F_{i+1}$ as a maximal element). So by inclusion-exclusion it follows that

\begin{align*}
\Upsilon_{\sigma_{i+1}} - \Upsilon_{\sigma_i} &= \Upsilon_{\Gamma} \cdot b + (\Upsilon_{\overline{\Gamma}} \cdot (a + b) - \Upsilon_{\partial \Gamma} \cdot aa) - \Upsilon_{\partial \Gamma} \cdot ab \\
&= \Upsilon_{\Gamma} \cdot b + \Upsilon_{\overline{\Gamma}} \cdot (a + b) - \Upsilon_{\partial \Gamma} \cdot a(a + b) \, . 
\end{align*}

\emph{Claim 2:} $\Upsilon_{P_1( \partial \sigma_{i+1})} - \Upsilon_{P_1( \partial \sigma_{i})} = (\Upsilon_{\partial \sigma_{i+1}} - \Upsilon_{\partial \sigma_i})(a+b)$.\\

 Any chain of $\partial \sigma_{i+1}$ that is not contained in $\partial \sigma_{i}$ must lie in the interior $\overline {\Gamma}$. The only chains of $\partial \sigma_{i}$ not contained in $\partial \sigma_{i+1}$ are those on the interior of $\Gamma$. Either of type of chain may be covered by the single maximal element or not. Since the chains contained strictly in the boundary are eliminated in the difference, we may write

$$\Upsilon_{P_1( \partial \sigma_{i+1})} - \Upsilon_{P_1( \partial \sigma_{i})} = (\Upsilon_{\overline{\Gamma}} - \Upsilon_{\Gamma} )(a + b) \, .$$ 

Combining these two claims gives that equations (13) and (14) of our argument are equal.\\

Now, we convert the result of equation (15) into the local $cd$-index

\begin{align*}
\localflag{\sigma_{i+1}} - \localflag{\sigma_i} &=\Upsilon_{\Gamma}(a + 2b) - \Upsilon_{\partial \Gamma} \cdot ab - \Upsilon_{\partial \Gamma} \cdot aa \\
&= \Upsilon_{\Gamma}(a + 2b) - \Upsilon_{\partial \Gamma} \cdot a(b+a)\\
\localab{\sigma_{i+1}} - \localab{\sigma_i}&= \Psi_\Gamma (a + b) - \Psi_{\partial \Gamma} \cdot (a - b)a\\
&= (\localab{\Gamma} + \Psi_{\partial \Gamma} \cdot a)(a + b) - \Psi_{\partial \Gamma} \cdot (a - b)a\\
&= \localab{\Gamma} (a + b) + \Psi_{\partial \Gamma}(ab + ba)\\
\localcd{\sigma_{i+1}} - \localcd{\sigma_i} &= \localcd{\Gamma} \cdot c + \Phi_{\partial \Gamma} \cdot d
\end{align*}

Since $\Gamma$ is near-Gorenstein, we have that $\localcd{\Gamma} \cdot c$ is non-negative. Similarly, since $\partial \Gamma$ is Gorenstein we have that $\Phi_{\partial \Gamma}\cdot d$ is non-negative. It follows that $\localcd(\sigma_i) \leq \localcd(\sigma_{i+1})$. \qed \\

This is a slightly stronger result than we set out to prove. The above equations provide an explicit description of how the local $cd$-index changes at each stage of a shelling.

\begin{theorem}
Let $\Pi$ be the face poset of a $k$-shellably triangulizable $(d-1)$-simplicial sphere. Then the $cd$-index of a stacked $d$-polytope with $kd + 1$ facets is an upper bound for the $cd$-index of $\Pi$. In other words, 

$$\Phi_{\Pi} \leq \Phi_{S_d(k)}$$
where $S_d(k)$ denotes a $k$-triangulizable stacked $d$-polytope.

\end{theorem}

\emph{Proof:} Let $T$ be the underlying triangulation of $\Pi$. Suppose $T$ has shelling order $F_1, F_2 \dots F_k$. For $1 \leq i \leq k$, let $T_i = F_1 \cup \dots \cup F_i$ be the subcomplex of $T$ induced by the first $i$ facets in the shelling order. How can we relate the $cd$-index of $\partial T_i$ and $\partial T_{i+1}$? Note that $\partial T_k = S$. 

We get from $T_i$ to $T_{i+1}$ by attaching exactly one $k$-simplex, $F_{i+1}$. If $F_{i+1}$ intersects $T_i$ in exactly one facet, then we are simply in the case of lemma 7.7 and have that  

$$\Phi_{\partial T_{i+1}} = \Phi_{\partial T_i} + \Phi_{\Delta_d} - \Phi_{\Delta_{d-1}} \cdot c \, . $$

Otherwise, $F_{i+1}$ intersects $T_i$ in some initial segment of a shelling of $\partial F_{i+1}$. Since $F_{i+1}$ is a simplex, we know from our previous propositions that $F_{i+1} \cap T_i \cong \text{star}_{\Delta_d}(\Delta_{d-r})$ for some $r$. Note that $\text{star}_{\Delta_d}(\Delta_{d-r})$ is a near-Gorenstein complex. Also note that the remaining facets of $F_{i+1}$ not in the intersection form a simplicial complex isomorphic to $\text{star}_{\Delta_d}(\Delta_{r-1})$, and that $\partial\stard{d-r} = \partial \stard{r-1}$.\\

\begin{figure}[h]
\centering
\begin{tikzpicture}

\draw [thick,red]  (0,0).. controls (1,0) and (2,0).. (2,1).. controls (2,2) and (1,2).. (0,2);

\draw (0,0).. controls (-1,0) and (-2,0).. (-2,1).. controls (-2,2) and (-1,2).. (0,2);

\draw[thick,blue, dashed] (0,0) -- (0,2);
 
\draw (-1,1) node{$\partial T_i$};
\draw (1,1) node{$F_{i+1}$};

\draw [fill=green] (0,0) circle [radius = 2pt];
\draw [fill=green] (0,2) circle [radius = 2pt];

\end{tikzpicture}
\caption{A generalized sketch of $\partial T_{i+1}$. The \textcolor{red}{red} line corresponds to $\stard{r-1}$, the \textcolor{green}{green} dots to $\partial \stard{r-1}$, and the \textcolor{blue}{blue} dashes to the faces in the interior of $\stard{d-r}$ that are in $\partial T_i$ but not $\partial T_{i+1}$.}
\end{figure}
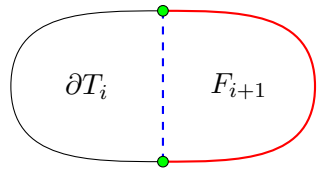

We now calculate the $ab$-index of $\partial T_{i+1}$ relative to $\partial T_i$. $F_{i+1}$ intersects $T_i$ at $\stard{d-r}$. We partition the chains of $\partial T_{i+1}$ into exactly two types: chains that are contained in both $\partial T_{i+1}$ and $\partial T_i$, and chains that are contained only in $\partial T_{i+1}$.\\

 All chains of $\partial T_i$ are contained in $\partial T_{i+1}$ with the exception of those whose maximal face lies in the interior of $\stard{d-r}$. Subtracting out exactly those chains gives the following quantity.

\begin{align}
&\Upsilon_{\partial T_i} - (\Upsilon_{\stard{d-r}} - \Upsilon_{ \partial \stard{d-r}} \cdot a) \nonumber \\
 &= \Upsilon_{\partial T_i} - (\Upsilon_{\stard{d-r}} - \Upsilon_{ \partial \stard{d-r}} \cdot (a+b) + \Upsilon_{ \partial \stard{d-r}} \cdot b) \nonumber \\
  &= \Upsilon_{\partial T_i} - (\Upsilon_{\stard{d-r}} - \Upsilon_{P_1(\partial \stard{d-r})} + \Upsilon_{ \partial \stard{d-r}} \cdot b) \nonumber \\
  &= \Upsilon_{\partial T_i} - (\localflag{\stard{d-r}} + \Upsilon_{ \partial \stard{d-r}} \cdot b ) \, .
\end{align}

The chains of $T_{i+1}$ that are not contained in $T_i$ are exactly those are contained in the interior of. $\text{star}_{\Delta_d}(\Delta_{r-1})$ which may be counted by

\begin{equation}
\Upsilon_{\stard{r-1}} - \Upsilon_{ \partial \stard{r-1}} \cdot a \,.
\end{equation}

Now, remembering that $\partial \stard{r-1} \cong \partial \stard{d-r}$, we sum (16) and (17) to compute $\Phi_{\partial T_{i+1}}$.

\begin{align*}
\Upsilon_{\partial T_{i+1}} &= \Upsilon_{\partial T_i} - (\localflag{\stard{d-r}} + \Upsilon_{ \partial \stard{d-r}} \cdot b ) 
+ \Upsilon_{\stard{r-1}} - \Upsilon_{ \partial \stard{r-1}} \cdot a \\
\Psi_{\partial T_{i+1}} &= \Psi_{\partial T_i} - (\localab{\stard{d-r}} + \Psi_{\partial \stard{d-r}} \cdot b ) 
+ \Psi_{\stard{r-1}} - \Psi_{ \partial \stard{r-1}} \cdot (a-b) \\
&= \Psi_{\partial T_i} - \localab{\stard{d-r}} + \Psi_{\stard{r-1}} - \Psi_{ \partial \stard{r-1}} \cdot a\\
&= \Psi_{\partial T_i} + \localab{\stard{r-1}} - \localab{\stard{d-r}}\\
\Phi_{T_{i+1}} &= \Phi_{\partial T_i} + \localcd{\stard{r-1}} - \localcd{\stard{d-r}}\, .
\end{align*}

In particular, consider the final terms $\localcd{\stard{r-1}} - \localcd{\stard{d-r}}$. If $\Pi$ is a stacked polytope, then $r=1$ since $F_{i+1}$ will intersect exactly one facet of $T_i$. Given a stacked $\Pi$, this holds for all $1 \leq i \leq k-1$.\\

\emph{Claim:} $\localcd{\stard{r-1}} - \localcd{\stard{d-r}}$ is maximized when $r=1$.\\

Note that if this claim holds, the theorem follows. For all $i < k$, the $cd$-index of $\partial T_{i+1}$ will increase from $\partial T_i$ by at most the amount it would if we were stacking a simplex onto $\partial T_i$.\\

By applying Lemma 7.23 to a shelling of a $d$-simplex, we have that for all  $1 \leq i \leq d-1$

$$\localcd{\stard{i}} \leq  \localcd{\stard{i-1}} \, . $$

That is,

$$\localcd{\stard{d-1}} \leq \localcd{\stard{d-2}} \leq \dots \leq \localcd{\stard{1}} \leq \localcd{\stard{0}} \, . $$

Since all of the above stars are near-Gorenstein complexes, their local $cd$-indexes are non-negative and it follows that for any integer $1 \leq r \leq d$ 

$$\localcd{\stard{r-1}} - \localcd{\stard{d-r}} \leq \localcd{\stard{r-1}} - \localcd{\stard{d-r}}\,. $$

It follows that the claim holds and our theorem is an immediate consequence. \qed\\

We may easily find a rough upper bound on the $cd$-index of a simplical $d$-polytope with $k$-facets by the following construction.

\begin{corollary}

Let $P$ be a simplicial polytope with $k$ facets. Then the $cd$ index of $P$ has the $cd$-index of $S_d(k)$ as an upper bound. That is,

$$\Phi_P \leq \Phi_{S_d(k)} .$$
\end{corollary} 

\emph{Proof:} Note that we may triangulate $P$ into some near-Gorenstein complex $\Gamma$ by placing a point $x$ on its interior, and then defining the faces of $\Gamma$ to be $$\Gamma := \mathcal{F}(P) \cup \{\text{conv}(F \cup x) \, | \, F \in \mathcal{F}(P)\}$$ where $\mathcal{F}(P)$ is the face set of $P$ and `conv' denotes the convex hull. The facets of $\Gamma$ are exactly the simplices conv$(F \cup x)$ for $F$ a facet of $P$. Since $P$ is a polytope, it has a shelling $F_1, F_2, \dots F_k$. It is simple to see that this induces a shelling $\text{conv}(F_1 \cup x), \text{conv}(F_2 \cup x) \dots \text{conv}(F_k \cup x)$ of $\Gamma$. So $P$ is $k$-shellably triangulizable and we can apply our previous theorem. \qed \\

This is the most trivial possible bound that arises from our theorem and is likely a much worse upper bound than a cyclic polytope. However, if one was studying a class of specific polytopes that triangulates nicely (i.e. into a number of simplices significantly smaller than the number of its facets) into a shellable complex, Theorem 7.24 may yet generate some useful bounds. 

\clearpage
\section{Morphisms between the $cd$-index, the toric $h$-vector and the local $h$-vector}

In the previous chapter, we discussed the classical $h$-vector for simplicial complexes. We now investigate the \emph{toric} and \emph{local} $h$-vector. These were introduced by Richard Stanley in \cite{stanleygeneralh} and \cite{stanleylocal} and have been the subject of active research ever since. The toric $h$-vector is an extension of the $h$-vector to arbitrary Eulerian posets that maintains most of its `nice' properties. The toric $h$-vector is a tool that bridges between  polyhedral combinatorics and algebraic geometry. The local $h$-vector is a is used in decomposing the $h$-vector relative to a subdivision. This decomposition behaves analogously to our $cd$-index subdivision decomposition. In \cite{bayerborg}, Bayer and Ehrenborg define an algebra morphism that maps the $cd$-index of a poset $P$ to the toric $h$-vector of $P$. Our ultimate goal is to demonstrate that this morphism gives a natural map between the $cd$-index subdivision decomposition and the $h$-vector subdivision decomposition, and use this to connect the local $cd$-index to the local $h$-vector.\\

We begin by defining the local $h$-vector with respect to a simplex. We eventually will generalize to Eulerian posets and the toric $h$-vector. Though doing so immediately would perhaps lead to less exposition, the theory of the toric $h$-vector and local $h$-vector for arbitrary Eulerian posets becomes exceedingly technical very quickly. We encourage the reader to study the simplical case to develop intuition before moving on to the general case. 

To align ourself with the notation of $\cite{stanleylocal}$, recall that given a vertex set $V$, we may denote the simplex with vertices labeled by $V$ as $2^V$. Furthermore we may identify faces uniquely by their vertex set, giving us a face for any $W \subseteq V$. This definition is consistent with that of an abstract simplicial complex. 

\begin{definition}

\cite{stanleylocal} Let $V$ be some vertex set with $|V| = d$. Given a simplicial subdivision $\Gamma$ of $2^V$, we define the polynomial $\ell_V (\Gamma, x) := \ell_0 + \ell_1 x + \dots + \ell_d x^d$ by the following equation: 

\begin{equation}
h(\Gamma, x) = \sum_{W \subseteq V} \ell_W (\Gamma_W, x) \, .
\end{equation}
\end{definition}

We call $\ell_V (\Gamma, x)$ the \emph{local h-polynomial} of $\Gamma$. Note that $\ell_V (\Gamma, x)$ is implicitly dependent on the base simplex $2^V$ of the subdivision. We analogously call $\ell_V (\Gamma) := [\ell_0, \ell_1, \dots, \ell_d]$ the \emph{local h-vector} of $\Gamma$. 

It may not be immediately obvious that this actually defines the local $h$-polynomial. However, this is actually its most concise expression. Inverting the previous formula using the method of inclusion-exclusion, we explicitly find that
\begin{equation}
\ell_V (\Gamma, x) = \sum_{W \subseteq V} (-1)^{|V \setminus W|} h(\Gamma_W, x) \, .
\end{equation}
We present some simple low-dimension examples from Stanley. \cite{stanleylocal}. These may be calculated relatively easily by recursive techniques. Things get more complicated in dimensions 3 and larger.
\begin{example}

Let $\Gamma$ be a subdivision of $2^V$.
\begin{itemize} 
\item $\ell_{\emptyset} (\emptyset, x) = 1$. 
\item For $|V| > 0$, consider the trivial subdivision $\Gamma = 2^V$. Then $\ell_V (2^V, x) = 0$.
\item For $|V| = 2$, $2^V$ is a line segment. Suppose $\Gamma$ adds $t$ interior points to the line. Then $\ell_V(\Gamma, x) = tx$. 
\item For $|V| = 3$, calculate $h(\Gamma, x) = h_0 + h_1x + h_2x^2 + h_3 x^3$. Then $\ell_V (\Gamma, x) = h_2x + h_2x^2$.
\end{itemize}
\end{example}
We use the above for reference to work through a local $h$-vector decomposition of the barycentric subdivision of a triangle. 

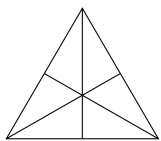
\begin{figure}[h]
\centering
\begin{tikzpicture}
\draw (0,0) -- (2,0) -- (1,1.73) --cycle;
\draw (1,0) -- (1,1.73);
\draw (0,0) -- (1.5, 0.865);
\draw (2,0) -- (0.5, 0.865);
\end{tikzpicture}
\caption{The barycentric subdivision of a triangle.}
\end{figure}

\begin{example} Consider a triangle $\Delta_2$ and let $\Gamma$ be the barycentric subdivision of it pictured in figure 30. Counting faces, we have $f(\Gamma) = [1, 7, 12, 6]$ (the empty face, seven vertices, 12 edges and 6 2-faces). We may compute that 
$$h(\Gamma, x) = 1 + 4x + x^2 \, .$$
\end{example}

Any vertex $V$ is not subdivided and has a local $h$-vector of 0. Note that any edge $E$ of $\Delta_2$ has exactly one vertex introduced into its interior by $\Gamma$. So $\ell_{E} (\Gamma_E, x) = x$. And from the final bullet of our previous example, for the single maximal face $F$ of $\Delta_2$ we have $\ell_F(\Gamma_F, x) = x + x^2$. Plugging these values into equation (18) gives

\begin{align*}
h(\Gamma, x) &= \ell_{\emptyset} (\Gamma_{\emptyset}, x) + 3 \ell_V(\Gamma_V, x) + 3\ell_E(\Gamma_E, x) + \ell_F(\Gamma_V, F)\\
&= 1 + 3 \cdot 0 + 3 \cdot x + (x + x^2) \\
&= 1 + 4x + x^2
\end{align*}

and verifies our earlier computation of $h(\Gamma, x)$.\\

Note that every local $h$-vector computed so far has been symmetrical. That is for any local $h$-vector $\ell(\Gamma) = [\ell_0, \ell_1, \dots, \ell_d]$ we have $\ell_i = \ell_{d-i}$ for all $0 \leq i \leq d$. This is not a coincidence. However, we refrain from stating this as a theorem until we reach the general case. \\

Restricting ourselves only to subdivisions of a simplex is unsatisfying. It turns out that the local $h$-vector may be applied to subdivisions simplicial complexes in general. Stanley presents the following theorem in \cite{stanleylocal}.

\begin{theorem}
Let $\Delta$ be a pure $d$-dimensional complex with a simplicial subdivision $\Gamma.$ Then the following equation holds. 

\begin{equation}
h(\Gamma, x) = \sum_{F \in \Delta} \ell_F(\Gamma_F, x)h(\text{link}_{\Delta} F, x) \, .
\end{equation}
\end{theorem}

Note that if $\Delta$ is a simplex, $h(\text{link}_{\Delta} F, x) = 1$ for all $F \in \Delta$ and the above theorem reduces to equation (16), our definition of the local $h$-vector. Recall that the face poset of $\text{link}_{\Delta} F$ is exactly the interval $[F, \widehat{1})_{\Delta}$ (that is, all the faces of $\Delta$ that contain $F$). If we replace every occurrence of the $h$-polynomial with the $cd$-index, this formula appears to be very similar to the $cd$-index subdivision decomposition. There exists a general version of this theorem (that we will see shortly) that makes the resemblance even more apparent. \\

We now introduce the toric and local $h$-vector in full generality. We mostly align ourselves with the conventions and notation of \cite{katzstapledon}.

\begin{definition} Let $P$ be an Eulerian poset of rank $n$. We define the \emph{g-polynomial} $g(P,x)$ of $P$ recursively as follows:

\begin{enumerate}
\item If $n=0$, $g(P,x) = 1$.
\item If $n > 0$, $g(P,x) \in \mathbb{Z}[x]$ is the unique polynomial with degree $ < n/2 $ that satisfies

\begin{equation}
x^n g(P, 1/x) = \sum_{\sigma \in P} g(\lowint{\sigma}, x) (x-1)^{n-\rho(\sigma)} \, .
\end{equation}

\end{enumerate}
\end{definition}

We define the \emph{g-vector} of $P$ by $g(P) := [g_0, g_1, \dots g_m]$ where $m = \floor{(n-1)/2}$.

\begin{definition}
 Let $P$ be a lower Eulerian poset of rank $n$. We then define the \emph{h-polynomial} $h(P,x)$ of $P$ to be the unique polynomial $h(P,x) \in \mathbb{Z}[x]$ satisfying

\begin{equation}
x^n h(P, 1/x) = \sum_{\sigma \in P} g(\lowint{\sigma}, x) (x-1)^{n-\rho(\sigma)} \, .
\end{equation}
\end{definition}

We define the \emph{h-vector} of $P$ by $h(P) := [h_0, h_1, \dots, h_n]$.\\

Note that if $P$ is Eulerian, we have $g(P,x) = h(P,x)$. Though notationally dense, it is easy to see that this $h$-polynomial reduces to the one previously defined for simplicial complexes (as in Definition 7.2) after applying the following proposition.

\begin{proposition}\cite{stanleyecvol1} For any rank $n$ boolean algebra $B_n$, we have $g(B_n,x) = 1$ .
\end{proposition}

If $P$ is the face poset of a simplicial complex, then $\lowint{\sigma} \cong B_k$ for any $\sigma \in P$ and for some $k \leq n$. So all $g$ terms disappear and we have 

\begin{align*}
x^n h(P, 1/x) &= \sum_{\sigma \in P} (x-1)^{n-\rho(\sigma)} \\
&= \sum_{i=0}^n f_{i-1} (x-1)^{n-i}
\end{align*}

which is exactly the definition of the $h$-polynomial presented in equation (10).

\begin{definition} Let $P$ be a Eulerian poset of rank $n > 0$. Then we call $h(P \setminus \{\widehat{1}\}, x)$ the \emph{toric h-vector} of $P$. 
\end{definition}

\begin{example} \cite{stanleyecvol1} Consider the boolean algebra on $d+1$ elements, $B_{d+1}$ where $d \geq 0$. Then

$$
h(B_{d+1} \setminus \{\widehat{1}\}) = 1 + x + \dots + x^d \, .$$
\end{example}

Recall that the $h$-vector of the boundary complex of a polytope holds a number of particularly significant combinatorial properties (we spent much of Chapter 7 discussing them!). The toric $h$-vector generalizes many of these properties to Eulerian posets. In particular, for any Eulerian $P$, $h(P \setminus \{\widehat{1}\}, x)$ is symmetric and we may generalize the Dehn-Sommerville relations beyond simplicial polytopes. 

\begin{theorem} \emph{(The Generalized Dehn-Sommerville Equations)} Let $P$ be an Eulerian poset of rank $d+1$ with $h(P \setminus\{\widehat{1}\}, x) = h_0 + h_1x + \dots h_d x^d$. Then for all $0 \leq i \leq d$, 
$$h_i = h_{d-i}\, .$$
\end{theorem}

\emph{Proof:} This follows from a Mobi\"us inversion argument. See Stanley's textbook \cite{stanleyecvol1}.\\

We also mentioned that the $h$-vector of the boundary complex of a simplicial polytope is both positive and unimodal in the lower bound theorem. Does this result hold for the toric $h$-vector? In the case of the face posets of rational polytopes, Stanley demonstrated the affirmative \cite{stanleygeneralh}. 

\begin{theorem} \emph{(Stanley)}Let $P$ be the face poset of a rational $d$-polytope (that is, a polytope such that all its vertices have rational coordinates). Consider its toric h-vector $h(\partial P) = [h_0, h_1, \dots, h_d]$. Then, for $m = \floor{d/2}$ we have

$$1 \leq h_0 \leq h_1 \leq \dots \leq h_m \, .$$
\end{theorem}

 This result comes from interpreting the $h_i$ coefficients as the dimensions of the even cohomology groups of a projective toric variety associated with the normal fan of $P$ \cite{katzstapledon}. Much of the modern interest in the theory of polytopes and fans is motivated by their connections to toric varieties and algebraic geometry. \\
 
 We now define yet another type of subdivision, originally presented by Katz and Stapledon \cite{katzstapledon}.
 
\begin{definition}
Let $\phi: \Gamma \rightarrow \Pi$ be an order preserving and rank increasing function between lower Eulerian posets $\Gamma$ and $\Pi$. We say that $\phi$ is \emph{strongly surjective} if
 
 \begin{itemize}
 \item $\phi$ is surjective.
 \item For all $z \in \Gamma$ and $x \in \Pi$ such that $\phi(z) \leq x$, there exists some $y \in \Gamma$ such that $z \leq y$, $\phi(y) \leq x$ and $\rho_{\Gamma}(y) = \rho_{\Pi}(x)$.
 \end{itemize}
 
 \end{definition}
 In other words, if $z$ is carried by $x$ in $\Pi$, there must be some $y$ in $\Gamma$ of equal rank to $x$ that is ordered above $z$.
 
\begin{definition}
 Let $\phi: \Gamma \rightarrow \Pi$ be an order preserving, rank increasing and strongly surjective function between lower Eulerian posets $\Gamma$ and $\Phi$. We say $\phi$ is a \emph{strong formal subdivision} if and only if the following condition holds for all $z \in \Gamma$ and $x \in \Pi$ such that $\phi(z) \leq x$:

\begin{equation}
\sum_{\substack{\phi(y) \leq x \\ z \leq y}} (-1)^{\rho_{\Pi}(x) - \rho_{\Gamma} (y)} = 
\begin{cases}
1 \text{ if } \phi(z) = x\\
0 \text{ otherwise.}
\end{cases}
\end{equation}
\end{definition}

This generalizes the idea of a formal subdivision as defined by Stanley in \cite{stanleylocal}.

\begin{lemma} Any strongly Eulerian subdivision is a strong formal subdivision.
\end{lemma}

\emph{Proof:} Recall the definition of a strongly Eulerian subdivision. That is, given a surjective order preserving function $\phi: \Gamma \rightarrow \Pi$ between rank $n$ lower Eulerian posets, for any $x \in \Pi$ we have the following: 

\begin{itemize}
\item $P_1(\phi^{-1}[\widehat{0}, x])$ is a near-Eulerian poset.
\item $\phi^{-1}[\widehat{0}, x]$ has rank equal to $\rho_{\Pi} (x)$.
\end{itemize}

The second condition immediately implies strong surjectivity. Now let us verify that equation (23) holds. \\

Let $x \in \Pi$ be an element of rank $n$, and let $z,w \in P_1(\phi^{-1} \lowint{x})$. Because $P_1(\phi^{-1} \lowint{x})$ is near-Eulerian, the following equation holds.
$$\sum_{z \leq y \leq w} (-1)^{\rho(w) - \rho(y)} =
\begin{cases}
1 \text{ if } w = \widehat{1} \text{ and } \phi(z) \neq x \\
0 \text{ otherwise.}
\end{cases}
 $$
 
 Take $w = \widehat{1}$ and we may calculate
 
\begin{align*}
\sum_{z \leq y \leq \widehat{1}} (-1)^{n+1 - \rho(y)} &=
\begin{cases}
1 \text{ if } \phi(z) \neq x \\
0 \text{ otherwise.}
\end{cases}\\
\sum_{z \leq y < \widehat{1}} (-1)^{n+1 - \rho(y)} &=
\begin{cases}
0 \text{ if } \phi(z) \neq x \\
-1 \text{ otherwise.}
\end{cases}
\end{align*}

However, elements $y \in P_1(\phi^{-1} \lowint{x})$ such that $z \leq y < \widehat{1}$ are exactly the elements of $\Gamma$ such that $z \leq y$ and $\phi(y) \leq x$. Noting that $\rho_{\Pi}(x) = n$ and dividing by $-1$ gives us the result that 

$$
\sum_{\substack{\phi(y) \leq x \\ z \leq y}} (-1)^{\rho_{\Pi}(x) - \rho_{\Gamma} (y)} = 
\begin{cases}
0 \text{ if } \phi(z) \neq x\\
1 \text{ otherwise.}
\end{cases}
$$
which is equivalent to equation (24) and completes the lemma. \qed\\

It immediately follows that Gorenstein* subdivisions are strong formal subdivisions. Additionally, any $CW$-subdivision (or polyhedral subdivision) is a strong formal subdivision \cite{stanleygeneralh}.\\

We now present the local $h$-vector in its full generality, as defined by \cite{katzstapledon}.

\begin{definition}

Let $\phi: \Gamma \rightarrow \Pi$ be a strong formal subdivision between lower Eulerian $\Gamma$ and Eulerian $\Pi$. We recursively define the \emph{local h-polynomial} $\ell_{\Pi} (\Gamma, x) \in \mathbb{Z}[x]$ to be the polynomial satisfying

\begin{equation}
h(\Gamma, x) = \sum_{\sigma \in \Pi} \ell_{\lowint{\sigma}}(\Gamma_{\sigma}, x) g(\upint{\sigma}, x) \, .
\end{equation}
\end{definition}

For $\ell_{\Pi} (\Gamma, x) = \ell_0 + \ell_1 x + \dots + \ell_n$, we write $\ell_{\Pi} (\Gamma) = [\ell_0 \dots \ell_n]$ and call $\ell_{\Pi} (\Gamma)$ the local $h$-vector of $\Gamma$ (relative to $\Pi$).\\

 If $\Pi$ is a simplex, all $g$ terms are $1$ and this reduces to our original definition of the local $h$-polynomial in equation (19). Much like that definition, this feels ultimately unsatisfying because we would like to describe $\ell_{\Pi} (\Gamma, x)$ explicitly. However, this presentation of the local $h$-polynomial is much cleaner than its explicit definition, which we now present for posterity \cite{katzstapledon}.

\begin{proposition} Given $\phi: \Gamma \rightarrow \Pi$, a strong formal subdivision between lower Eulerian $\Gamma$ and Eulerian $\Pi$, we may define $\ell_{\Pi} (\Gamma, x)$ by 

$$ \ell_{\Pi} (\Gamma, x)  = \sum_{\sigma \in \Pi} h(\Gamma_{\sigma}, x) (-1)^{\rho(\sigma, \widehat{1})} g([\sigma, \widehat{1}]^*, x)
$$ 

where $[\sigma, \widehat{1}]^*$ is the dual poset of $[\sigma, \widehat{1}]$ (that is, the poset on the same underlying set with all order relations reversed).
\end{proposition}

Demonstrating why these two formulations of the local $h$-vector are equivalent is highly technical. We have been ignoring much of the underlying theory whenever possible in an effort to be concise. At its most general, the $h$-vector may be studied from the perspective of kernels and acceptable functions on the incidence algebra of a poset. Kernels and acceptability operators are a broad theory generalizing the study of elements of an algebra that are invariant with respect to some involution. For a detailed treatment, see Katz and Stapledon \cite{katzstapledon} and/or Stanley \cite{stanleylocal}.\\

Recall that we earlier discussed the symmetry of the local $h$-vector. We now formalize that statement \cite{katzstapledon}. 

\begin{theorem} Let $\phi: \Gamma \rightarrow \Pi$ be a strong formal subdivision, where $\Gamma$ is a rank $n$ lower Eulerian poset and $\Pi$ is an Eulerian poset. Then 

$$\ell_{\Pi} (\Gamma, x) = x^n \ell_{\Pi} (\Gamma, 1/x) \, . $$ 
\end{theorem}

So we have a symmetry in the local $h$-vector analogous to the Dehn-Sommerville equations. Can we say anything about non-negativity or unimodality of the local $h$-vector? In some cases, yes. Given a polyhedral complex $\Pi$, we say that $\Gamma$ is a \emph{regular subdivision} if it arises from projecting the `lower faces' of some convex polytope $P$ on to $\Pi$ (we refrain from a formal definition). Stanley demonstrated the following \cite{stanleylocal}.

\begin{theorem} \emph{(Stanley)} Let $P$ be a rational $d$-polytope. Suppose $\Gamma$ is some rational polytopal subdivision of P (that, is all vertices of the subdivision may be realized with rational coordinates). Let $\ell_{\Pi} (\Gamma, x) = \ell_0 + \ell_1 x + \dots + \ell_d$ Then the following holds:

\begin{enumerate}
\item $\ell_i \geq 0$ for all $i$. 
\item If additionally $\Gamma$ is a regular subdivision, we have $0 \leq \ell_0 \leq \ell_1 \leq \dots \leq \ell_m$ where $m = \floor{d/2}$.
\end{enumerate}
\end{theorem}

Like his proof of the unimodality of the toric $h$-vector for rational polytopes, this is not a very combinatorial proof and instead goes through the intersection cohomology of toric varieties. \\

We now aim to demonstrate a correspondence between the local $cd$-index and the local $h$-vector via an algebra morphism developed by Bayer and Ehrenborg in \cite{bayerborg}. The following coproduct structure is described in depth in \cite{borgreaddy}. Note that we differ from their notation slightly (namely, we hesitate to denote a coproduct as $\Delta$ because we have spent most of this document using that to refer to simplices). 

\begin{definition}
Let $\mathbb{Z}\langle a, b \rangle$ denote the algebra of non-commutative indeterminates $a$ and $b$ over $\mathbb{Z}$. 
\end{definition}

\begin{itemize}
\item Define $M: \mathbb{Z}\langle a, b \rangle \otimes \mathbb{Z}\langle a, b \rangle \rightarrow \mathbb{Z}\langle a, b \rangle$ to multiply the tensors together. That is, for $u,v \in \mathbb{Z}\langle a, b \rangle$ we have 

$$M(u \otimes v) = uv \, .$$

\item We define a coproduct $C: \mathbb{Z}\langle a, b \rangle \rightarrow\ \mathbb{Z}\langle a, b \rangle \otimes \mathbb{Z}\langle a, b \rangle$ by the following. For any monomial $w = w_1 w_2 \dots w_n$ we define

$$C(w_1 \dots w_n) = \sum_{i=1}^n w_1 \dots w_{i-1} \otimes w_{i+1} \dots w_n $$

and extend to arbitrary elements of $\mathbb{Z}\langle a, b \rangle$ by linearity.

\end{itemize} 
 
For our coproduct, we may use Sweedler notation to improve readability at the expense of utility and instead write for any element $w \in \mathbb{Z}\langle a, b \rangle$,

$$ C(w) = \sum_{w} w_{(1)} \otimes w_{(2)} \, .$$

For any poset $P$, we have the $ab$-index $\Psi_P \in \mathbb{Z}\langle a, b \rangle$.\\

We first need to establish some notation and define a few more maps before we can express Bayer and Ehrenborg's morphism \cite{bayerborg}. 

\begin{itemize}
\item Let $U_{\leq m} [\cdot] : \mathbb{Z}[x] \rightarrow \mathbb{Z}[x]$ denote the map that truncates all terms of degree larger than $m$. That is, for $f(x) = \sum_{i=0}^{k} a_i x^i$ we have

$$U_{\leq m}[f(x)] = 
\sum_{i=0}^{\min{\{m,k\}}} a_i x^i 
$$

\item We define the algebra map $\kappa: \mathbb{Z}\langle a, b \rangle \rightarrow \mathbb{Z}[x]$ by $\kappa(a) = x - 1$, $\kappa (b) = 0$ and extend linearly for all $ab$-words and elements of $\mathbb{Z}\langle a, b \rangle$. Note that this implies for any poset $P$ of rank $n$, $\kappa (\Psi_P) = (x-1)^{n-1}$.
\end{itemize}

We may now present a morphism between the $ab$-index and the toric $h$-vector \cite{bayerborg}.

\begin{theorem} \emph{(Bayer and Ehrenborg)} Let $f$ and $g$ be linear maps $f,g : \mathbb{Z}\langle a, b \rangle \rightarrow \mathbb{Z}[x]$ defined by the following intertwined recursive definition:

\begin{itemize}
\item For any monomial $v$, 
\begin{align*}
f(v) &= \kappa(v) + M \circ (g \otimes \kappa) \circ C(v) \\
&= \kappa(v) + \sum_v g(v_{(1)}) \kappa(v_{(2)})
\end{align*}

and extend linearly.

\item For any monomial $v$ of degree $n$ with $m = \floor{n/2}$ set 

$$g(v) = U_{\leq m} [(1-x)f(v)] $$

and extend linearly.
\end{itemize}

Then, for any graded poset $P$ with $\widehat{0}$ and $\widehat{1}$, $f(\Psi_P) = h(P \setminus \{\widehat{1} \}, x)$ and $g(\Psi_P) = g(P, x)$ where $h(P \setminus \{\widehat{1} \}, x)$ is the toric $h$-polynomial of $P$ and $g(P,x)$ is the $g$-polynomial of $P$.
\end{theorem}

\emph {Proof:} This may be proved by a simple inductive argument. See \cite{bayerborg} for details.\\

We claim that these morphisms are well behaved with respect to our $ab$-index subdivision decomposition.

\begin{lemma} Let $\phi: \Pihat \rightarrow \Pi$ be a strong Eulerian subdivision between Eulerian posets. For any  $\sigma \in \Pi$ we have

$$f(\localab{\sigmahat}) = \ell_{\sigma} (\sigmahat, x) \, .$$
\end{lemma}

\emph{Proof:} Recall that $\localab{\sigmahat} = \Psi_{\sigmahat} - \Psi_{P_1( \partial \sigma )}$. So by applying our previous theorem 8.20, we have 

\begin{align*}
f(\localab{\sigmahat}) &= f(\Psi_{\sigmahat}) - f(\Psi_{P_1( \partial \sigmahat )})\\
&= h(\sigmahat \setminus \{\widehat{1}\}, x) - h(\partial \sigmahat, x) \, .
\end{align*}

Since strong formal subdivisions restricted to order ideals are still strong formal subdivision \cite{katzstapledon}, it follows that by restricting $\phi$ to $\sigma$ we have that $\sigmahat \setminus \{\widehat{1}\}$ and $\partial \sigmahat$ are both strong formal subdivisions of $\sigma$. So we may apply a local $h$-vector decomposition relative to the face poset of $\sigma$, giving us

\begin{equation}
f(\localab{\sigmahat}) = \sum_{F \in \sigma} \ell_F(\sigmahat_F, x) \cdot h([F, \widehat{1}]_{\sigma},x) - \sum_{F \in \sigma} \ell_F((\partial \sigmahat)_F, x) \cdot h([F, \widehat{1}]_{\sigma},x) \, .
\end{equation}

Note that for any $F \in \sigma$ such that $F \neq \widehat{1}$, we have $F \in \partial \sigma$ and it follows that $\sigmahat_F \cong \partial \sigmahat_F$. So for $F \neq \widehat{1}$ we have 

$$\ell_F(\sigmahat_F, x)  = \ell_F((\partial \sigmahat)_F, x) $$

So all terms with $F \neq \widehat 1$ cancel in equation (26). Also, $\widehat{1}_{\sigma}$ is only trivially subdivided with respect to $\partial \sigmahat$ so it follows that $\ell_{\widehat{1}} ((\partial \sigmahat)_{\widehat{1}}, x) = 0$.

The only term left standing is 

$$
f(\localab{\sigmahat})  = \ell_{\widehat{1}}(\sigmahat_{\widehat {1}}, x) \cdot h([\widehat{1}, \widehat{1}]_{\sigma},x) =  \ell_{\sigma} (\sigmahat, x) \cdot 1 
$$

and the result follows. \qed.\\

Since $\Pi$ is Eulerian, it also follows that for any $\sigma \in \Pi$, $f(\Phi_{\upint{\sigma}}) = h([\sigma, \widehat{1}), x)$. We may summarize the result into the following corollary.

\begin{corollary} Let $\phi: \Pihat \rightarrow \Pi$ be a strong Eulerian subdivision, with $f: \mathbb{Z} \langle a,b \rangle \rightarrow \mathbb{Z}[x]$ defined as above. Then the following diagram holds:
\begin{figure}[h!]
\centering
\begin{tikzpicture}
\draw (0,4) node{$\Psi_{\Pihat}$};
\draw (2,4) node{$=$};
\draw (4,4) node{$\sum_{\sigma \in \Pi}$};
\draw (6,4) node{$\localab{\sigmahat}$};
\draw (7,4) node{$\cdot$};
\draw (8,4) node{$\Psi_{\upint{\sigma}}$};

\draw (0,0) node{$h(\Pihat \setminus \{\widehat{1} \}, x)$};
\draw (2,0) node{$=$};
\draw (4,0) node{$\sum_{\sigma \in \Pi \setminus \{\widehat{1} \} }$};
\draw (6,0) node{$\ell_{\sigma}(\sigmahat, x)$};
\draw (7,0) node{$\cdot$};
\draw (8,0) node{$h([\sigma, \widehat{1}), x)$};

\draw [->] (0,3.5) -- (0,.5);
\draw [->] (6,3.5) -- (6,.5);
\draw [->] (8,3.5) -- (8,.5);

\draw (0,2) node[anchor=east]{$f$};
\draw (6,2) node[anchor=east]{$f$};
\draw (8,2) node[anchor=east]{$f$};
\end{tikzpicture}

\end{figure}
\end{corollary}

Consider restricting $\mathbb{Z} \langle a,b \rangle$ to the subalgebra $\mathbb{Z} \langle c,d \rangle$ where $c = a + b$ and $d = ab + ba$. We then have $\kappa(c) = a + b$ and $\kappa(d) = 0$ and $f$ and $g$ remain well defined. So the above morphisms are also valid for $\Phi_{\Pihat}$, $\localcd{\sigmahat}$ and $\Phi_{\upint{\sigma}}$ \cite{bayerborg}. The morphism $f$ gives a direct correspondence between the $cd$-index subdivision decomposition and the local $h$-vector subdivision decomposition. 

\clearpage
\section{Closing Remarks}
We finish with a few open questions. In general, the $cd$-index contains strictly more combinatorial data than the $h$-polynomial and we can not recover $\Phi_P$ from $h(P,x)$. However, when $P$ is simplicial its flags are entirely determined by its $h$-vector \cite{stanleylocal}. Defining some morphism from $\mathbb{Z}[x]$ to $\mathbb{Z} \langle a,b \rangle$  or  $\mathbb{Z} \langle c,d \rangle$ that takes the $h$-polynomial of a simplicial poset to the $ab$ or $cd$-index while behaving well relative to local decomposition of a subdivision would be a nice construction.

In \cite{katzstapledon}, Katz and Stapledon define the \emph{mixed} $h$-polynomial of a strong formal subdivision $\phi: \Pihat \rightarrow \Pi$. This is a two variable specialization of the $h$-polynomial that encodes the $h$-polynomial of $\Pihat$ along while tracking the \emph{excess} of each element of $\Pihat$ relative to the subdivision $\Phi$. Finding a map from the $cd$-index to the mixed $h$-polynomial, or developing a mixed analog of the $cd$-index that behaves analogously with respect to the local subdivision decomposition would be a natural continuation of the work of this thesis.

\clearpage
\bibliography{bibliography}
\bibliographystyle{plain}

\end{document}